\documentclass[a4paper,12pt]{article}
\begin{document}

\title{Line integration and second order partial differential equations
over Cayley-Dickson algebras.}
\author{Ludkovsky S.V.}
\date{10 August 2010}
\maketitle
\begin{abstract}
Line integration of generalized functions is studied. Second order
partial differential equations with piecewise continuous and
generalized variable coefficients over Cayley-Dickson algebras are
investigated. Formulas for integrations of such equations are
deduced. For this purpose a non-commutative line integration is
used. Examples of solutions of partial differential equations are
given.
\end{abstract}

\section{Introduction.}
\par It is well-known, that differential equations have many-sided applications
in different sciences including physics, mechanics, other natural
sciences, techniques, economics, etc. The differential equations
also are very important for mathematics
\cite{grubb,hoermpd,rubinstb,kneschkeb, matveevb,vladumf}.
Predominantly differential equations are considered over fields such
as real, complex, or with non-archimedean norms. Recently they are
also begun to be studied over Clifford algebras
\cite{guhaspr,guespr,guesprqa}.
\par Such algebras have a long history, because quaternions were first introduced by W.R.
Hamilton in 1843. He had planned to use them for problems of
mechanics and mathematics \cite{hamilt,rothe}. Their generalization
known as the octonion algebra was introduced by J.T. Graves and A.
Cayley in 1843-45. Then Dickson had investigated more general
algebras known now as the Cayley-Dickson algebras
\cite{baez,dickson,kansol}.
\par The Cayley-Dickson algebras, particularly, octonions and quaternions are
widely used in physics, but mainly algebraically. Already Maxwell
had utilized quaternions to derive his equations of electrodynamics,
but then he had rewritten them in real coordinates.
\par In the 50-th of the 20-th century Yang and Mills had used them in quantum field theory,
but theory of functions over octonions and quaternions in their
times was not sufficiently developed to satisfy their needs.
Discussing that situation they have formulated the problem of
developing analysis over octonions and quaternions \cite{guetze}.
This is natural, because quantum fields are frequently non-abelian
\cite{solov}.  Dirac had used complexified quaternions to solve the
Klein-Gordon hyperbolic differential equation with constant
coefficients. \par This work continues previous articles of the
author. In those articles (super)-differentiable functions of
Cayley-Dickson variables and their non-commutative line integrals
were investigated
\cite{ludfov,ludoyst,luladfcdv,lujmsnfpcd,lujmsnfpcd}. In the papers
\cite{lutsltjms,ludeoc} differential equations and their systems
over octonions and quaternions were studied. \par The Cayley-Dickson
algebras ${\cal A}_r$ have the even generator $i_0=1$ and the purely
imaginary odd generators $i_1,...,i_{2^r-1}$, $2\le r$, $i_k^2=-1$
and $i_0i_k=i_k$ and $i_ki_l= - i_li_k$ for each $1\le k\ne l$. For
$3\le r$ the multiplication of these generators is generally
non-associative, so they form not a group, but a non-commutative
quasi-group with the property of alternativity $i_k (i_ki_l) =
(i_k^2)i_l$ and $(i_li_k)i_k=i_l(i_k^2)$ instead of associativity.
Ordinary super-analysis operates with graded algebras over Abelian
groups. Therefore, super-analysis over the Cayley-Dickson algebras
is in some respect more complicated than usual super-analysis, for
example, over the Grassman algebras. \par The aim of this paper is
in developing of Dirack's approach on partial differential equations
with variable piecewise continuous or generalized coefficients.
\par The technique presented there is developed here below for solutions
of partial differential equations of the second order of arbitrary
signatures and with variable coefficients which may also be
piecewise continuous or generalized functions. Moreover, signatures
may change piecewise in a domain. Formulas for integrations of such
equations are deduced. For this purpose a non-commutative line
integration of generalized functions is developed. Examples of
partial differential equations are given. Moreover the approach of
\S \S 2-25 over the Cayley-Dickson algebras ${\cal A}_v$ gives the
fundamental solution of any first and second order linear partial
differential equation with variable $z$-differentiable ${\cal
A}_v$-valued coefficients, $z\in U\subset {\cal A}_v$, where $U$ is
a domain in ${\cal A}_v$ satisfying some mild convexity conditions
described below. These results can be used for solutions of concrete
partial differential equations or their systems of different orders
with piecewise continuous or generalized coefficients, for example,
of Helmholtz' or Klein-Gordon's types, which are important in optics
of composite materials or quantum field theory. Finally solutions of
some types of non-linear partial differential equations over
Cayley-Dickson algebras are studied.
\par Main results of this paper are obtained for the first time.

\section{Partial differential equations of the second order.}
\par {\bf 1. Remarks and notations.} For a subset $U$
in either the quaternion skew field ${\bf H} = {\cal A}_2$ or in the
octonion algebra ${\bf O}={\cal A}_3$ or the Cayley-Dickson algebra
${\cal A}_r$, $r\ge 4$, we put $\pi _{{\sf s},{\sf p},{\sf t}}(U):=
\{ {\sf u}: z\in U, z=\sum_{{\sf v}\in \bf b}w_{\sf v}{\sf v},$
${\sf u}=w_{\sf s}{\sf s}+w_{\sf p}{\sf p} \} $ for each ${\sf s}\ne
{\sf p}\in \bf b$, where ${\sf t}:=\sum_{{\sf v}\in {\bf b}\setminus
\{ {\sf s}, {\sf p} \} } w_{\sf v}{\sf v} \in {\cal A}_{r,{\sf
s},{\sf p}}:= \{ z\in {\cal A}_r:$ $z=\sum_{{\sf v}\in \bf b} w_{\sf
v}{\sf v},$ $w_{\sf s}=w_{\sf p}=0 ,$ $w_{\sf v}\in \bf R$ $\forall
{\sf v}\in {\bf b} \} $, where ${\bf b} := \{ i_0,i_1,...,i_{2^r-1}
\} $ is the family of standard generators of the algebra ${\cal
A}_r$ so that $i_j^2 = -1$, for each $j\ge 1$, $i_ji_k = - i_ki_j$
for each $j\ne k \ge 1$, $i_0=1$. Geometrically the domain $\pi
_{{\sf s},{\sf p},{\sf t}}(U)$ means the projection on the complex
plane ${\bf C}_{{\sf s},{\sf p}}$ of the intersection $U$ with the
plane ${\tilde \pi }_{{\sf s},{\sf p},{\sf t}}\ni {\sf t}$, ${\bf
C}_{{\sf s},{\sf p}} := \{ a{\sf s}+b{\sf p}:$ $a, b \in {\bf R} \}
$, since ${\sf s}{\sf p}^*\in {\hat b}:={\bf b}\setminus \{ 1 \} $.
Recall that in \S \S 2.5-7 \cite{ludfov} for each continuous
function $f: U\to {\cal A}_r$ it was defined the operator ${\hat f}$
by each variable $z\in {\cal A}_r$. If a function $f$ is
$z$-differentiable by the Cayley-Dickson variable $z\in U\subset
{\cal A}_r$, $2\le r$, then ${\hat f}(z) =dg(z)/dz$, where
$(dg(z)/dz).1 = f(z)$.
\par A Hausdorff topological space $X$ is said to be $n$-connected
for $n\ge 0$ if each continuous map $f: S^k\to X$ from the
$k$-dimensional real unit sphere into $X$ has a continuous extension
over $\bf R^{k+1}$ for each $k\le n$  (see also \cite{span}). A
$1$-connected space is also said to be simply connected.
\par  It is supposed further, that a domain $U$ in ${\cal A}_r^m$
has the property that \par $(D1)$ each projection ${\bf
p}_j(U)=:U_j$ is $(2^r-1)$-connected; \par $(D2)$ $\pi _{{\sf
s},{\sf p},{\sf t}}(U_j)$ is simply connected in $\bf C$ for each
$k=0,1,...,2^{r-1}$, ${\sf s}=i_{2k}$, ${\sf p}=i_{2k+1}$, ${\sf
t}\in {\cal A}_{r,{\sf s},{\sf p}}$ and $u\in {\bf C}_{{\sf s},{\sf
p}}$, for which there exists $z={\sf u}+{\sf t}\in U_j$, \\ where
$e_j = (0,...,0,1,0,...,0)\in {\cal A}_r^m$ is the vector with $1$
on the $j$-th place, ${\bf p}_j(z) = \mbox{ }^jz$ for each $z\in
{\cal A}_r^m$, $z=\sum_{j=1}^m\mbox{ }^jz e_j$, $\mbox{ }^jz\in
{\cal A}_r$ for each $j=1,...,m$, $m\in {\bf N} := \{ 1,2,3,... \}
$. Frequently we take $m=1$. Henceforward, we consider a domain $U$
satisfying Conditions $(D1,D2)$ if any other is not outlined.
\par The family of all ${\cal A}_r$ locally analytic functions $f(x)$
on $U$ with values in ${\cal A}_r$ is denoted by ${\cal H}(U,{\cal
A}_r)$. It is supposed that a locally analytic function $f(x)$ is
written in the $x$-representation $\nu (x)$, also denoted by $\nu
=\nu ^f$. The latter is equivalent to the super-differentiability of
$f$ (see \cite{ludoyst,ludfov, lutsltjms}). Each such $f$ is
supposed to be specified by its phrase $\nu $.
\par For each super-differentiable function $f(x)$ its
non-commutative line integral $\int_{\gamma } f(x)dx$ in $U$ is
defined along a rectifiable path $\gamma $ in $U$. It is the
integral of a differential form ${\hat f}(x).dx$, where \par $(I1)$
${\hat f}(x) =dg(x)/dx$, \par $(I2)$ $[dg(x)/dx)].1=f(x)$ for each
$x\in U$. \par A branch of the non-commutative line integral can be
specified with the help of either the left or right algorithm (see
\cite{ludoyst,ludfov, lutsltjms}). We take further for definiteness
the left algorithm if something another will not described. For
$f\in {\cal H}(U,{\cal A}_r)$ and a rectifiable path $\gamma :
[a,b]\to {\cal A}_r$ the integral $\int_{\gamma }f(x)dx$ depends
only on an initial $\alpha =\gamma (a)$ and final $\beta =\gamma
(b)$ points due to the non-commutative analog of the homotopy
theorem in $U$, where $a<b\in {\bf R}$. When initial and final
points or a path are not marked we denote the operation of the
non-commutative line integration in the domain $U$ simply by $\int
f(x)dx$ analogously to the indefinite integral.
\par To rewrite a function from real variables $z_j$
in the $z$-representation the following identities are used:
\par $(1)$ $z_j=(-zi_j+ i_j(2^r-2)^{-1} \{ -z
+\sum_{k=1}^{2^r-1}i_k(zi_k^*) \} )/2$ \\ for each
$j=1,2,...,2^r-1$, $$(2)\quad z_0=(z+ (2^r-2)^{-1} \{ -z +
\sum_{k=1}^{2^r-1}i_k(zi_k^*) \} )/2,$$  where $2\le r\in \bf N$,
$z$ is a Cayley-Dickson number decomposed as
\par $(3)$ $z=z_0i_0+...+z_{2^r-1}i_{2^r-1}\in {\cal
A}_r$, $z_j\in \bf R$ for each $j$, $i_k^* = {\tilde i}_k = - i_k$
for each $k>0$, $i_0=1$, since $i_k(i_0i_k^*)=i_0=1$, $i_k(i_ji_k^*)
= - i_k(i_k^*i_j) = - (i_ki_k^*) i_j = - i_j$ for each $k\ge 1$ and
$j\ge 1$ with $k\ne j$ (shortly $k\ne j\ge 1$), $i_k (i_ki_k^*)
=i_k$ for each $k\ge 0$.
\par As usually $C^0(U,{\cal A}_v)$ denotes the $\bf R$-linear space
of all continuous ${\cal A}_v$-valued functions $f: U\to {\cal
A}_v$. More generally $C^n(U,{\cal A}_v)$ denotes the $\bf R$-linear
space of all $n$ times continuously differentiable by real variables
$z_0,...,z_{2^v-1}$ functions $f: U\to {\cal A}_v$, where $n\in {\bf
N}$. Certainly, $C^n(U,{\cal A}_v)$ can be supplied with the
structure of left- and right-module over the Cayley-Dickson algebra
${\cal A}_v$ using point-wise multiplication of functions $f(z)$ on
Cayley-Dickson numbers from the left and the right.
\par {\bf 2. Factorization and integration of equations.}
\par We consider the second order partial differential equation:
\par $(1)$ $Af=g$, where
$$A= \sum_{l,m=1}^k {\bf a}_{l,m}\partial ^2/\partial
\tau _l\partial \tau _m + \sum_{l=1}^k \alpha _l\partial /\partial
\tau _l $$ is a partial differential operator of the second order.
Let us suppose that the quadratic form $$a(\tau ) := \sum_{l,m} {\bf
a}_{l,m} \tau _l\tau _m$$ is non-degenerate and is not always
negative, because otherwise we can consider $-A$. Moreover, let a
matrix of coefficients be real and symmetric ${\bf a}_{l,m}(\tau
)={\bf a}_{m,l}(\tau )\in {\bf R}$, $\alpha _l, \tau _l\in {\bf R}$
for each $l, m =1,...,k$. Then we reduce this form $a(\tau )$ by an
invertible $\bf R$ linear operator $C=C(\tau )$ to the sum of
squares. This means, that
$$(2)\quad A = \sum_{l=1}^k {\bf b}_l
\partial ^2/ \partial s_l^2 + \sum_{l=1}^k \beta _l \partial /\partial s_l ,$$ where
$\partial s_j/\partial \tau _l = C_{l,j}(\tau ) $, $C= (C_{l,j})$,
with real-valued functions ${\bf b}_l$ and $\beta _l$ for each $l$.
Here  $${\bf b}_l \delta _{j,l} = \sum_{p,m} {\bf a}_{p,m} C_{p,j}
C_{m,l}\mbox{  and}$$  $$\beta _j = \sum {\bf a}_{p,m} (\partial
C_{p,j} /\partial \tau _m) + \sum_{v=1}^k \alpha _v C_{v,j} $$ for
all $j, l=1,...,k$. In the case when coefficients of $A$ are
constant, using a multiplier of the type $\exp (\sum_l \epsilon _l
s_l)$ it is possible to reduce this equation to the case so that if
${\bf b}_l\ne 0$, then $\beta _l = 0$ (see \S 3, Chapter 4 in
\cite{rubinstb}).  Therefore, one can as usually simplify the
operator with the help of such change of coordinates and consider
that only $\beta _1$ may be non-zero if ${\bf b}_1=0$.
\par Thus one can choose an invertible real matrix $(c_{h,m})_{h,m
=1,...,k}$ corresponding to $C=C(\tau )$ so that ${\bf b}_l\le 0$
for $p+1\le l\le k$ and ${\bf b}_l\ge 0$ for $0<l\le p$, where
$0<p\le k$, $ ~ q:=k-p$. When $q=0$ and $\beta _l=0$ for each $l$
the operator is elliptic, for $q=0$ and $\beta _1\ne 0$ the operator
is parabolic, for $0<p<k$ and $\beta _l=0$ for each $l$ the operator
is hyperbolic. Sometimes the matrix $C$ can be chosen constant on a
domain, where the signature $(p, q)$ of the quadratic form $a(\tau
)$ is constant. We suppose that the sums $\sum_{l=1}^p {\bf
b}_l^2(x)>0$ and $\sum_{l=p+1}^k {\bf b}_l^2(x)>0$ are positive
$\lambda $-almost everywhere on a domain $U$, where $\lambda $ is
the measure induced by the Lebesgue measure on the real shadow of
the  Cayley-Dickson algebra. Generally the natural number $k-p=q(x)$
may either be constant or change while crossing the surface $ \{
x\in U: ~ \sum_{l=1}^k {\bf b}_l^2(x)=0 \} $, when the domain $U$
satisfies Conditions 1$(D1,D2)$.
\par We consider elliptic and hyperbolic partial differential operators
reduced to the sum of squares
$$(3)\quad A = [\sum_{l=0}^k {\bf b}_l(x)
\partial ^2/ \partial x_l^2] ,$$
where ${\bf b}_l(x)\in {\bf R}$ for all
$x=x_0i_0+...+x_{2^r-1}i_{2^r-1}$ in the open domain $U\subset {\cal
A}_r$ satisfying Conditions 1$(D1,D2)$ in the Cayley-Dickson algebra
${\cal A}_r$, $1\le k \le 2^r-1$, $2\le r\le 3$. Practically the
coefficient ${\bf b}_l$ can depend only on $x_0,...,x_k$ remaining
$z$-differentiable in definite $z$-representations due to Formulas
1$(1-3)$ for each $l$.
\par More generally we can consider partial differential operators of the form
\par $(4)$ $A = c_1B_1 +...+ c_mB_m$, where $c_jB_jf=c_j(B_jf)$,
while each
$$(4')\quad B_j=\sum_{k=m_1+...+m_{j-1}+1}^{m_1+...+m_j} {\bf
b}_k(x)\partial ^2/\partial x_k^2$$  is an elliptic partial
differential operator of the second order by variables
$x_{m_1+...+m_{j-1}+1}$,...,$x_{m_1+...+m_j}$; $c_j \in {\cal A}_r$
with $Re (c_j)\ge 0$ for each $1\le j\le l$, $Re (c_j)<0$ for every
$j>l$, with $|c_j|=1$ for each $j=1,...m$, where $1\le r $, $1\le
l<m$, $m_0=0$.
\par We remind, that Dirac had used complexified bi-quaternions to solve Klein-Gordon's hyperbolic
partial differential equation with constant coefficients appearing
in spin problems. That is, he had decomposed d'Alembert's operator
$\partial ^2/\partial t^2 - \nabla ^2$ as the product ${\bf i}^*
\sigma {\bf i} \sigma $ over the complexified bi-quaternion algebra
${\bf H}_{\bf C}$ with the first order differential operator $\sigma
$.
\par If follow this approach one takes the complexified
Cayley-Dickson algebra
\par $(5)$ $({\cal A}_r)_{\bf C} = {\cal A}_r \oplus {\cal A}_r{\bf
i}$, \\ where ${\bf i}$ is taken to be commuting with $i_j$ for each
$j=0,...,2^r-1$. Now the algebra $({\cal A}_r)_{\bf C}$ is already
not the division algebra even for $2\le r\le 3$, that is two non
zero elements with zero product occur in it. Then each element $z =
(z_1,0)$ in $({\cal A}_r)_{\bf C}$ can be written in the $2\times 2$
matrix form ${ {z_1 ~ 0} \choose{0 ~ z_1}}$ and $z=(0,z_2)$ can be
written in the form ${ {0 ~ ~ ~ z_2}\choose{ - z_2 ~ 0 }}$, where
entries $z_1, z_2\in {\cal A}_r$ are Cayley-Dickson numbers, ${\bf
i} = { {0 ~ ~ 1}\choose{ - 1 ~ 0 }}$.
\par Let each coefficient $c_j$ be written in the polar form \par
$(6)$ $c_j = \exp (i_{\kappa (j)}\gamma _j)$ \\ with $0\le |\gamma
_j| \le \pi $, $j=1,...,m$, $1\le r $, $1\le \kappa (j)\le \kappa
(j+1)$ for each $j$. Put $p=p_1+...+p_m$, where $p_j= 0$ for either
$\gamma _j=0$ or $\kappa (j)=\kappa (j-1)$, while $p_j=1$ for
$\gamma _j\ne 0$ and $\kappa (j)\ne \kappa (j-1)$. Up to an
isomorphism we take the Cayley-Dickson algebra ${\cal A}_v$ with
$v\ge r$ satisfying inequalities $2^{v-1} < 2^p (m+1) \le 2^v$.
Further we make the complexification $({\cal A}_v)_{\bf C}$ of the
Cayley-Dickson algebra ${\cal A}_v$.
\par Take two non-negative integer numbers $0\le r$ and $v$ with $r\le v\in {\bf
Z}$. We consider the quotient algebra over the real field ${\cal
A}_v/{\cal A}_r =: {\cal A}_{r,v}$. For $r=v$ this algebra is
isomorphic with the real field $\bf R$. For $r<v$ the algebra ${\cal
A}_{r,v}$ is isomorphic with $\bigoplus_{k=0}^{2^{v-r}-1}{\bf R}
i_{2^rk}$. The latter algebra is produced by generators $ \{
i_{2^r\beta }: ~ \beta =2^{\gamma }-1; ~ \gamma =0,1,...,v-r \} $
and their finite ordered products, that gives the generators set $
\{ i_{2^rk}: ~ k=0,...,2^{v-r}-1 \} $, where generators satisfying
the numbering rule $i_j i_{2^s} = i_{j+2^s}$ for each $1\le s$,
$j=0, 1,...,2^s-1$ can be taken up to an isomorphism of the
Cayley-Dickson algebra ${\cal A}_{s+1}$. Therefore, the algebra
${\cal A}_{r,v}$ is isomorphic with the Cayley-Dickson algebra
${\cal A}_{v-r}$, since the doubling procedure can be started from
another suitable purely imaginary Cayley-Dickson numbers such as
generators \cite{baez,kansol}. But we consider in ${\cal A}_{r,v}$
its specific generators basis $ \{ i_{2^rk}: ~ k=0,...,2^{v-r}-1 \}
$.
\par For each Cayley-Dickson
numbers $x, y \in {\cal A}_r$ we define the real-valued scalar
product
\par $(RS)$ $(x,y) = (x,y)_r := Re (x{\tilde y})$, \\ where ${\tilde z} = z^*$
denotes the conjugated number, while $Re (y) := (y+ y^*)/2$ denotes
the real part of $y$.
\par  The real scalar product $(.,.)_r$ in ${\cal A}_r$ we
extend on the algebra ${\cal A}_{r,v} $ as
$$(SP)\quad  <x,y>_{r,v} = x{\tilde y} =
\sum_{j,k=0}^{2^{v-r}-1} x_{2^rj}y_{2^rk} i_{2^rj} i_{2^rk}^*$$ for
each $x, y \in {\cal A}_{r,v}$, $x= \sum_{j=0}^{2^{v-r}-1} x_{2^rj}
i_{2^rj}$, $x_{2^rj}\in {\bf R}$ for each $j=0,...,2^{v-r}-1$.
Particularly, one gets $<x,y>_{0,v} = <x,y>_v$. In the case of the
complexified algebra $({\cal A}_{r,v})_{\bf C}$ the scalar product
is:
\par $$(SPC)\quad  <(a,b),(c,d)>_{r,v} = <(a,b),(c,d)> = (<a,c> - <b,d> , <a,d> + <b,c>) ,$$
for all $(a,b)$ and $(c,d)\in ({\cal A}_{r,v})_{\bf C}$.
\par We recall the doubling
procedure for the Cayley-Dickson algebra ${\cal A}_{r+1}$ from
${\cal A}_r$. Each Cayley-Dickson number $z\in {\cal A}_{r+1}$ is
written in the form $z=\xi + \eta {\bf l}$, where ${\bf l}^2=-1$,
${\bf l}\notin {\cal A}_r$, $\xi , \eta \in {\cal A}_r$. The
addition of such numbers is componentwise. The conjugate of any
Cayley-Dickson number $z$ is given by the formula:
\par $(M1)$ $z^* := \xi ^* - \eta {\bf l}$. \\
The multiplication in ${\cal A}_{r+1}$ is defined by the following
equation:
\par $(M2)$ $(\xi + \eta {\bf l})(\gamma +\delta {\bf l})=(\xi \gamma
-{\tilde {\delta }}\eta )+(\delta \xi +\eta {\tilde {\gamma }}){\bf l}$ \\
for each $\xi $, $\eta $, $\gamma $, $\delta \in {\cal A}_r$, $z :=
\xi +\eta {\bf l}\in {\cal A}_{r+1}$, $\zeta :=\gamma +\delta {\bf
l} \in {\cal A}_{r+1}$.
\par Using Formula $(M2)$  we get:
$(bi_{2^rk})(i_{2^rk}b)^* = (bi_{2^rk}) (b^*i_{2^rk}^*) = b^2 = (b^2
i_{2^rk}) i_{2^rk}^* = i_{2^rk} (i_{2^rk} b^2)^*$ for each $k\ge 1$
and $b\in {\cal A}_r$, since $i_j^*=-i_j$ for each $j\ge 1$. Another
useful identity is the following: $(i_si_{2^rj})i_{2^rk}^* = - (i_s
i_{2^rk}) i_{2^rj}^*$ for each $0\le s\le 2^r-1$ and $k\ne j$ with
$k\ge 1$ and $j\ge 1$, since $(i_si_{2^rj})i_{2^rk} = (i_s
i_{2^rk}^*) i_{2^rj}$. Certainly also the equality $(i_si_0)i_j^* +
(i_si_j)i_0^* =0$ holds for each $j\ge 1$ and $1\le s \le 2^r-1$,
since $i_0=1$. Therefore, Formulas $(SP)$ and $(4')$ together with
the latter identities imply:
$$(6)\quad <cBy,y>_v = \sum_{j=0}^{2^{v-r}-1} c <By_{2^rj},
y_{2^rj}>_{r,v}$$ for each $c\in {\cal A}_r$ and a twice
differentiable function $y$ with values in ${\cal A}_{r,v}$.
\par Relative to the complex scalar product given by Equality $(SPC)$
we decompose the operator $A$ (see $(4,4')$ above) in the form \par
$(7)$ $A = ({\bf i} \sigma ) ({\bf i} \sigma _1) + Q= - \sigma
\sigma _1 +Q$, \\ where $\sigma $, $\sigma _1$ and $Q$ are partial
differential operators of the first order, each complex number
$\alpha \in \bf C$ is presented as a real $2\times 2$ matrix.
Particularly, ${\bf i} = {{0~ ~ 1}\choose{ -1 ~ 0 }}$, ${\bf i}^*
={{0~ -1}\choose{ 1 ~ ~ 0 }}$.  Each subalgebra ${\sf g}_{j,k}$
constructed from two generators $i_j\ne i_k$ is associative,
consequently, $(wi_k)(w^*i_k^*) = w^2$ and $w((wi_k)i_k) = - w^2$
for each $w=w_0+w_ji_j$ with $w_0, w_j\in {\bf R}$. Therefore, we
can take $$(8)\quad \sigma f(z) = \sum_{j=1}^m
\sum_{k=m_1+...+m_{j-1}+1}^{m_1+...+m_j} {\bf a}_k(z)(\partial
f/\partial z_{2^rk}) [w_j^* i_{2^rk}^*]\mbox{  and}$$  $$(9)\quad
\sigma _1 f(z) = \sum_{j=1}^m
\sum_{k=m_1+...+m_{j-1}+1}^{m_1+...+m_j} {\bf a}_k(z)(\partial
f/\partial z_{2^rk})[w_ji_{2^rk}^*]$$  on the space of ${\cal
A}_r$-valued (super-)differentiable functions $f$ for $r\le 2$ or
real-valued functions $f$ for $3\le r$ of the ${\cal A}_{r,v}$
variable, since $({\bf i}i_k)^2= {\bf i}^2i_k^2=(-1)^2=1$ for each
$k\ge 1$, where $w_j^2=c_j$ for all $j$ and ${\bf a}_k^2(x)={\bf
b}_k(x)$ for each $k$ and $x$, $w_j\in {\cal A}_r$, ${\bf a}_k(x)\in
{\bf R}$ for all $k$ and $x$, $i_{2^rk}\in {\cal A}_{r,v}$,
$z_{2^rk} =x_k$, $z=\sum_k z_{2^rk} i_{2^rk} \in {\cal A}_{r,v}$,
$\partial f(z)/\partial z_{2^rk} = (df(z)/dz).i_{2^rk}$. For $b=
\partial ^2f/\partial z_{2^rk}^2$ and ${\bf l} = i_{2^rk}$ and $w\in
{\cal A}_r$ one has the identities: $(b(w{\bf l})) (w^*{\bf l}) =
((wb){\bf l})(w^*{\bf l}) = - w(wb) = - w^2b$ and $(((b{\bf
l})w^*){\bf l})w = (((bw){\bf l}){\bf l})w = - (bw)w = - bw^2$ in
the considered here cases. The operator $Q$ is given by the
equality: $$(10)\quad Qf(z) =  \sum_{j_1, j_2=1}^m
\sum_{k_1=m_1+...+m_{j_1-1}+1}^{m_1+...+m_{j_1}}
\sum_{k_2=m_1+...+m_{j_2-1}+1}^{m_1+...+m_{j_2}}$$  $$ {\bf
a}_{k_1}(z) \{ (\partial {\bf a}_{k_2}(z)/\partial z_{2^r{k_1}})
(\partial f/\partial z_{2^r{k_2}}) [w_{j_2} i_{2^r{k_2}}^*] \}
[w_{j_1}^* i_{2^r{k_1}}^*],$$  since ${\bf i} i_k = i_k {\bf i}$ in
the complexified Cayley-Dickson algebra $({\cal A}_v)_{\bf C}$ for
each $k$. The latter equality $(10)$ shows, that the differential
operator $Q$ is non-zero, when ${\bf a}_k(z)$ are non-constant
coefficients. \par If use $i_0=1$ and $\partial /\partial z_0$ also
one can write out d'Alembert's operator in our notation \par $(11)$
$\partial ^2/\partial z_0^2 - \sum_{j=1}^3 \partial ^2/\partial
z_j^2= ({\bf i}^*\partial /\partial z_0 + i_1\partial /\partial z_1
+ i_2 \partial /\partial z_2 + i_3 \partial /\partial z_3) ({\bf
i}\partial /\partial z_0 + i_1\partial /\partial z_1 + i_2 \partial
/\partial z_2 + i_3 \partial /\partial z_3)$. \par We recall, that
the Cayley-Dickson algebra ${\cal A}_r$ is power associative, that
is $z^kz^l=z^{k+l}$ for all natural numbers $k$ and $l$. But the
complexified Cayley-Dickson algebra $({\cal A}_r)_{\bf C}$ is not
power associative for $r\ge 3$, since the Cayley-Dickson algebra
${\cal A}_r$ is not associative for $r\ge 3$. Therefore, we do not
widely use the complexified Cayley-Dickson algebras, but we utilize
the Cayley-Dickson algebras ${\cal A}_v$ over the real field $\bf
R$, when something other will not be specified.
\par With these decomposition  of operators given by Equations $(7-9,11)$ the differential equation
$(1)$ can be integrated with the help of the non-commutative line
integration. We consider at first the partial differential equation
\par $(12)$ $\Upsilon f =g$ \\ on an open domain $U$ in ${\cal A}_v$, where
\par $(13)$ $\Upsilon f= \sum_{j=0}^{2^v-1}  (\partial f/\partial
z_j) [i_j^* {\psi }_j(z)]$, \\ $f$ and $g$ and ${\psi }_j(z)$ are
${\cal A}_v$-valued functions on the domain $U$ satisfying
Conditions 1$(D1,D2)$, where $g, ~ \psi _j\in C^0(U,{\cal A}_v)$ for
each $j$, particularly they may be ${\cal A}_v$
(super-)differentiable functions.
\par {\bf 3. Line
integration over Cayley-Dickson algebras.} Take any phrase
\par $(1)$ $\mu (z) =\sum_m \{ c_m, z^m \}_{q(m)} $ \\ corresponding to the function
$f$, where \par $ \{ c_m, z^m \}_{q(m)} = \{
c_{1,m_1}z^{m_1}...c_{k,m_k}z^{m_k} \}_{q(m)} $, \\ $q(m)$ is a
vector indicating on an order of multiplications in the curled
brackets, $c_{j,m_j}\in {\cal A}_v$ for each $j$, $m=(m_1,...,m_k)$,
$k\in {\bf N}$, $0\le m_j\in {\bf Z}$ for each $j$,
$z^k=(...((zz)z)...)z$.  We put for convenience $z^0=1$ in the
considered phrases. Though the symbol $z^0$ can be retained when
necessary to specify a branch of the line integral over the
Cayley-Dickson algebra ${\cal A}_r$ (see
\cite{ludfov,ludoyst,lutsltjms}). Using the shift $z\mapsto
(z-\mbox{}_0z)$ we can consider such series with the center at a
point $\mbox{}_0z$ instead of zero.
Then the derivative of the phrase is: \\
$(2)$ $d\mu (z)/dz = $ \\ $\sum_{m,l,j} \{
c_{1,m_1}z^{m_1}...c_{j-1,m_{j-1}}z^{m_{j-1}}
c_{j,m_j}((z^{m_j-l-1}I)z^l) c_{j+1,m_{j+1}}z^{m_{j+1}} ...
c_{k,m_k}z^{m_k} \}_{q(m)} $, \\ where $I$ denotes the unit
operator, so that $d\mu /dz$ is the operator valued derivative
function, $0\le l \le m_j-1$, $ ~ j=1,...,k$. From Equality $(2)$ it
follows that
\par $(3)$ $(d\mu (i_p x)/dx).1 = (d\mu (z)/dz).i_p = \partial \mu
(z)/\partial z_p$ for $z=i_p x$.
\par If $\gamma : [a,b]\to {\cal A}_r$ is a function, then
\par $V_a^b\gamma := \sup_P  |\gamma (t_{j+1})-\gamma (t_j)|$ \\ is
called the variation of $\gamma $ on the segment $[a,b]\subset {\bf
R}$, where the supremum is taken by all finite partitions $P$ of the
segment $[a,b]$, $P= \{ t_0=a<t_1<...<t_n =b \} $, $n\in {\bf N}$. A
continuous function $\gamma : [a,b]\to {\cal A}_r$ with the finite
variation $V_a^b\gamma <\infty $ is called a rectifiable path. It is
convenient to take the unit segment $[a,b]=[0,1]$ using a suitable
reparametrization.
\par We say that a function $\nu $ on $U$ is absolutely continuous
on $U$ if for each rectifiable path $\gamma : [0,1]\to {\cal A}_v$
for each $\epsilon >0$ and each $\tau \in [0,1]$ a positive number
$\delta >0$ exists so that $V_{\tau }^{\min (1, \tau +\delta )}\nu
(\gamma )<\epsilon $ and $V_{\max (0, \tau -\delta )}^{\tau } \nu
(\gamma )<\epsilon $. \par We call a function $\nu $ of bounded
variation on $U$ if for each rectifiable path $\gamma : [a,b]\to U$
the variation $V_a^b\nu (\gamma )<\infty $ is finite. The family of
all functions $\nu : U\to {\cal A}_v$ of bounded variation will be
denoted by ${\cal V}(U,{\cal A}_v)$.
\par The non-commutative line integral $\int_{\gamma }f(z)d\nu (z)$ along a
rectifiable path $\gamma : [0,1]\to U\subset {\cal A}_v$ for a
phrase $\mu $ and a given function $\nu $ of bounded variation  is
the limit by partitions $P=\{ 0=\tau _0<\tau _1<...<\tau _n=1 \} $
with their diameter $\delta (P)=\sup_j |\tau _{j+1}-\tau _j|$
tending to zero of integral sums
\par $\int_{\gamma }f(z)d\nu (z) := \lim_{\delta (P)\to 0} \sum_j
(d\kappa (z)/dz)|_{z= \gamma (\tau _j)}.[\nu (\gamma (\tau _{j+1}))
- \nu (\gamma (\tau _j))]$, \\ where $(d\kappa (z)/dz).1= \mu (z)$
for all $z\in U$. The notation \par ${\hat f}(z)= dg(z)/dz$ and
${\hat {\mu }} (z) = d\kappa (z)/dz$ \\ is also used, where $g(z)$
is a super-differentiable function to which the phrase $\kappa $
corresponds. \par If $f$ is a continuous function we fix for it a
sequence $f^n(z)$ of super-differentiable functions and their
phrases $\mu ^n (z)$ such that $f^n(z)$ converges to $f(z)$ on each
compact subset of the domain $U$, where $n\in {\bf N}$. The
non-commutative line integral has a continuous extension on the $\bf
R$-linear space, left and right ${\cal A}_v$ module, of continuous
functions $C^0(U,{\cal A}_v)$ for a marked function $\nu (z)$ of
bounded variation and a given rectifiable path $\gamma $:
\par $\int_{\gamma }f(z)d\nu (z) := \lim_{n\to \infty }
\int_{\gamma }f^n(z)d\nu (z)$. \par This means that the $\bf R$
homogeneous ${\cal A}_v$ additive operator ${\hat f}(z)$ is defined
for the continuous function $f$ in the sense of distributions:
\par $({\hat f};\nu ,\gamma ) := \int_{\gamma }f(z)d\nu (z)$ \\ for each
rectifiable path $\gamma $ in $U$ and every function $\nu (z)$ of
bounded variation. Particularly, $\nu (z)=id(z)=z$ on $U$ can also
be taken.
\par If $\nu $ and $f$ are super-differentiable functions such that
the derivative $d\nu (z)/dz$ of $\nu $ is the invertible $\bf R$
homogeneous ${\cal A}_v$ additive operator for each $z\in U$, then a
super-differentiable solution of the differential equation
\par $(dg(z)/dz).(d\nu (z)/dz) = dq((\nu (z))/dz$ \\ on $U$ exists, since
$dz/d\nu =(d\nu /dz)^{-1}$. That is, $(dg(z)/dz).d\nu (z) = (dq(\nu
(z))/dz).dz$. Therefore,
\par $\int_{\gamma } f(z)d\nu (z) = \lim_{\delta (P)\to 0}
\sum_j (dq(\nu )/d\nu )|_{\nu =\nu (\gamma (\tau _j))} .[\nu (\gamma
(\tau _{j+1})) -
\nu (\gamma (\tau _j))] = \int_{\nu (\gamma )} p(y)dy $, \\
where $p(\nu ) = (dq(\nu )/d\nu ).1$ (see also Theorems 2.11 and
2.13 in \cite{luviralgoclanl}).
\par A function $\nu : U\to {\cal A}_v$ is called piecewise
continuous or differentiable or super-differentiable on a domain $U$
in the Cayley-Dickson algebra $U$, if a family of open or canonical
closed subsets $U_j$ of $U$ exists so that each restriction $\nu
|_{U_j}$ is continuous or differentiable or super-differentiable
respectively, where $U=\bigcup_jU_j$ and $U_j\cap U_k=\partial
U_j\cap \partial U_k$ for each $j\ne k$, $\partial U_j := cl
(U_j)\setminus Int (U_j)$, $cl (U_j)$ denotes the closure of $U_j$
in ${\cal A}_v$ and $Int (U_j)$ denotes the interior of $U_j$ in
${\cal A}_v$.
\par  If $f$ is
a continuous function and $\nu $ is a function of bounded variation
for which the limits $\lim_nf^n=f$ and $\nu =\lim_n\nu ^n$ uniformly
converge on each compact subset of $U$ and phrases $\mu ^n$ of $f^n$
and $\epsilon ^n$ of $\nu ^n$ are specified, where $f^n$ are
super-differentiable functions and $\nu ^n$ are piecewise
super-differentiable functions on $U$  so that $\lim_n V_0^1(\nu
^n(\gamma )-\nu (\gamma )) =0$ for each rectifiable path $\gamma $
in $U$, then
\par $\int_{\gamma } f(z)d\nu (z) = \lim_n \int_{\gamma } f^n(z)d\nu
^n(z) = \lim_n \int_{\nu ^n (\gamma )} p ^n (y)dy = \int_{\nu
(\gamma )} p(y)dy$, \\ where $p(y)=\lim_n p^n(y)$. This means that
under rather general conditions the line integral of the type
$\int_{\gamma }f(z)d\nu (z)$ relative to the function $\nu $ of
bounded variation reduces to the usual non-commutative line integral
$\int_{\eta } p(y)dy$, where $\eta = \nu (\gamma )$.
\par  Take the
branch of the non-commutative line integral prescribed by the left
algorithm (see \S 2 in \cite{ludoyst,ludfov}). The real algebra
${\sf g}_{k,l,s}$ formed from the generators $i_j$, $i_k$ and $i_s$
is alternative. Each rectifiable path can be presented as the limit
of rectifiable paths consisting of joined segments parallel to the
straight lines $i_j{\bf R}$ with respective $j$. We certainly have
$(i_qi_p)i_p = - i_q$ for each $p\ge 1$ and $(i_qi_0)i_0=i_q$ for
each $q\ge 0$. \par For each $j=0,...,2^r-1$ the $\bf R$- linear
projection operator $\pi _j : {\cal A}_r\to {\bf R}i_j$ exists due
to Formulas 1$(1-3)$ so that $\pi _j(z)=i_jz_j=z_ji_j$: \par $(P1)$
$\pi _j(z) = (- i_j (zi_j) - (2^r-2)^{-1} \{ -z
+\sum_{k=1}^{2^r-1}i_k(zi_k^*) \} )/2$
\\ for each $j=1,2,...,2^r-1$, $$(P2)\quad \pi _0(z) = (z+
(2^r-2)^{-1} \{ -z + \sum_{k=1}^{2^r-1}i_k(zi_k^*) \} )/2,$$  where
$2\le r\in \bf N$.

\par {\bf 4. Line anti-derivatives over Cayley-Dickson algebras.}
\par {\bf Theorem.} {\it Let a first order partial differential
operator $\Upsilon $ be given by Equation 2$(13)$ with real-valued
continuous functions $\psi _j(z)\in C^0(U,{\cal A}_v)$ for each $j$
such that $\psi _j(z)\ne 0$ for each $z\in U$ and each $j=0,...,n$,
where a domain $U$ satisfies Conditions 1$(D1,D2)$, $ ~ \mbox{}_0z$
is a marked point in $U$, $1<n< 2^v$, $2\le v$. Then a line integral
${\cal I}_{\Upsilon }: C^0(U,{\cal A}_v)\to C^1(U,{\cal A}_v)$,
${\cal I}_{\Upsilon }f(z) := \mbox{}_{\Upsilon }\int_{\mbox{}_0z}^z
f(y)dy$ on $C^0(U,{\cal A}_v)$ exists so that
\par $(1)$ $\Upsilon {\cal I}_{\Upsilon }f(z) = f(z)$ \\ for each
$z\in U$; this anti-derivative is $\bf R$-linear (or $\bf
H$-left-linear when $v=2$):
\par $(2)$ ${\cal I}_{\Upsilon }[af(z)+ bg(z)] = a{\cal
I}_{\Upsilon }f(z) + b{\cal I}_{\Upsilon }g(z)$ \\
for any real constants $a, b\in {\bf R}$ (or $a, b \in {\bf H}$ for
$v=2$) and continuous functions $f, g \in C^0(U,{\cal A}_v)$. If
there is a second anti-derivative ${\cal I}_{\Upsilon ,2}f(z)$, then
${\cal I}_{\Upsilon }f(z)- {\cal I}_{\Upsilon ,2}f(z)$ belongs to
the kernel $ker (\Upsilon )$ of the operator $\Upsilon $. }
\par {\bf Proof.}  Using  the multiplication on the marked doubling generator $i_{2^v}$
from the right we have $$(3)\quad [\sum_{j=0}^{k-1} i_j(\partial
g(z)/\partial z_j)]i_{2^v}\psi _j(z) = \sum_{j=0}^{k-1}(\partial
g(z)/\partial z_j)i_{j+2^v}\psi _j(z),$$ where $i_ji_{2^v} =:
i_{j+2^v}$ for each $0\le j\le 2^v-1$, $2\le v$. On the other hand,
$\sum_{j=0}^{k-1} i_j(\partial g(z)/\partial z_j)\psi _j(z) =
[\sum_{j=0}^{k-1} (\partial g(z)/\partial z_j)^*i_j^*\psi _j(z)]^*$,
since $\psi _j(z)$ is real for each $j$ and $z$. Therefore, it is
sufficient to consider the first-order partial differential operator
of the form:
\par $(4)$  $\Upsilon g (z) = \sum_{j=1}^n  (\partial g/\partial
z_j) i_j^*\psi _j(z) $ \\ on the $\bf R$-linear space $C^1(U,{\cal
A}_v)$ of all continuously differentiable functions $g: U\to {\cal
A}_v$ by real variables $z_0,...,z_{2^v-1}$, where $0<n\le 2^v-1$.
The space of super-differentiable functions is everywhere dense in
$C^0(U,{\cal A}_v)$ and the line integral has the continuous
extension on $C^0(U,{\cal A}_v)$ along any continuous rectifiable
path in $U$. Therefore, we take the space of super-differentiable
functions and then take the continuous extension of ${\cal
I}_{\Upsilon }$ on $C^0(U,{\cal A}_v)$ such that \par $\lim_l{\cal
I}_{\Upsilon }f^l= {\cal I}_{\Upsilon }\lim_lf^l={\cal I}_{\Upsilon
}f$ \\ for a sequence $f^l$ of super-differentiable functions
uniformly converging to $f$ on compact sub-domains $V$ in $U$, where
${\cal I}_{\Upsilon }f^l$ is described below. Each function $\psi
_j(z)$ is continuous and each function $\nu _j(z)$ is continuously
differentiable on $U$ (see also below), consequently, the integral
$\int_{\gamma } f^l(y)d\nu _j(y)$ is continuously differentiable by
$z=\gamma (1)$ (i.e. by each real variable $z_k$) and their sequence
by $l$ uniformly converges on each compact sub-domain $V$ in $U$.
Therefore, from $\Upsilon {\cal I}_{\Upsilon }f^l=f^l$ for each
natural number $l\in {\bf N}$ we get
\par $\Upsilon {\cal I}_{\Upsilon }f=\Upsilon {\cal I}_{\Upsilon
}\lim_l f^l = \lim_l \Upsilon {\cal I}_{\Upsilon }f^l =
\lim_lf^l=f$, \\ since the sequence $\{ {\cal I}_{\Upsilon
}f^l(z)|_V: ~ l \}$ is fundamental in $C^1(V,{\cal A}_v)$ for each
compact sub-domain $V$ in $U$ and ${\cal I}_{\Upsilon }f(z)\in
C^1(U,{\cal A}_v)$.
\par Consider the left algorithm of a calculation of the line
integral over the Cayley-Dickson algebra ${\cal A}_v$ (see \S 3 and
references therein). We shall seek an anti-derivative in the form:
$$(5)\quad \mbox{}_{\Upsilon }\int_{\mbox{}_0z}^z
f(y)dy := n^{-1} \sum_{j=1}^n (\int_{\mbox{}_0z}^z f(y)d\nu
_j(y))i_j $$ and use the homotopy theorem in the domain $U$
satisfying conditions 1$(D1,D2)$ so that $\gamma $ is a continuous
rectifiable path joining points $\mbox{}_0z=\gamma (0)$ and
$z=\gamma (1)$ (see \cite{ludfov,ludoyst,ludeoc}). Moreover, a
branch of the anti-derivative operator ${\cal I}_{\Upsilon }f(z)$
can be chosen such that it can be expressed with the help of a
non-commutative line integral.
\par   In view of Theorem 2.11
\cite{luviralgoclanl} and \S 3 we get
\par $(6)$ $(\partial  (\int_{\mbox{}_0z}^z f(y)d\nu _j(y)) /\partial z_k) =
({\hat f} (z).[d\nu _j(z)/dz_k]) $ \\
(see also the chain rule over the Cayley-Dickson algebra in
\cite{ludfov,ludoyst,ludeoc}).
\par Next we need some identities in the Cayley-Dickson algebra.
Each Cayley-Dickson number has the decomposition: $z = z_0i_0 +...+
z_{2^v-1}i_{2^v-1}$, where $z_0,...,z_{2^v-1}\in {\bf R}$, $z\in
{\cal A}_v$. To establish the identity \par $(7)$ $(ay)z^* + (az)y^*
= a 2 Re (yz^*)$ \\ for any $a, y, z\in {\cal A}_v$ it is sufficient
to prove it for any three basic generators of the Cayley-Dickson
algebra ${\cal A}_v$, since the real field ${\bf R}$ is its center,
while the multiplication in ${\cal A}_v$ is distributive
$(a+y)z=az+yz$ and $((\alpha a) (\beta y)) (\gamma z^*) = (\alpha
\beta \gamma ) ((ay)z^*)$ for all $\alpha , \beta , \gamma \in {\bf
R}$ and $a, y, z \in {\cal A}_v$. If $a=i_0$, then $(7)$ is evident,
since $yz^* + zy^* = yz^* + (yz^*)^* = 2 Re (yz^*)$. If either
$y=i_0$, then $(ay)z^*+ (az)y^*= az^* + az = a 2 Re (z)= a 2 Re
(yz^*)$. Analogously for $z=i_0$. For three purely imaginary
generators $i_p, i_s, i_k$ consider the minimal Cayley-Dickson
algebra $\Phi = alg_{\bf R} (i_p, i_s, i_k)$ over the real field
generated by them. If it is associative, then it is isomorphic with
either the complex field $\bf C$ or the quaternion skew field $\bf
H$, so that $(ay)z^* + (az)y^* = a(yz^*+zy*) = a 2 Re (yz^*)$. If
the algebra $\Phi $ is isomorphic with the octonion algebra, then we
use Formulas 2$(M1,M2)$ for either $a, y\in {\bf H}$ and $z={\bf l}$
or $a, z\in {\bf H}$ and $y={\bf l}$. This gives $(7)$ in all cases,
since the algebra $alg_{\bf R} (i_p,i_s)$ with two basic generators
$i_p$ and $i_s$ is always associative. Particularly, if $y=i_s\ne
z=i_k$, then the result is zero.
\par Using $(7)$ we get more generally, that
\par $(8)$ $((ay)z^*)b^* + ((az)y^*)b^* = (a 2 Re (yz^*))b^*=
(ab^*) 2 Re (yz^*) ,$ \\ consequently,
\par $(9)$ $((ay)z^*)b^* + ((az)y^*)b^* + ((by)z^*)a^* + ((bz)y^*)a^*=
4 Re (ab^*) Re (yz^*) $ \\
for any Cayley-Dickson numbers $a, b, y, z\in {\cal A}_v$.
\par We shall take unknown functions $\nu _j(z)\in {\cal A}_v$ as solutions of
the system of linear partial differential equations by real
variables $z_k$:
\par $(10)$ $\partial \nu _j(z)/\partial z_j = 1/\psi _j(z)$ for all
$1\le j\le n$ and $z\in U$;
\par $(11)$ $\psi _k(z) \partial \nu _j(z)/\partial z_k = \psi
_j(z)\partial \nu _k(z)/\partial z_j $ for all $1\le j<k\le n$ and
$z\in U$. Each function $\nu _j(z)$ can be written as $\nu _j(z)=
\sum_{l=0}^{2^v-1} \nu _{j,l}(z)i_l$ with real-valued components
$\nu _{j,l}(z)$. Practically, it is sufficient to consider non-zero
$\nu _{j,l}(z)$ for $l=1,...,n$. Thus using the generators
$i_0,....,i_{2^v-1}$ the system can be written in the real form.
 This system has a non-trivial $C^1$ solution $\nu _j(z)$ for each $j$
(see \S 12.2 \cite{matveevb}, particularly, in the class of
super-differentiable functions for super-differentiable $\psi _j(z)$
see also \cite{lujmsnfpcd,ludeoc}). In System $(10,11)$ functions
$\psi _j$ are real and coordinates are real, consequently, a
solution $\{ \nu _j(z): ~ j \} $ may be chosen real-valued.
\par From Identities 3$(2,3)$ and  $(6,9-11)$ we infer that $$(12)\quad
\sum_{j\ne k\ge 1} [(\partial  (\int_{\mbox{}_0z}^z f(y) d\nu _j(y))
/\partial z_k)i_j]i_k^*\psi _k(z) =$$ $$\sum_{1\le j< k\le n} \{
[({\hat f}(z).(\partial \nu _j(z)/\partial z_k))i_j]i_k^*\psi _k(z)
+ [({\hat f}(z).(\partial \nu _k(z)/\partial z_j))i_k]i_j^*\psi
_j(z) \} =0\mbox{  and}$$ $$(13)\quad \sum_{j=1}^n [(\partial
(\int_{\mbox{}_0z}^z f(y) d\nu _j(y) )/\partial z_j)i_j]i_j^* \psi
_j(z) = n f(z),$$  since $\sum_{j=1}^n i_ji_j^*=n$ and $n$ is some
fixed natural number for the domain $U$, $ ~ {\hat f}(z).x=f(z)x$
for each real number $x$, $ ~ (zi_j)i_j^* = z$ for each $z\in {\cal
A}_v$, where $\hat f$ is the operator corresponding to $d\kappa
(z)/dz$, when $f$ is in the $z$-representation $\mu $ (see the
notation in \S 3). Using Formulas $(4,5,12,13)$ we get Formula
$(1)$.
\par From the identity $\int_{\gamma }a\mu dz=a\int_{\gamma }\mu
dz$ for a suitable branch of the line integral given by the left
algorithm and for each non-trivial phrase $\mu $ and constants $a,
b\in {\bf R}$ for $v\ge 3$ or $a, b \in {\bf H}$ for $v=2$ (see the
rules in \cite{ludfov,ludoyst,lujmsnfpcd,lutsltjms}) we get Formula
$(2)$.
\par Since $\Upsilon ({\cal I}_{\Upsilon }f(z)- {\cal I}_{\Upsilon
,2}f(z))=0$, the difference $({\cal I}_{\Upsilon }f(z)- {\cal
I}_{\Upsilon ,2}f(z))$ belongs to the kernel $ker (\Upsilon ) =
\Upsilon ^{-1}(0)$, where $\Upsilon : C^1(U,{\cal A}_v)\to
C^0(U,{\cal A}_v)$.

\par {\bf 4.1. Example.} If $\psi _j$ depends only on $z_j$ for
each $j$, there exists a $C^1$ differentiable change of variables
$\zeta = \zeta (z)$ so that $\partial g(\zeta )/\partial \zeta _j =
(\partial g(z)/\partial z_j)\psi _j(z)$ for each differentiable
function $g: U\to {\cal A}_v$ by real variables $z_0,...,z_{2^v-1}$
on $U$, where
\par $(1)$ $(\partial z_k/\partial \zeta _j) = \delta_{j,k} \psi _k(z)$
\\ for all $j$ and $k$, $ ~ \delta _{j,j}=1$, while $\delta
_{j,k}=0$ for each $j\ne k$. We take new functions $\mbox{}_jg$
satisfying the equation:
\par $(2)$ $\mbox{}_jg(i_j^*z)=g(z)$ for each $z\in U$ and all $j$.
We also put \par $(3)$ $\eta _j(z)=i_j^*z$. \\
The multiplication of generators implies that $i_j^*(i_jz)=z$ for
all $j=0,...,2^v-1$ and $z\in {\cal A}_v$. Therefore, from Equations
$(1,2)$ we deduce that
\par $(4)$ $(dg(z)/dz).i_j = (d\mbox{}_jg(\eta _j)/d\eta _j).
[(d\eta _j/dz).i_j]=(d\mbox{}_jg(\eta _j)/d\eta
_j).1=(d\mbox{}_kg(\eta _k)/d\eta _k).[i_k^*i_j],$ \\ since $(d\eta
_j/dz).i_j = i_j^*i_j=1$ for each $j$. Then we take the integral
$$(5)\quad \mbox{}_{\Upsilon }\int_{\mbox{}_0z}^z g(y)dy := n^{-1} \sum_{j=1}^n
\int_{\mbox{}_0z}^z \mbox{}_jg(\eta _j(y))i_j d\eta _j(y),$$  since
$\int_{\mbox{}_0z}^z \mbox{}_jg(\eta _j(y))i_j d\eta _j(y)
=(\int_{\mbox{}_0z}^z \mbox{}_jg(\eta _j(y)) d\eta _j(y))i_j.$ \\
\par Mention that generally $\Upsilon (f(z)b)$ may be not equal to
$(\Upsilon f(z))b$ for a constant $b\in {\cal A}_v\setminus {\bf R}$
and a function $f\in C^1(U,{\cal A}_v)$ with $v\ge 2$, since the
Cayley-Dickson algebra is non-commutative.
\par This theorem can be generalized in the following manner
encompassing wider class of partial differential operators of the
first order over Cayley-Dickson algebras.
\par {\bf 5. Theorem.} {\it Suppose that the first order
partial differential operator $\Upsilon $ is given by the formula
\par $(1)$ $\Upsilon f= \sum_{j=0}^n  (\partial f/\partial
z_j) {\phi }_j^*(z)$, \\ where ${\phi }_j(z)\ne \{ 0 \}$ for each
$z\in U$ and $\phi _j(z)\in C^0(U,{\cal A}_v)$ for each $j=0,...,n$
such that $Re (\phi _j(z)\phi _k^*(z))=0$ for each $z\in U$ and each
$0\le j\ne k\le n$, where a domain $U$ satisfies Conditions
1$(D1,D2)$, $ ~ \mbox{}_0z$ is a marked point in $U$, $1<n< 2^v$,
$2\le v$. Suppose also that the system $ \{ \phi _0(z),...,\phi
_n(z) \} $ is for $n=2^v-1$, or can be completed by Cayley-Dickson
numbers $\phi _{n+1}(z),...,\phi _{2^v-1}(z)$, such that $(\alpha )$
$ ~ alg_{\bf R} \{ \phi _j(z),\phi _k(z), \phi _l(z) \} $ is
alternative for all $0\le j, k, l \le 2^v-1$ and $(\beta )$
$alg_{\bf R} \{ \phi _0(z),...,\phi _{2^v-1}(z) \} = {\cal A}_v$ for
each $z\in U$. Then a line integral ${\cal I}_{\Upsilon }:
C^0(U,{\cal A}_v)\to C^1(U,{\cal A}_v)$, ${\cal I}_{\Upsilon }f(z)
:= \mbox{}_{\Upsilon }\int_{\mbox{}_0z}^z f(y)dy$ on $C^0(U,{\cal
A}_v)$ exists so that
\par $(2)$ $\Upsilon {\cal I}_{\Upsilon }f(z) = f(z)$ \\ for each
$z\in U$; this anti-derivative is $\bf R$-linear (or $\bf
H$-left-linear when $v=2$). If there is a second anti-derivative
${\cal I}_{\Upsilon ,2}f(z)$, then ${\cal I}_{\Upsilon }f(z)- {\cal
I}_{\Upsilon ,2}f(z)$ belongs to the kernel $ker (\Upsilon )$ of the
operator $\Upsilon $. }
\par {\bf Proof.} We shall demonstrate that a branch of
the anti-derivative operator ${\cal I}_{\Upsilon }f(z)$ can be
chosen such that it can be expressed with the help of a
non-commutative line integral from \S 3. Using the technique of \S 4
we can consider the case of purely imaginary $\phi _j(z)$ for all
$z\in U$ and $j=0,...,n$. We seek an anti-derivative operator in the
form: $$(3)\quad \mbox{}_{\Upsilon }\int f(z)dz = (n+1)^{-1}
\sum_{j=0}^n \int_{\mbox{}_0z}^z q(z)d\nu _j(z).$$ For finding
unknown functions $q$ and $\nu _j$, $j=0,...,n$ we impose the
following conditions:
\par $(4)$ $({\hat q}(z).[\partial \nu _j(z)/\partial z_j])\phi
_j^*(z) = f(z)$ for each $j=0,...,n$ and
\par $(5)$ $({\hat q}(z).[\partial \nu _j(z)/\partial z_k])\phi
_k^*(z) + ({\hat q}(z).[\partial \nu _k(z)/\partial z_j])\phi
_j^*(z) =0$ for all $0\le j<k\le n$. \par As in \S 4 it is
sufficient to consider the case of a locally analytic
(super-differentiable) function $f$ using the limit transition. The
function $f$ is given on $U$ and it defines the operator $\hat f$ on
$U$, i.e. its phrase $\hat \mu $ is prescribed by the left algorithm
for a given phrase $\mu $ of $f$ (see
\cite{ludfov,ludoyst,lujmsnfpcd,lutsltjms}). The operator $\hat q$
means that a function $g$ and a phrase $\kappa $ of $g$ exist such
that
\par ${\hat q} (z) = dg(z)/dz$, $ ~ {\hat q}(z).1=q(z)$ for each $z\in U$.
\par In accordance with the
conditions of this theorem the algebra $~ alg_{\bf R} (\phi
_j(z),\phi _k(z))$ is alternative for all $0\le j\le k\le n$ and
$z\in U$. Therefore, due to Condition $(\beta )$ Equations $(4,5)$
take the form:
\par $(6)$ $(dg(z)/dz).[\partial \nu _j(z)/\partial z_j] = f(z)(1/\phi
_j^*(z))$ for each $j=0,...,n$ and
\par $(7)$ $((dg(z)/dz).[\partial \nu _j(z)/\partial z_k])\phi
_k^*(z) + ((dg(z)/dz).[\partial \nu _k(z)/\partial z_j])\phi _j^*(z)
=0$ for all $0\le j<k\le n$. \par Solutions of this system exist
(see \cite{lujmsnfpcd,ludeoc}). To be more concrete we impose
additional relations:
\par $(8)$ $\partial \nu _j(z)/\partial z_j = \phi _j(z)$ for all
$j=0,....,n$ and $z\in U$, \\ consequently, the system of partial
differential equations $(6)$ becomes: \par $(9)$ $(dg(z)/dz).\phi
_j(z)=f(z)(1/\phi _j^*(z))$ for each $j=0,...,n$, \\ since $ ~
alg_{\bf R} \{ \phi _j(z),\phi _k(z), \phi _l(z) \} $ is alternative
for all $0\le j, k, l \le 2^v-1$ and $alg_{\bf R} \{ \phi
_0(z),...,\phi _{2^v-1}(z) \} = {\cal A}_v$ for each $z\in U$ so
that each Cayley-Dickson number $\xi \in {\cal A}_v$ has the
decomposition $\xi = \xi _0\phi _0(z)+...+\xi _{2^v-1}\phi
_{2^v-1}(z)$ with real coefficients $\xi _0,...,\xi _{2^v-1}\in {\bf
R}$.
\par Solving the latter system $(9)$ one gets the function $g(z)$ on
$U$. Substituting the known function $g$ in System $(6,7)$ one gets
a $C^1$ solution $\nu _0(z)$,...,$\nu _n(z)$ on $U$; or a
super-differentiable solution, when $\phi _j(z)$ for each $j$ and
$f(z)$ are super-differentiable on $U$. Mention that the function
$g$ depends $\bf R$-linearly on $f$, since the system of equations
which was considered above is linear by $f$ and $g$. Thus the
operator $\hat q$ depends $\bf R$-linearly on $f$.
\par Using Formulas $(4,5)$ and 4$(6,9)$ we deduce that $$(10)\quad
\sum_{j\ne k\ge 0} [\partial  (\int_{\mbox{}_0z}^z q(y) d\nu _j(y))
/\partial z_k]\phi _k^*(z) =$$ $$\sum_{0\le j< k\le n} \{ [{\hat
q}(z).(\partial \nu _j(z)/\partial z_k)]\phi _k^*(z) + [{\hat
q}(z).(\partial \nu _k(z)/\partial z_j)]\phi _j^*(z) \} =0\mbox{
and}$$ $$(11)\quad \sum_{j=0}^n [\partial (\int_{\mbox{}_0z}^z q(y)
d\nu _j(y) )/\partial z_j]\phi _j^*(z) =\sum_{j=0}^n [ {\hat
q}(z).(\partial \nu _j(z)/\partial z_j)]\phi _j^*(z)= (n+1) f(z),$$
since $Re (\phi _j(z)\phi _k^*(z))=0$ for each $z\in U$ and each
$0\le j\ne k\le n$.
\par  The rest of the proof is analogous to that of Theorem 4.
\par {\bf 6. Corollary.} {\it Let suppositions of Theorem 5 be satisfied
so that $\phi _j(z)=\omega (z;i_j)\psi _j(z)$ for each $z\in U$,
where $\omega $ is an $\bf R$-linear automorphism $\omega : {\cal
A}_v\to {\cal A}_v$ mapping the standard base of generators $\{ i_j
\} $ into a base of generators $\{ \omega (z;i_j): j=0,...,2^v-1 \}
$, $|\omega (z;i_j)| =1$, where $\psi _j(z)$ satisfies conditions of
theorem 4 for each $j=0,...,n$. Then the first order differential
operator 5$(1)$ has an anti-derivative ${\cal I}_{\Upsilon }$ on
$C^0(U,{\cal A}_v)$. Two anti-derivatives of Theorems 4 and 5 under
these suppositions are related with the help of the automorphism
$\omega $.}
\par {\bf Proof.} This follows immediately from Theorem 5.
It remains to find a relation between two anti-derivatives for two
different partial differential operators:
\par $(1)$ $\Upsilon _{\omega } f= \sum_{j=0}^n  (\partial f/\partial
z_j) {\phi }_j^*(z)$ \\ and $\Upsilon $ given by equation 2$(13)$.
\par  For each Cayley-Dickson number
$z=z_0i_0+...+z_{2^v-1}i_{2^v-1}\in {\cal A }_v$ its image is
$\omega (y;z) = z_0N_0+z_1N_1+...+z_{2^v-1}N_{2^v-1}$, consequently,
$\omega (y;z^*) = [\omega (y;z)]^*$, where $z_j\in {\bf R}$, $ ~
N_j=N_j(y) := \omega (y;i_j)$ for each $j$. Particularly, $N_0=i_0$,
since $i_0i_j=i_j$ and $\omega (y;i_j)=\omega (y;i_0i_j)=\omega
(y;i_0)\omega (y;i_j)$ for each $j$ and $y$. Therefore, $\omega
(y;x)=x$ for each real number $x\in {\bf R}$, since $\omega (y;1)=1$
and the mapping $\omega (y;*)$ is $\bf R$-linear by the second
argument, $1=i_0$. Therefore, applying the automorphism $\omega $ we
deduce that
\par $(2)$ $\Upsilon _{\omega }f(z) = \omega (z;\Upsilon s(z))$, \\
where $\omega (z;s(z))=f(z)$ for each $z\in U$, that is $s(z)=\omega
_2^{-1}(z;f(z))$, $ ~ \omega _2^{-1}(z;*)$ denotes the inverse
automorphism by the second argument for $z\in U$. Let us take the
function $f(z)= \mbox{}_{\Upsilon _{\omega }} \int_{\mbox{}_0z}^z
g(y)dy$, where $g(z)$ is a continuous function. Then $\Upsilon
_{\omega } f(z)=g(z)$ for each $z\in U$ and from $(2)$ and 5$(1,2)$
one gets
\par $(3)$ $\omega _2^{-1}(z;g(z)) = \Upsilon \omega _2^{-1}(z;
\mbox{}_{\Upsilon _{\omega }}\int_{\mbox{}_0z}^z g(y)dy)=\Upsilon
\mbox{}_{\Upsilon }\int_{\mbox{}_0z}^z \omega _2^{-1}(y;g(y))dy,$
consequently, applying $\mbox{}_{\Upsilon} \int$ and $\omega (z;*)$
one also gets
\par $(4)$ $\mbox{}_{\Upsilon _{\omega }}\int_{\mbox{}_0z}^z g(y)dy =
\omega (z; \mbox{}_{\Upsilon }\int_{\mbox{}_0z}^z \omega _2^{-1}
(y;g(y))dy$ \\ for each continuous function $g$ on $U$.
\par {\bf 6.1. Remark.} If in Theorem 5 drop Conditions $(\alpha
,\beta )$, then partial differential equations 5$(4,5)$ will be hard
to resolve. \par To specify the anti-derivative operator ${\cal
I}_{\Upsilon }$ in Theorems 4 and 5 more concretely it is possible
to choose a family of rectifiable continuous paths (or $C^1$ paths)
$\{ \gamma ^z: ~ z\in U \} $ such that $\gamma ^z(0)=\mbox{}_0z$ and
$\gamma ^z(1)=z$ and $\lim_{z\to y} \sup_{\tau \in [0,1]} |\gamma
^z(\tau ) - \gamma ^y(\tau )|=0$.

\par Another more rigorous procedure is in providing a foliation
of a domain $U$ by locally rectifiable paths $\{ \gamma ^{\alpha }:
~ \alpha \in \Lambda \} $, where $\Lambda $ is a set. We take for
definiteness a canonical closed domain $U$ in ${\cal A}_v$
satisfying Conditions 1$(D1,D2)$. \par A path $\gamma : <a,b>\to U$
is called locally rectifiable, if it is rectifiable on each compact
segment $[c,e]\subset <a,b>$, where $<a,b> =[a,b] := \{ t\in {\bf
R}: ~ a\le t \le b \} $ or $<a,b> =[a,b) := \{ t\in {\bf R}: ~ a\le
t < b \} $ or $<a,b> =(a,b] := \{ t\in {\bf R}: ~ a< t \le b \} $ or
$<a,b> =(a,b) := \{ t\in {\bf R}: ~ a< t < b \} $. \par A domain $U$
is called foliated by rectifiable paths $\{ \gamma ^{\alpha }: ~
\alpha \in \Lambda \} $ if $\gamma : <a_{\alpha },b_{\alpha }> \to
U$ for each $\alpha $ and it satisfies the following three
conditions: \par $(F1)$ $\bigcup_{\alpha \in \Lambda } \gamma
_{\alpha }(<a_{\alpha },b_{\alpha }>) =U$ and \par $(F2)$ $\gamma
_{\alpha }(<a_{\alpha },b_{\alpha }>)\cap \gamma ^{\beta }(<a_{\beta
},b_{\beta }>) = \emptyset $ for each $\alpha \ne \beta \in \Lambda
$. \\ Moreover, if the boundary $\partial U = cl (U)\setminus Int
(U)$ of the domain $U$ is non-void then \par $(F3)$ $\partial U =
(\bigcup_{\alpha \in \Lambda _1} \gamma ^{\alpha }(a_{\alpha }))\cup
(\bigcup_{\beta \in \Lambda _2} \gamma ^{\beta }(b_{\beta }))$, \\
where $\Lambda _1 = \{ \alpha \in \Lambda : <a_{\alpha },b_{\beta
}>= [a_{\alpha },b_{\beta }> \} $, $\Lambda _2 = \{ \alpha \in
\Lambda : <a_{\alpha },b_{\beta }>= <a_{\alpha },b_{\beta }] \} $.
For the canonical closed subset $U$ we have $cl (U)=U=cl (Int(U))$,
where $cl (U) $ denotes the closure of $U$ in ${\cal A}_v$ and $Int
(U)$ denotes the interior of $U$ in ${\cal A}_v$. For convenience
one can choose $C^1$ foliation, i.e. each $\gamma ^{\alpha }$ is of
class $C^1$. When $U$ is with non-void boundary we choose a
foliation family such that $\bigcup_{\alpha \in \Lambda } \gamma
(a_{\alpha }) =\partial U_1$, where a set $\partial U_1$ is open in
the boundary $\partial U$ and so that $w|_{\partial U_1}$ would be a
sufficient initial condition to characterize a unique branch of an
anti-derivative $w={\cal I}_{\Upsilon }f$.
\par When $\partial U\ne \emptyset $
a marked point $\mbox{}_0z$ can be chosen on the boundary $\partial
U$ and each point on the boundary can be joined by a rectifiable
path in $U$ with $\mbox{}_0z$. This foliation is justified by the
formula:
\par $\int_{\gamma }f(z)d\nu
(z)=\int_{\gamma ^1}f(z)d\nu (z)+ \int_{\gamma ^2}f(z)d\nu (z)$\\
for each continuous function $f$ on $U$ and each function $\nu $ of
bounded variation on $U$, for any rectifiable paths $\gamma ^1:
[a_1,b_1]\to U$ and $\gamma ^2: [a_2,b_2]\to U$ so that
$a=a_1<b_1=a_2<b_2=b$ while $\gamma : [a,b]\to U$ is given piecewise
as $\gamma (t)=\gamma ^1(t)$ for each $t\in [a_1,b_1]$ and $\gamma
(t)=\gamma ^2(t)$ for each $t\in [a_2,b_2]$. Thus instead of
$\int_{\mbox{}_0z}^zf(z)d\nu (z)$, i.e. $\int_{\gamma }f(z)d\nu (z)$
with $\gamma (a)=\mbox{}_0z$ and $\gamma (b)=z$, we take
$\int_{\gamma ^{\alpha }|_{[c,e]}}f(z)d\nu (z)$ for any
$[c,e]\subset <a_{\alpha },b_{\alpha }>$. If $\lim_{c\to a_{\alpha
}, e\to b_{\alpha }} \int_{\gamma ^{\alpha }|_{[c,e]}}f(z)d\nu (z)$
converges we denote it by $\int_{\gamma ^{\alpha }}f(z)d\nu (z)$ and
take instead of the family $\{ \int_{\gamma ^{\alpha
}|_{[c,e]}}f(z)d\nu (z): ~ [c,e]\subset <a_{\alpha },b_{\beta }> \}
$. Therefore, a branch of the anti-derivation operator prescribed by
the family $\{ (\int_{\gamma ^{\alpha }} \sum_j q(y) d\nu _j(y)): ~
\alpha \in \Lambda \} $ or $\{ (\int_{\gamma ^{\alpha }|_{[c,e]}}
\sum_j q(y) d\nu _j(y)): ~ \alpha \in \Lambda ; ~ [c,e]\subset
<a_{\alpha },b_{\beta }> \} $ is defined up to a function defined on
the boundary $\partial U$ when it is non-void or by convergence to a
definite limit at infinity along paths, when $U$ is unbounded in
certain directions ${\bf R}\eta $ in the Cayley-Dickson algebra
${\cal A}_v$, $\eta \in {\cal A}_v$. \par Clearly, boundary
conditions are necessary for specifying a concrete solution or a
branch of an anti-derivative, since in the definition of the line
integral $\int_{\gamma } f(z)d\nu (z)$ the operator ${\hat f}$ is
restricted to the condition ${\hat f}(z).1=f(z)$ for each $z\in U$
so it is defined up to a function of $2^v-1$ independent real
variables (see also \S 3). In accordance with the formulas of \S \S
4 and 5 the anti-derivation operators are defined up to functions of
$2^v-1$ real variables after a suitable change of variables. For
example, $\sum_{j=0}^n(\partial g(z)/\partial z_j) i_j^*=0$ for
$g(z)= nz_0+z_1i_1^*+...+z_ni_n^*$, or $\sum_{j=0}^n(\partial
q(z)/\partial z_j) i_j^*=0$ on the plane $z_0-z_1-...-z_n=0$ for
$q(z)= z_0^2+z_1^2i_1^*+...+z_n^2i_n^*$. These functions can be
written in the $z$-representation due to Formulas 1$(1-3)$.
\par For concrete domains some concrete boundary conditions can be
chosen (see also below). Mention, that a minimal necessary correct
boundary conditions may be not on the entire boundary, but on its
part. Otherwise, they may be on some hyper-surface $S$ in $U$ of
real dimension $2^v-1$ depending on the domain, for example, for an
infinite cylinder $\sf C$ in both directions along its axis with $S$
being the intersection of $\sf C$ with a hyper-plane perpendicular
to its axis.
\par Mention that the homotopy theorem for domains satisfying
Conditions 1$(D1,D2)$ is accomplished for super-differentiable
functions on $U$ (see \cite{ludoyst,ludfov}), but for a continuous
function $f$ on $U$ it may certainly be not true. This is caused by
several reasons. If a family of locally analytic functions $f^n$
converges to $f$ uniformly on a compact sub-domain $V$ in $U$ a
radius $r^n_x$ of local convergence of a power series of $f^n$ in a
neighborhood of a point $x\in V$ may tend to zero with $n$ tending
to the infinity. Phrases $\mu ^n$ in the $z$-representation
corresponding to $f^n$ may be inconsistent on the intersection
$V_x\cap V_y$ of open neighborhoods $V_x$ and $V_y$ of different
points $x, y \in V$, when $V_x\cap V_y\ne \emptyset $ . Functions
$f^n$ or their phrases $\mu ^n$ may be with branching points in the
domain $U$. That is functions $f^n$ accomplishing the approximation
of $f$ may have several branches on $U$ and a slit of $U$ by a
$2^v-1$ dimensional sub-manifold $S^n$ over $\bf R$ may be necessary
to specify branches of $f^n$. But the family $S^n$ with different
$n$ may be inconsistent and $S^n$ may depend of $n$. \par For
super-differentiable functions $f^n$ operator valued functions
${\hat f}^n$ are also super-differentiable. If $f$ is only
continuous non super-differentiable function on the domain $U$, then
the operator valued function ${\hat f}$ is defined only in the sense
of distributions $[{\hat f},\gamma ;\nu ) = \int_{\gamma } f(z)d\nu
(z)$ for any rectifiable path $\gamma $ in $U$ and each function
$\nu $ of bounded variation on $U$.
 Moreover, the homotopy theorem may be non true for generalized functions (see
below).

\par {\bf 7. Particular case.}
We consider a phrase $\nu $ which can be presented as \par
$(P3)$ $\nu = \rho (\mu )$ with a right ${\cal A}_v$-linear
(super)-differentiable phrase $\mu $ and a projection operator $\rho
$ being an $\bf R$-linear combination of the projection operators
$\pi _j$. Particularly, $\rho $ may be the identity operator or one
of the $\pi _j$.
\par For any $z$-differentiable phrase $\psi $ and constants $a, b\in {\cal A}_v$
we have $\int_{\gamma } a(\psi (z)b)dz = a((\int_{\gamma }\psi
(z)dz)b)$ and $\int_{\gamma } (a\psi (z))b dz = (a(\int_{\gamma
}\psi (z)dz))b$. Then in view of the homotopy theorem
\cite{ludoyst,ludfov} Equation 3$(2)$ implies for any such $\nu
=\rho (\mu )$ that
$$(1)\quad \int_{\gamma } \Upsilon (\nu (z)) dz = \rho (\int_{\gamma }
[d\mu (z)/dz]. \{ \sum_j [(dz/dz).i_j](i_j^* {\psi }_j(z) dz) \} )
=$$
$$ \rho (\int_{\gamma }[d\mu (z)/dz]. \sum_{j} \{ i_j(i_j^* \psi
_j(z))dz \} ) = \rho (\int_{\gamma } [d\mu (z)/dz].[a(z) dz])$$ $$ =
\rho (\mu (z_a(\beta ))) - \rho (\mu (z_a(\alpha  ))) = \nu
(z_a(\beta )) - \nu (z_a(\alpha )),$$ since each ${\psi }_j(z)$ is
the ${\cal A}_v$-valued function, where
\par $(2)$ $z_a(x) = \int_{\alpha }^x  a(t)dt +\phi _a(x',\alpha ') $,
\par $(3)$ $a(z) := \sum_{j=0}^{2^v-1}{\psi }_j(z) ~ $, $\gamma
(0)=\alpha $, $\gamma (1)=\beta $. In particular, if each function
${\psi }_j$ is identically constant, then
\par $(4)$ $\int_{\gamma } \sum_j [(dz/dz).i_j] [(i_j^* {\psi }_j(z) dz] =
t\beta - t\alpha  -t \phi _1 (\beta ', \alpha ')$, \\ where
$t=\sum_j{\psi }_j$.
\par For non right ${\cal A}_v$-linear $z$-differentiable phrase
$\mu $ Formulas $(1-3)$ may already be not valid. Certainly common
line integrals of $z$-differentiable phrases (functions) can be
calculated by the general algorithms (see
\cite{ludoyst,ludfov,lutsltjms,luviralgoclanl}). A result of the
line integration along a rectifiable path $\gamma $ in the domain
$U$ we denote as the composition of two functions
$$(5)\quad \sum_j \int_{\gamma } [(d\mu (z)/dz).i_j] [
i_j^* \psi _j(z)]dz_0 = \int_{\gamma } (d\nu (\xi )/d\xi ).d\xi $$
$$ = \lambda (\xi (\beta )) - \lambda (\xi (\alpha )),$$ where $\lambda $ and
$\xi $ are two $z$-differentiable functions on their domains,
$\gamma (0)=\alpha $, $\gamma (1)=\beta $. Frequently one can use a
Cayley-Dickson subalgebra $\cal G$ isomorphic with either the
quaternion skew field $\bf H$ or the octonion algebra $\bf O$ so
that $\gamma (1)-\gamma (0)\in {\cal G}$ and use the homotopy
theorem.  On the other hand, each rectifiable continuous path
$\gamma $ in the domain $U$ in the Cayley-Dickson algebra ${\cal
A}_v$ can be presented as a uniform limit of rectifiable continuous
paths $\gamma ^n$ in $U$ composed of segments parallel to axes ${\bf
R}i_k$, $k=0,...,2^v-1$. Therefore,
$$\int_{\gamma } f(z)dz = \lim_{n\to \infty } \int_{\gamma ^n} f(z)
dz$$ for any continuous function on $U$ (see \cite{ludoyst,ludfov}).
The functions $\lambda $ and $\xi $ depend on $\psi $ so in more
details we denote them by $\lambda = \lambda _{\psi }$ and $\xi =
\xi _{\psi }$.
\par Thus the general integral of Equation 2$(12)$ is:
$$(6)\quad  \lambda _{\psi }(\xi _{\psi }(x)) = - \phi _{\lambda
'} (Im ~ \xi _{\psi }(x)) + \int_{\alpha }^x g_1(z)dz + \phi
_{g_1}(x'),$$ where $Im (z) := z- Re (z)$, $Re (z) := (z+z^*)/2$.
The term $\phi _{\lambda '}(Im (\xi ))$ takes into account the
non-commutativity for $2\le v$ and non-associativity for $3\le v$ of
the Cayley-Dickson algebra ${\cal A}_v$, since its center is the
real field ${\bf R} = Z({\cal A}_v)$ for any $v\ge 2$. There is the
bijective correspondence between $\lambda _{\psi }(\xi _{\psi })$
and $f$ which will be specified below.
\par {\bf 8. Transformation of the first order partial differential
operator over the Cayley-Dickson algebras.}
\par To simplify the operator $\Upsilon $ and its particular variant $\sigma $
one can use a change of variables. We consider this operator in the
form:
\par $(1)$ $\Upsilon f= \sum_{j=0}^{2^v-1}  (\partial f/\partial
z_j) {\eta }_j(z)$, \\ with either  $\eta _j(z) = i_j^*\psi _j(z)$
or $\eta _j(z)=\phi _j^*(z)\in {\cal A}_v$ for each $j$ (see
Theorems 4 and 5 above). For it we seek the change of variables
$x=x(z)$ so that
\par $(2)$ $\sum_{j=0}^{2^v-1} (\partial x_l/\partial z_j) \omega _j(z) = t_l$, \\ where
$t_l\in {\cal A}_v$ is a constant for each $l$, for $\eta _j$ not
being identically zero, while $\omega _j$ is chosen arbitrarily also
$z$-differentiable so that the resulting matrix $\Omega $ will not
be degenerate, i.e. its rows are real-independent as vectors (see
below). Certainly $(\partial x_l/\partial z_j)\in {\bf R}$ are real
partial derivatives, since $x_l$ and $z_j$ are real coordinates. We
suppose that the functions $\eta _j(z)$ are linearly independent
over the real field for each $z$ in the domain $U$. Using the
standard basis of generators $\{ i_j: ~ j=0,...,2^v -1 \} $ of the
Cayley-Dickson algebra ${\cal A}_v$ and the decompositions $\omega
_j = \sum_k \omega _{j,k} i_k$ and $t_j = \sum_k t_{j,k} i_k$ with
real elements $\omega _{j,k}$ and $t_{j,k}$ for all $j$ and $k$ we
rewrite System $(2)$ in the matrix form:
\par $(3)$ $(\partial x_l/\partial z_j)_{l, j =0,...,2^v-1} \Omega = T$, \\ where $\Omega
= (\omega _{j,k})_{j, k =0,...,2^v-1} $, $T=(t_{j,k})_{j, k
=0,...,2^v-1} $. Suppose that the functions $\omega _j(z)$ are
arranged into the family $\{ \omega _j: j=0,...,2^v-1 \} $ as above
and are such that the matrix $\Omega (z)$ is non-degenerate for all
$z$ in the domain $U$. For example, this is always the case, when
$|\omega _j(z)|>0$ and $Re [\omega _j(z) \omega _k(z)^*] =0$ for
each $j\ne k$ for each $z\in U$. Here particularly $\omega
_j(z)=\eta _j(z)$ can also be taken for all $j=0,...,2^v-1$ and
$z\in U$. Therefore, Equality $(3)$ becomes equivalent to
\par $(4)$ $(\partial x_l/\partial z_j)_{l, j =0,...,2^v-1} = T \Omega
^{-1}$. \\ We take the real matrix $T$ of the same rank as the real
matrix $(\omega _{j,k})_{j,k =0,...,2^v-1}$. Thus $(4)$ is the
linear system of partial differential equations of the first order
over the real field. It can be solved by the standard methods
\cite{matveevb}.
\par We remind how each linear partial differential equation $(3)$ or
$(4)$ can be resolved. Write it in the form:
\par $(5)$ $X_1(x_1,...,x_n,u) \partial u/\partial x_1
+...+X_n(x_1,...,x_n,u) \partial u /\partial x_n = R(x_1,...,x_n,u)$
\\ with $u$ and $x_1,...,x_n$ here instead of $x_l$ and $z_0,...,z_{2^v-1}$
in $(3)$ seeking simultaneously suitable $R$ corresponding to
$t_{l,k}$. A function $u=u(x_1,...,x_n)$ defined and continuous with
its partial derivatives $\partial u/\partial x_1$, ...,$\partial
u/\partial x_n$ in some domain $V$ of variables $x_1,...,x_n$ in
${\bf R}^n$ making $(5)$ the identity is called a solution of this
linear equation. If $R=0$ identically, then the equation is called
homogeneous. A solution $u=const $ of the homogeneous equation
\par $(6)$ $X_1(x_1,...,x_n,u) \partial u/\partial x_1
+...+X_n(x_1,...,x_n,u) \partial u /\partial x_n = 0$ \\ is called
trivial. Then one composes the equations:
\par $(7)$ $dx_1/X_1(x)=dx_2/X_2(x)=...=dx_n/X_n(x)$, \\ where
$x=(x_1,...,x_n)$. This system is called the system of ordinary
differential equations in the symmetric form corresponding to the
homogeneous linear equation in partial derivatives. It is supposed
that the coefficients $X_1,...,X_n$ are defined and continuous
together with their first order partial derivatives by $x_1,...,x_n$
and that $X_1,...,X_n$ are not simultaneously zero in a neighborhood
of some point $x^0$. Such point $x^0$ is called non singular. For
example when the function $X_n$ is non-zero System $(7)$ can be
written as:
\par $(8)$ $dx_1/dx_n = X_1/X_n$,...,$dx_{n-1}/dx_n = X_{n-1}/X_n$.
\\ This system satisfies conditions of the theorem about an
existence  of integrals of the normal system. A system of $n$
differential equations
\par $(9)$ $dy_k/dx=f_k(x,y_1,...,y_n)$, $k=1,...,n$, \\
is called normal of the $n$-th order. It is called linear if all
functions $f_k$ depend linearly on $y_1,...,y_n$. Any family of
functions $y_1,...,y_n$ satisfying $(9)$ in some interval $(a,b)$ is
called its solution. A function $g(x,y_1,...,y_n)$ different from a
constant identically and differentiable in a domain $D$ and such
that its partial derivatives $\partial g/\partial y_1$,...,$\partial
/\partial y_n$ are not simultaneously zero in $D$ is called an
integral of System $(9)$ in $D$ if the complete differential $dg =
(\partial g/\partial x)dx + (\partial g/\partial y_1)dy_1 +...+
(\partial g/\partial y_n)dy_n$ becomes identically zero, when the
differentials $dy_k$ are substituted on their values from $(9)$,
that is $(\partial g(x,y)/\partial x) + (\partial g/\partial
y_1)f_1(x,y) +...+ (\partial g(x,y)/\partial y_n)f_n(x,y)=0$ for
each $(x,y)\in D$, where $y=(y_1,...,y_n)$. The equality
$g(x,y)=const $ is called the first integral of System $(9)$.
\par It is supposed that each function $f_k(x,y)$ is continuous on
$D$ and satisfies the Lipschitz conditions by variables
$y_1,...,y_n$: \par $(L)$ $|f_k(x,y)-f_k(x,z)|\le C_k |y-z|$ \\ for
all $(x,y)$ and $(x,z)\in D$, where $C_k$ are constants. Then System
$(9)$ has exactly $n$ independent integrals in some neighborhood
$D^0$ of a marked point $(x^0,y^0)$ in $D$, when the Jacobian
$\partial (g_1,...,g_n)/\partial (y_1,...,y_n)$ is not zero on $D^0$
(see Section 5.3.3 \cite{matveevb}). \par In accordance with Theorem
12.1,2 \cite{matveevb} each integral of System $(7)$ is a
non-trivial solution of Equation $(6)$ and vice versa each
non-trivial solution of Equation $(6)$ is an integral of $(7)$. If
$g_1(x_1,...,x_n)$,...,$g_{n-1}(x_1,...,x_n)$ are independent
integrals of $(7)$, then the function \par $(10)$ $u=\Phi
(g_1,...,g_{n-1})$, \\ where $\Phi $ is an arbitrary function
continuously differentiable by $g_1,...,g_{n-1}$, is the solution of
$(6)$. Formula $(10)$ is called the general solution of Equation
$(6)$. \par To the non-homogeneous Equation $(5)$ the system
\par $(11)$ $dx_1/X_1=...=dx_n/X_n=du/R$ \\ is posed. System $(11)$
gives $n$ independent integrals $g_1,...,g_n$ and the general
solution \par $(12)$ $\Phi
(g_1(x_1,...,x_n,u),...,g_n(x_1,...,x_n,u))=0$ \\ of $(5)$, where
$\Phi $ is any continuously differentiable function by
$g_1,...,g_n$. If the latter equation is possible to resolve
relative to $u$ this gives the solution of $(5)$ in the explicit
form $u=u(x_1,...,x_n)$ which generally depends on $\Phi $ and
$g_1,...,g_n$. Therefore, Formula $(12)$ for different $R$ and $u$
and $X_j$ corresponding to $t_{l,k}$ and $x_l$ and $\omega _{j,k}$
respectively can be satisfied in $(3)$ or $(4)$, the variables $x_j$
are used in $(12)$ instead of $z_j$ in $(3,4)$, where
$k=0,...,2^v-1$.
\par  Thus after the change
of the variables the operator $\Upsilon $ takes the form:
\par $(13)$ $\Upsilon f= \sum_{j=0}^{2^v-1}  (\partial f/\partial
x_j) t_j$ \\ with constants $t_j\in {\cal A}_v$. Undoubtedly, also
the operator $\Upsilon $ with $j=0,...,n$, $2^{v-1}\le n \le 2^v-1$
instead of $2^v-1$ can also be reduced to the form $\Upsilon f =
\sum_{j=0}^n (\partial f/\partial x_j) t_j$, when the rank is $rank
(\omega _{j,k})= n+1$ in a basis of generators $N_0,...,N_n$, where
$N_0$,...,$N_{2^v-1}$ is a generator basis of the Cayley-Dickson
algebra ${\cal A}_v$.
 Particularly, if
the rank is $rank (\omega _{j,k})= m\le 2^v$ and $T$ has the unit
upper left $m\times m$ block and zeros outside it, then $t_j=N_j$
for each $j=0,...,m-1$ can be chosen.
\par One can mention that direct algorithms of
Theorems 4 and 5 may be simpler for finding the anti-derivative
operator ${\cal I}_{\Upsilon }$, than this preliminary
transformation of the partial differential operator $\Upsilon $ to
the standard form $(13)$.

\par {\bf 9. Definitions.}
\par Let $X$ and $Y$ be two $\bf R$ linear normed spaces which are
also left and right ${\cal A}_r$ modules, where $1\le r$. Let $Y$ be
complete relative to its norm.  We put $X^{\otimes k} := X\otimes
_{\bf R} ... \otimes _{\bf R} X$ is the $k$ times ordered tensor
product over $\bf R$ of $X$. By $L_{q,k}(X^{\otimes k},Y)$ we denote
a family of all continuous $k$ times $\bf R$ poly-linear and ${\cal
A}_r$ additive operators from $X^{\otimes k}$ into $Y$. Then
$L_{q,k}(X^{\otimes k},Y)$ is also a normed $\bf R$ linear and left
and right ${\cal A}_r$ module complete relative to its norm. In
particular, $L_{q,1}(X,Y)$ is denoted also by $L_q(X,Y)$.
\par We present $X$ as the direct sum $X=X_0i_0\oplus ... \oplus
X_{2^r-1} i_{2^r-1}$, where $X_0$,...,$X_{2^r-1}$ are pairwise
isomorphic real normed spaces. If $A\in L_q(X,Y)$ and $A(xb)=(Ax)b$
or $A(bx)=b(Ax)$ for each $x\in X_0$ and $b\in {\cal A}_r$, then an
operator $A$ we call right or left ${\cal A}_r$-linear respectively.
\par An $\bf R$ linear space of left (or right) $k$ times ${\cal A}_r$
poly-linear operators is denoted by $L_{l,k}(X^{\otimes k},Y)$ (or
$L_{r,k}(X^{\otimes k},Y)$ respectively).
\par As usually a support of a function $g: S\to {\cal A}_r$
on a topological space $S$ is by the definition $supp (g) = cl \{
t\in S: ~ g(t)\ne 0 \} $, where the closure is taken in $S$.
\par We consider a space of test function ${\cal D} := {\cal D}({\bf
R}^n,Y)$ consisting of all infinite differentiable functions $f:
{\bf R}^n\to Y$ on ${\bf R}^n$ with compact supports. A sequence of
functions $f_n\in {\cal D}$ tends to zero, if all $f_n$ are zero
outside some compact subset $K$ in the Euclidean space ${\bf R}^n$,
while on it for each $k=0,1,2,...$ the sequence $ \{ f^{(k)}_n: ~
n\in {\bf N} \} $ converges to zero uniformly. Here as usually
$f^{(k)}(t)$ denotes the $k$-th derivative of $f$, which is a $k$
times $\bf R$ poly-linear symmetric operator from $({\bf
R}^n)^{\otimes k}$ to $Y$, that is $f^{(k)}(t).(h_1,...,h_k)=
f^{(k)}(t).(h_{\sigma (1)},...,h_{\sigma (k)})\in Y$ for each
$h_1,...,h_k\in {\bf R}^n$ and every transposition $\sigma : \{
1,...,k \} \to \{ 1,...,k \}$, $\sigma $ is an element of the
symmetric group $S_k$, $t\in {\bf R}^n$. For convenience one puts
$f^{(0)}=f$. In particular, $f^{(k)}(t).(e_{j_1},...,e_{j_k})=
\partial ^kf(t)/\partial t_{j_1}...\partial t_{j_k}$ for all
$1\le j_1,...,j_k\le n$, where $e_j = (0,...,0,1,0,...,0)\in {\bf
R}^n$ with $1$ on the $j$-th place.
\par Such convergence in $\cal D$ defines closed subsets in this space $\cal D$, their
complements by the definition are open, that gives the topology on
$\cal D$. The space ${\cal D}$ is $\bf R$ linear and right and left
${\cal A}_r$ module.

\par By a generalized function of class ${\cal D}' := [{\cal D}({\bf R}^n,Y)]'$ is called a continuous
$\bf R$-linear ${\cal A}_r$-additive function $g: {\cal D} \to {\cal
A}_r$. The set of all such functionals is denoted by ${\cal D}'$.
That is, $g$ is continuous, if for each sequence $f_n\in \cal D$,
converging to zero, a sequence of numbers $g(f_n)=: [g,f_n) \in
{\cal A}_r$ converges to zero for $n$ tending to the infinity. \par
A generalized function $g$ is zero on an open subset $V$ in ${\bf
R}^n$, if $[g,f)=0$ for each $f\in {\cal D}$ equal to zero outside
$V$. By a support of a generalized function $g$ is called the
family, denoted by $supp (g)$, of all points $t\in {\bf R}^n$ such
that in each neighborhood of each point $t\in supp (g)$ the
functional $g$ is different from zero. The addition of generalized
functions $g, h$ is given by the formula: \par $(1)$ $[g+h,f):=
[g,f)+ [h,f)$.
\par The multiplication $g\in {\cal D}'$ on an infinite differentiable
function $w$ is given by the equality:
\par $(2)$ $[gw,f)=[g, wf)$ either for $w: {\bf R}^n\to {\cal A}_r$ and
each test function $f\in \cal D$ with a real image $f({\bf
R}^n)\subset {\bf R}$, where $\bf R$ is embedded into $Y$; or $w:
{\bf R}^n\to {\bf R}$ and $f: {\bf R}^n\to Y$. \par  A generalized
function $g'$ prescribed by the equation: \par $(3)$ $[g',f):= -
[g,f')$ is called a derivative $g'$ of a generalized function $g$,
where $f' \in {\cal D}({\bf R}^n,L_q({\bf R}^n,Y))$, $g'\in [{\cal
D}({\bf R}^n,L_q({\bf R}^n,Y))]'$.
\par Another space ${\cal B} := {\cal B}({\bf R}^n,Y)$ of test functions consists of all
infinite differentiable functions  $f: {\bf R}^n\to Y$ such that the
limit $\lim_{|t|\to +\infty } |t|^m f^{(j)}(t)=0$ exists for each
$m=0,1,2,...$, $j=0,1,2,...$. A sequence $f_n\in \cal B$ is called
converging to zero, if the sequence $|t|^mf_n^{(j)}(t)$ converges to
zero uniformly on ${\bf R}^n\setminus B({\bf R}^n,0,R)$ for each $m,
j=0,1,2,...$ and each $0<R< + \infty $, where $B(Z,z,R) := \{ y\in
Z: ~ \rho (y,z)\le R \} $ denotes a ball with center at $z$ of
radius $R$ in a metric space $Z$ with a metric $\rho $. The family
of all $\bf R$-linear and ${\cal A}_r$-additive functionals on $\cal
B$ is denoted by ${\cal B}'$.
\par In particular we can take $X={\cal A}_r^{\alpha }$, $Y= {\cal
A}_r^{\beta }$ with $1\le \alpha , \beta \in \bf Z$. Analogously
spaces ${\cal D}(U,Y)$, $[{\cal D}(U,Y)]'$, ${\cal B}(U,Y)$ and
$[{\cal B}(U,Y)]'$ are defined for domains $U$ in ${\bf R}^n$. For
definiteness we write ${\cal B}(U,Y) = \{ f|_U: ~ f\in {\cal B}({\bf
R}^n,Y) \} $ and ${\cal D}(U,Y) = \{ f|_U: ~ f\in {\cal D}({\bf
R}^n,Y) \} $.
\par A function $g: U\to {\cal A}_v$ is called locally integrable, if it is
absolutely integrable on each bounded $\lambda $ measurable
sub-domain $V$ in $U$, i.e. \par $\int_V |g(z)|\lambda (dz)<\infty
$, where $\lambda $ denotes the Lebesgue measure on $U$.
\par A generalized function $f$ is called regular if locally integrable
functions $\mbox{}_{j,k}f^1, \mbox{}_lf^2: U\to {\cal A}_v$ exist
such that
\par $[f,\omega ) = \int_U \{ \sum_{j,k,l} \mbox{}_{j,k}f^1(z)
\mbox{}_k\omega (z) \mbox{}_jf^2(z) \}_{q(3)} \lambda (dz)$, \\ for
each test function $\omega \in {\cal B}(U,Y)$ or $\omega \in {\cal
D}(U,Y)$ correspondingly, where $ ~ \omega = (\mbox{}_1\omega
,...,\mbox{}_{\beta }\omega )$, $ ~ \mbox{}_k\omega (z)\in {\cal
A}_v$ for each $z\in U$ and all $k$, $~ q(3)$ is a vector indicating
on an order of the multiplication in the curled brackets and it may
depend on the indices $j, l=1,...,\alpha $, $~ k=1,...,\beta $.
\par We supply the space ${\cal B}({\bf R}^n,Y)$
with the countable family of semi-norms
\par $(4)$ $p_{\alpha , k} (f) := \sup_{x\in {\bf R}^n} |(1+|x|)^k \partial ^{\alpha
}f(x)|$ \\ inducing its topology, where $k=0, 1, 2,...$; $\alpha =
(\alpha _1,...,\alpha _n)$, $0\le \alpha _j\in {\bf Z}$. On this
space we take the space ${\cal B}'({\bf R}^n,Y)_l$ of all $Y$ valued
continuous generalized functions (functionals) of the form
\par $(5)$ $f=f_0i_0+...+f_{2^v-1}i_{2^v-1}$ and
$g=g_0i_0+...+g_{2^v-1}i_{2^v-1}$, where $f_j$ and $g_j\in {\cal
B}'({\bf R}^n,Y)$,  with restrictions on ${\cal B}({\bf R}^n,{\bf
R})$ being real- or ${\bf C}_{\bf i}={\bf R}\oplus {\bf i}{\bf R}$-
valued generalized functions $f_0,...,f_{2^v-1}$,
$g_0,...,g_{2^v-1}$ respectively. Let $\phi = \phi _0i_0+...+\phi
_{2^v-1}i_{2^v-1}$ with $\phi _0,...,\phi _{2^v-1}\in {\cal B}({\bf
R}^n,{\bf R})$, then
\par $(6)$ $[f,\phi ) = \sum_{k,j=0}^{2^v-1} [f_j,\phi _k) i_ki_j$.
Let their convolution be defined in accordance with the formula:
\par $(7)$ $[f*g,\phi ) = \sum_{j,k=0}^{2^v-1} ([f_j*g_k,\phi )
i_j)i_k$ \\ for each $\phi \in {\cal B}({\bf R}^n,Y)$. Particularly,
\par $(8)$ $(f*g)(x) = f(x-y)*g(y) = f(y)*g(x-y)$ \\ for all $x, y \in {\bf R}^n$
due to $(7)$, since the latter equality is satisfied for each pair
$f_j$ and $g_k$.

\par {\bf 10. The decomposition theorem of partial differential operators
over the Cayley-Dickson algebras.}
\par We consider a partial
differential operator of order $u$: $$(1)\quad Af(x)= \sum_{|\alpha
|\le u} {\bf a}_{\alpha }(x)\partial ^{\alpha } f(x),$$ where
$\partial ^{\alpha } f=\partial ^{|\alpha |}f(x)/\partial
x_0^{\alpha _0}...\partial x_n^{\alpha _n}$, $x=x_0i_0+...x_ni_n$,
$x_j\in {\bf R}$ for each $j$, $1\le n=2^r-1$, $\alpha = (\alpha
_0,...,\alpha _n)$, $|\alpha |=\alpha _0+...+\alpha _n$, $0\le
\alpha _j\in {\bf Z}$. By the definition this means that the
principal symbol
$$(2)\quad A_0 := \sum_{|\alpha |= u} {\bf a}_{\alpha }(x)\partial
^{\alpha }$$ has $\alpha $ so that $|\alpha |=u$ and ${\bf
a}_{\alpha }(x)\in {\cal A}_r$ is not identically zero on a domain
$U$ in ${\cal A}_r$. As usually $C^k(U,{\cal A}_r)$ denotes the
space of $k$ times continuously differentiable functions by all real
variables $x_0,...,x_n$ on $U$ with values in ${\cal A}_r$, while
the $x$-differentiability corresponds to the super-differentiability
by the Cayley-Dickson variable $x$.
\par Speaking about locally constant or locally differentiable
coefficients we shall undermine that a domain $U$ is the union of
sub-domains $U^j$ satisfying conditions 15$(D1,i-vii)$ and $U^j\cap
U^k = \partial U^j\cap \partial U^k$ for each $j\ne k$. All
coefficients ${\bf a}_{\alpha }$ are either constant or
differentiable of the same class on each $Int (U^j)$ with the
continuous extensions on $U^j$. More generally it is up to a $C^u$
or $x$-differentiable diffeomorphism of $U$ respectively.
\par  If an operator $A$ is of the odd order $u=2s-1$, then an operator $E$ of
the even order $u+1=2s$ by variables $(t,x)$ exists so that \par
$(3)$ $Eg(t,x)|_{t=0}=Ag(0,x)$ for any $g\in C^{u+1}([c,d]\times
U,{\cal A}_r)$, where $t\in [c,d]\subset {\bf R}$, $c\le 0<d$, for
example, $Eg(t,x) =
\partial (tAg(t,x))/\partial t$.
\par Therefore, it remains the case of the operator $A$ of the even order $u=2s$.
Take $z=z_0i_0+...+z_{2^v-1}i_{2^v-1}\in {\cal A}_v$, $z_j\in {\bf
R}$. Operators depending on a less set $z_{l_1},...,z_{l_n}$ of
variables can be considered as restrictions of operators by all
variables on spaces of functions constant by variables $z_s$ with
$s\notin \{ l_1,...,l_n \} $.
\par {\bf Theorem.} {\it Let $A=A_u$ be a partial differential operator
of an even order $u=2s$ with locally constant or variable (locally)
$C^{s'}$ or $x$-differentiable on $U$ coefficients ${\bf a}_{\alpha
}(x)\in {\cal A}_r$ such that it has the form \par $(4)$ $Af =
c_{u,1}(B_{u,1}f) +...+ c_{u,k}(B_{u,k}f)$, where each
\par $(5)$ $B_{u,p}=B_{u,p,0}+Q_{u-1,p}$ \\ is a partial differential operator
by variables
$x_{m_{u,1}+...+m_{u,p-1}+1}$,...,$x_{m_{u,1}+...+m_{u,p}}$ and of
the order $u$, $m_{u,0}=0$, $c_{u,k}(x)\in {\cal A}_r$ for each $k$,
its principal part \par $(6)$ $B_{u,p,0}= \sum_{|\alpha |=s} {\bf
a}_{p,2\alpha }(x)\partial ^{2\alpha }$ \\ is elliptic with real
coefficients ${\bf a}_{p,2\alpha }(x)\ge 0$, either $0\le r\le 3$
and $f\in C^u(U,{\cal A}_r)$, or $r\ge 4$ and $f\in C^u(U,{\bf R})$.
Then three partial differential operators $\Upsilon ^s$ and
$\Upsilon _1^s$ and $Q$ of orders $s$ and $p$ with $p\le u-1$ with
locally constant or variable (locally) $C^{s'}$ or
$x$-differentiable correspondingly on $U$ coefficients with values
in ${\cal A}_v$ exist and coefficients of the third operator $Q$ may
be generalized functions, when coefficients of $A$ are discontinuous
locally constant or $C^{s'}$ discontinuous on the entire $U$ or when
$s'<s$, $r\le v$, such that
\par $(7)$ $Af=\Upsilon ^s(\Upsilon _1^sf) +Qf$.}
\par {\bf Proof.}  Certainly we have $ord Q_{u-1,p}\le u-1$,
$ord (A-A_0) \le u-1$. We choose the following operators:
$$(8)\quad \Upsilon ^s f(x) = \sum_{p=1}^k \sum_{|\alpha |\le s, ~
\alpha _q = 0 \forall q<(m_{u,1}+...+m_{u,p-1}+1) \mbox{ and }
q>(m_{u,1}+...+m_{u,p})} (\partial ^{\alpha } f(x)) [w_p^* \psi _{p,
\alpha }]\mbox{ and}$$
$$(9)\quad \Upsilon ^s_1 f(x) =
\sum_{p=1}^k \sum_{|\alpha |\le s, ~ \alpha _q = 0 \forall
q<(m_{u,1}+...+m_{u,p-1}+1) \mbox{ and } q>(m_{u,1}+...+m_{u,p})}
(\partial ^{\alpha } f(x)) [w_p\psi _{p,\alpha }^*],$$ where
$w_p^2=c_{u,p}$ for all $p$ and ${\psi }_{p,\alpha }^2(x)= - {\bf
a}_{p,2\alpha }(x)$ for each $p$ and $x$, $w_p\in {\cal A}_r$,
${\psi }_{p,\alpha }(x)\in {\cal A}_{r,v}$ and ${\psi }_{p,\alpha
}(x)$ is purely imaginary for ${\bf a}_{p,2\alpha }(x)>0$ for all
$\alpha $ and $x$, $Re (w_p Im (\psi _{p,\alpha }))=0$ for all $p$
and $\alpha $, $Im (x) = (x-x^*)/2$, $v>r$. Here ${\cal A}_{r,v} =
{\cal A}_v/{\cal A}_r$ is the real quotient algebra. The algebra
${\cal A}_{r,v}$ has the generators $i_{j2^r}$, $j=0,...,2^{v-r}-1$.
A natural number $v$ so that $2^{v-r} -1\ge \sum_{p=1}^k
\sum_{q=0}^u {{m_p+q-1}\choose q}$ is sufficient, where ${m\choose
q} = m!/(q!(m-q)!)$ denotes the binomial coefficient,
${{m+q-1}\choose q}$ is the number of different solutions of the
equation $\alpha _1+...+\alpha _m =q$ in non-negative integers
$\alpha _j$. We have either $\partial ^{\alpha + \beta }f\in {\cal
A}_r$ for $0\le r\le 3$ or $\partial ^{\alpha + \beta }f\in {\bf R}$
for $r\ge 4$. Therefore, we can take $\psi _{p,\alpha }(x) \in
i_{2^rq}{\bf R}$, where $q=q(p,\alpha )\ge 1$, $ ~ ~ q(p^1,\alpha
^1)\ne q(p,\alpha )$ when $(p,\alpha )\ne (p^1,\alpha ^1)$.
\par Thus Decomposition $(7)$ is valid due to the following. For $b=
\partial ^{\alpha +\beta }f(z)$ and ${\bf l} = i_{2^rp}$ and $w\in
{\cal A}_r$ one has the identities: \par $(10)$ $(b(w{\bf l}))
(w^*{\bf l}) = ((wb){\bf l})(w^*{\bf l}) = - w(wb) = - w^2b$ and
\par $(11)$ $(((b{\bf l})w^*){\bf l})w = (((bw){\bf l}){\bf l})w = - (bw)w =
- bw^2$ in the considered here cases, since ${\cal A}_r$ is
alternative for $r\le 3$ while ${\bf R}$ is the center of the
Cayley-Dickson algebra (see Formulas 2$(M1,M2)$).
\par This decomposition of the operator $A_{2s}$ is generally up to
a partial differential operator of order not greater, than $(2s-1)$:
\par $(12)$ $Qf(x) = \sum_{p=1}^k c_{u,p} Q_{u-1,p} +$
\par $\sum_{|\alpha |\le s, |\beta | \le s; \gamma \le \alpha , \epsilon
\le \beta , |\gamma + \epsilon |>0}[\prod _{j=0}^{2^v-1} {{\alpha
_j}\choose {\gamma _j}} {{\beta _j}\choose {\epsilon _j}} ]
(\partial ^{\alpha +\beta - \gamma - \epsilon } f(x))$ \\ $
[(\partial ^{\gamma } {\eta }_{\alpha }(x)) ((\partial ^{\epsilon }
{\eta }_{\beta }^1(x)]$, \\ where operators $\Upsilon ^s$ and
$\Upsilon ^s_1$ are already written in accordance with the general
form \par $(13)$ $\Upsilon ^sf(x) = \sum_{|\alpha |\le s} (\partial
^{\alpha }f(x)) \eta _{\alpha }(x)$; \par  $(14)$ $\Upsilon ^s_1
f(x) = \sum_{|\beta |\le s} (\partial ^{\beta }f(x)) \eta _{\beta
}^1(x)$. \par The coefficients of $Q$ may be generalized functions,
since they are calculated with the participation of partial
derivatives of the coefficients of the operator $\Upsilon ^s_1$, but
the coefficients of the operators $\Upsilon ^s$ and $\Upsilon ^s_1$
may be discontinuous locally constant or $C^{s'}$ discontinuous on
the entire $U$ or $s'<s$ when for the initial operator $A$ they are
such.
\par When $A$ in $(3)$ is with constant coefficients, then the
coefficients $w_p$ and $\psi _{p,\alpha }$ for $\Upsilon ^m$ and
$\Upsilon ^m_1$ can also be chosen constant and $Q - \sum_{p=1}^k
c_{u,p} Q_{u-1,p} =0$.
\par {\bf 11. Corollary.} {\it Let suppositions of Theorem 10 be
satisfied. Then a change of variables locally affine or variable
$C^1$ or $x$-differentiable on $U$ correspondingly on $U$ exists so
that the principal part $A_{2,0}$ of $A_{2}$ becomes with constant
coefficients, when ${\bf a}_{p,2\alpha }>0$ for each $p$, $\alpha $
and $x$.}
\par {\bf 12. Corollary.} {\it If two operators $E=A_{2s}$ and $A=A_{2s-1}$
are related by Equation 10$(3)$, and $A_{2s}$ is presented in
accordance with Formulas 10$(4,5)$, then three operators $\Upsilon
^s$, $\Upsilon ^{s-1}$ and $Q$ of orders $s$, $s-1$ and $2s-2$ exist
so that
\par $(1)$ $A_{2s-1}=\Upsilon ^s\Upsilon ^{s-1} +Q$.}
\par {\bf Proof.} It remains to verify that $ord (Q)\le 2s-2$ in
the case of $A_{2s-1}$, where $Q= \{ \partial (tA_{2s-1})/ \partial
t - \Upsilon ^s\Upsilon _1^s \} |_{t=0}$. Indeed, the form $\lambda
(E)$ corresponding to $E$ is of degree $2s-1$ by $x$ and each
addendum of degree $2s$ in it is of degree not less than $1$ by $t$,
consequently, the product of forms $\lambda (\Upsilon _s) \lambda
(\Upsilon ^s_1)$ corresponding to $\Upsilon ^s$ and $\Upsilon ^s_1$
is also of degree $2s-1$ by $x$ and each addendum of degree $2s$ in
it is of degree not less than $1$ by $t$. But the principal parts of
$\lambda (E)$ and $\lambda (\Upsilon _s) \lambda (\Upsilon ^s_1)$
coincide identically by variables $(t,x)$, hence $ord ( \{ E -
\Upsilon ^s\Upsilon _1^s \} |_{t=0}) \le 2s-2$. Let $a(t,x)$ and
$h(t,x)$ be coefficients from $\Upsilon ^s_1$ and $\Upsilon ^s$.
Using the identities
\par $a(t,x) \partial _t \partial ^{\gamma } tg(x) = a(t,x)
\partial ^{\gamma } g(x)$ and \par $h(t,x) \partial ^{\beta } \partial _t
[a(t,x) \partial ^{\gamma } g(x)] = h(t,x)
\partial ^{\beta } [(\partial _t a(t,x))
\partial ^{\gamma } g(x)] $ \\ for any functions
$g(x)\in C^{2s-1}$ and $a(t,x)\in C^s$, \par $ord [(h(t,x) \partial
^{\beta }), (a(t,x) \partial ^{\gamma })]|_{t=0}\le 2s-2$, \\ where
$\partial _t =\partial /\partial t$, $|\beta |\le s-1$, $|\gamma
|\le s$, $[A,B] := AB-BA$ denotes the commutator of two operators,
we reduce $(\Upsilon ^s\Upsilon ^s_1 + Q_1)|_{t=0}$ from Formula
10$(7)$ to the form prescribes by equation $(1)$.
\par \par {\bf 13.} We consider operators of the form:
\par $(1)$ $(\Upsilon ^k + \beta I_r) f(z) :=
\{ \sum_{0<|\alpha |\le k} (\partial ^{\alpha } f(z) {\eta }_{\alpha
}(z) \} + f(z)\beta (z)$, \\ with $\eta _{\alpha }(z)\in {\cal
A}_v$, $\alpha = (\alpha _0,...,\alpha _{2^r-1})$, $0\le \alpha
_j\in {\bf N}$ for each $j$, $|\alpha |=\alpha _0+...+\alpha
_{2^r-1}$, $\beta I_r f(z) := f(z) \beta $,
\par $\partial ^{\alpha }f(z) := \partial ^{|\alpha |} f(z)/\partial
z_0^{\alpha _0}...\partial z_{2^r-1}^{\alpha _{2^r-1}}$, $2\le r \le
v<\infty $, $\beta (z)\in {\cal A}_v$, $z_0,...,z_{2^r-1}\in {\bf
R}$, $z=z_0i_0+...+z_{2^r-1} i_{2^r-1}$.
\par {\bf Proposition.} {\it The operator $(\Upsilon ^k+\beta
)^*(\Upsilon ^k+\beta )$ is elliptic on the space $C^{2k}({\bf
R}^{2^r},{\cal A}_v)$, where $(\Upsilon ^k+\beta )^*$ denotes the
adjoint operator (i.e. with adjoint coefficients).}
\par {\bf Proof.} In view
of Formulas $(1)$ and 4$(8)$ the form corresponding to the principal
symbol of the operator $(\Upsilon ^k+\beta )^*(\Upsilon ^k+\beta )$
is with real coefficients, of degree $2k$ and non-negative definite,
consequently, the operator $(\Upsilon ^k+\beta )^*(\Upsilon ^k+\beta
)$ is elliptic.

\par {\bf 14. Example.} Let $\Upsilon ^*$ be the adjoint
operator defined on differentiable ${\cal A}_v$ valued functions $f$
given by the formula: $$(1)\quad (\Upsilon +\beta )^*f=
[\sum_{j=0}^n (\partial f(z)/\partial z_j) {\phi }_j(z)] +f(z)\beta
(z)^*.$$ Thus we can consider the operator
\par $(2)$ $\Xi _{\beta }:= (\Upsilon + \beta )(\Upsilon + \beta
)^*$. \\ From Proposition $(13)$ we have that the operator $\Xi
_{\beta }$ is elliptic as classified by its principal symbol with
real coefficients. Put $\Xi = \Xi _0$. In the $x$ coordinates from
\S 8 it has the simpler form:
$$(3)\quad (\Upsilon +\beta ) (\Upsilon +\beta )^* f=
\sum_{j=0}^n (\partial ^2f/\partial x_j^2) |t_j|^2$$  $$ + 2
\sum_{0\le j<k \le n} (\partial ^2f/\partial x_j\partial x_k) Re (
t_jt_k^*)+ 2 \sum_{j=0}^n (\partial f/\partial x_j) Re ( t_j^* \beta
) + \{ f |\beta |^2 + \sum_{j=0}^n [f (\partial \beta ^* /\partial
x_j)] t_j \} ,$$ because the coefficients $t_j$ are already
constant. After a change of variables reducing the corresponding
quadratic form to the sum of squares $\sum_j \epsilon _j s_j^2$ we
get the formula:
\par $(4)$ $\Upsilon  \Upsilon ^* f= \sum_{j=1}^m  (\partial ^2f/\partial
s_j^2) \epsilon _j$, \\ where  $s_j\in {\bf R}$, $\epsilon _j= 1$
for $1\le j\le p$ and $\epsilon _j = - 1$ for each $p < j\le m$,
$m\le 2^v$, $1\le p\le m$ depending on the signature $(p,m-p)$.
\par Generally (see Formula 5$(1)$) we have
\par $(5)$ $A= (\Upsilon + \beta ) (\Upsilon _1+ \beta ^1) f(z) =
B_0f(z) +Qf(z)$, where \\ $(6)$ $B_0f(z) = \sum_{j,k} [(\partial
^2f(z)/\partial z_j\partial z_k)\phi ^1_j(z)^*]\phi _k^*(z)
+[f(z)\beta ^1(z)]\beta (z)$ and \par $(7)$ $Qf(z) = \sum_{j,k}
[(\partial f(z)/\partial z_j)(\partial \phi ^1_j(z)^*/\partial
z_k)]\phi _k^*(z)
 +\sum_j [(\partial f(z)/\partial z_j) \phi ^1_j(z)^*]\beta (z)
 + \sum_k [f(z)(\partial \beta ^1(z)/\partial
z_k)]\phi _k^* (z)$,
 \par $(8)$ $(\Upsilon _1+\beta ^1)f(z)=
[\sum_j (\partial f(z)/\partial z_j)\phi ^1_j(z)^*] + f(z)\beta
^1(z)$.
\par The latter equations show that coefficients
of the operator $Q$ may be generalized functions, when $\phi
^1_j(z)$ for some $j$ or $\beta ^1(z)$ are locally $C^0$ or $C^0$ or
locally $C^1$ functions, while $\phi _k(z)$ for each $k$ and $\beta
(z)$ are locally $C^0$ or $C^0$ functions on $U$. We consider this
in more details in the next section.
\par {\bf 15. Partial differential operators with generalized
coefficients.}
\par Let an operator $Q$ be given by Formula 14$(7)$ on a domain $U$.
Initially it is considered as a domain in the Cayley-Dickson algebra
${\cal A}_v$. But in the case when $Q$ and $f$ depend on smaller
number of real coordinates $z_0,...,z_{n-1}$ we can take the real
shadow of $U$ and its sub-domain $V$ of variables
$(z_0,...,z_{n-1})$, where $z_k$ are marked for example being zero
for all $n\le k\le 2^v-1$. Thus we take a domain $V$ which is a
canonical closed subset in the Euclidean space ${\bf R}^n$,
$2^{v-1}\le n\le 2^v-1$, $v\ge 2$. \par ~ A canonical closed subset
$P$ of the Euclidean space $X={\bf R}^n$ is called a quadrant if it
can be given by the condition $P := \{ x\in X: q_j(x)\ge 0 \} $,
where $(q_j: j\in \Lambda _P)$ are linearly independent elements of
the topologically adjoint space $X^*$. Here $\Lambda _P\subset \bf
N$ (with $card (\Lambda _P)=k\le n$) and $k$ is called the index of
$P$. If $x\in P$ and exactly $j$ of the $q_i$'s satisfy $q_i(x)=0$
then $x$ is called a corner of index $j$.

\par That is $P$
is affine diffeomorphic with $P^n = \prod_{j=1}^n [{\sf a}_j,b_j]$,
where $-\infty \le {\sf a}_j <b_j \le \infty $, $[{\sf a}_j,b_j] :=
\{ x\in {\bf R}: ~ {\sf a}_j\le x \le b_j \} $ denotes the segment
in $\bf R$. This means that there exists a vector $p\in {\bf R}^n$
and a linear invertible mapping $C$ on ${\bf R}^n$ so that $C(P)-p =
P^n$.  We put $t^{j,1} := (t_1,...,t_j,...,t_n: ~ t_j={\sf a}_j)$,
$t^{j,2} := (t_1,...,t_j,...,t_n: ~ t_j=b_j)$. Consider
$t=(t_1,...,t_n)\in P^n$. \par This permits to define a manifold $M$
with corners. It is a metric separable space modelled on $X={\bf
R}^n$ and is supposed to be of class $C^s$, $1\le s$. Charts on $M$
are denoted $(U_l, u_l, P_l)$, that is, $u_l: U_l\to u_l(U_l)
\subset P_l$ is a $C^s$-diffeomorphism for each $l$, $U_l$ is open
in $M$, $u_l\circ {u_j}^{-1}$ is of $C^s$ class of smoothness from
the domain $u_j(U_l\cap U_j)\ne \emptyset $ onto $u_l(U_l\cap U_j)$,
that is, $u_j\circ u_l^{-1}$ and $u_l\circ u_j^{-1}$ are bijective,
$\bigcup_jU_j=M$.

\par A point $x\in M$ is called a corner of index $j$
if there exists a chart $(U,u,P)$ of $M$ with $x\in U$ and $u(x)$ is
of index $ind_M(x)=j$ in $u(U)\subset P$. A set of all corners of
index $j\ge 1$ is called a border $\partial M$ of $M$, $x$ is called
an inner point of $M$ if $ind_M(x)=0$, so $\partial M=\bigcup_{j\ge
1}\partial ^jM$, where $\partial ^jM:=\{ x\in M:  ind_M(x)=j \} $
(see also \cite{michor}). We consider that
\par $(D1)$ $V$ is a canonical closed subset in the Euclidean space
${\bf R}^n$, that is $V = cl ( Int (V))$, where $Int (V)$ denotes
the interior of $V$ and $cl (V)$ denotes the closure of $V$. \par
Particularly, the entire space ${\bf R}^n$ may also be taken.
\par Let a manifold $W$ be satisfying the following conditions $(i-v)$.
\par $(i)$. The manifold $W$ is continuous and
piecewise $C^{\alpha }$, where $C^l$ denotes the family of $l$ times
continuously differentiable functions. This means by the definition
that $W$ as the manifold is of class $C^0\cap C^{\alpha }_{loc}$.
That is $W$ is of class $C^{\alpha }$ on open subsets $W_{0,j}$ in
$W$ and $W\setminus (\bigcup_j W_{0,j})$ has a codimension not less
than one in $W$.
\par $(ii)$. $W=\bigcup_{j=0}^m W_{j}$, where $W_{0} = \bigcup_k
W_{0,k}$, $W_{j}\cap W_{k} = \emptyset $ for each $k\ne j$, $m =
dim_{\bf R} W$, $dim_{\bf R} W_{j} = m-j$, $W_{j+1}\subset
\partial W_{j}$. \par $(iii)$. Each $W_{j}$ with $j=0,...,m-1$ is an
oriented $C^{\alpha }$-manifold, $W_{j}$ is open in $\bigcup_{k=j}^m
W_{k}$. An orientation of $W_{j+1}$ is consistent with that of
$\partial W_{j}$ for each $j=0,1,...,m-2$. For $j>0$ the set $W_{j}$
is allowed to be void or non-void.
\par $(iv)$. A sequence $W^k$ of $C^{\alpha }$ orientable manifolds
embedded into ${\bf R}^n$, $\alpha \ge 1$, exists such that $W^k$
uniformly converges to $W$ on each compact subset in ${\bf R}^n$
relative to the metric $dist$. \par For two subsets $B$ and $E$ in a
metric space $X$ with a metric $\rho $ we put
\par $(1)\quad dist
(B,E) := \max \{ \sup_{b\in B} dist (\{ b \} ,E), \sup_{e\in E} dist
(B,\{ e \} ) \} ,$ where \par $dist (\{ b \} ,E) := \inf_{e\in E}
\rho (b,e)$, $dist (B, \{ e \} ) := \inf_{b\in B} \rho (b,e)$, $b\in
B$, $e\in E$.
\par Generally, $dim_{\bf R} W=m\le n$. Let
$(e_1^k(x),...,e_m^k(x))$ be a basis in the tangent space $T_xW^k$
at $x\in W^k$ consistent with the orientation of $W^k$, $k\in {\bf
N}$.
\par We suppose that the sequence of orientation frames $(e^k_1(x_k),...,e_m^k(x_k))$
of $W^k$ at $x_k$ converges to $(e_1(x),...,e_m(x))$ for each $x\in
W_0$, where $\lim_kx_k = x\in W_0$, while $e_1(x)$,...,$e_m(x)$ are
linearly independent vectors in ${\bf R}^n$.
\par $(v)$. Let a sequence of Riemann volume elements $\lambda _k$
on $W^k$ (see \S XIII.2 \cite{zorich}) induce a limit volume element
$\lambda $ on $W$, that is, $\lambda (B\cap W) = \lim_{k\to \infty }
(B\cap W^k)$ for each compact canonical closed subset $B$ in ${\bf
R}^n$, consequently, $\lambda (W\setminus W_0)=0$. We shall consider
surface integrals of the second kind, i.e. by the oriented surface
$W$ (see $(iv)$), where each $W_j$, $j=0,...,m-1$ is oriented (see
also \S XIII.2.5 \cite{zorich}).
\par  Suppose that a
boundary $\partial U$ of $U$ satisfies Conditions $(i-v)$ and \par
$(vii)$ let the orientations of $\partial U^k$ and $U^k$ be
consistent for each $k\in {\bf N}$ (see Proposition 2 and Definition
3 \cite{zorich}).
\par Particularly, the Riemann volume element $\lambda _k$ on
$\partial U^k$ is consistent with the Lebesgue measure on $U^k$
induced from ${\bf R}^n$ for each $k$. This induces the measure
$\lambda $ on $\partial U$ as in $(v)$. This consideration certainly
encompasses the case of a domain $U$ with a $C^{\alpha }$ boundary
$\partial U$ as well.
\par Suppose that $U_1$,...,$U_l$ are domains in ${\bf R}^n$
satisfying conditions $(D1,i-vii)$ and such that $U_j\cap
U_k=\partial U_j\cap \partial U_k$ for each $j\ne k$,
$U=\bigcup_{j=1}^lU_j$. Consider a function $g: U\to {\cal A}_v$
such that each its restriction $g|_{U_j}$ is of class $C^s$, but $g$
on the entire domain $U$ may be discontinuous or not $C^k$, where
$0\le k\le s$, $1\le s$. If $x\in \partial U_j\cap \partial U_k$ for
some $j\ne k$, $x\in Int (U)$, such that $x$ is of index $m\ge 1$ in
$U_j$ (and in $U_k$ also). Then there exists canonical $C^{\alpha }$
local coordinates $(y_1,...,y_n)$ in a neighborhood $V_x$ of $x$ in
$U$ such that $S := V_x\cap \partial ^mU_j = \{ y: ~ y\in V_x; ~
y_1=0,...,y_m = 0 \} $. Using locally finite coverings of the
locally compact topological space $\partial U_j\cap \partial U_k$
without loss of generality we suppose that $C^{\alpha }$ functions
$P_1(z),...,P_m(z)$ on ${\bf R}^n$ exist with $S = \{ z: ~ z\in {\bf
R}^n; ~ P_1(z)=0,...,P_m(z) = 0 \} $. Therefore, on the surface $S$
the delta-function $\delta (P_1,...,P_m)$ exists, for $m=1$ denoting
them $P=P_1$ and $\delta (P)$ respectively (see \S III.1
\cite{gelshil}). It is possible to choose $y_j=P_j$ for $j=1,...,m$.
Using generalized functions with definite supports, for example
compact supports, it is possible without loss of generality consider
that $y_1,...,y_n\in {\bf R}$ are real variables. Let $\theta (P_j)$
be the characteristic function of the domain ${\cal P}_j := \{ z:
P_j(z) \ge 0 \} $, $\theta (P_j):= 1$ for $P_j\ge 0$ and $\theta
(P_j)=0$ for $P_j<0$. Then the generalized function $\theta
(P_1,...,P_m) := \theta (P_1)...\theta (P_m)$ can be considered as
the direct product of generalized functions $\theta
(y_1)$,...,$\theta (y_m)$, $1(y_{m+1},...,y_n)\equiv 1$, since
variables $y_1,...,y_n$ are independent. Then in the class of
generalized functions the following formulas are valid:
\par $(2)$ $\partial \theta (P_j)/\partial z_k=\delta (P_j)(\partial
P_j/\partial z_k)$ for each $k=1,...,n$, consequently,
\par $(3)$ $grad [\theta (P_1,...,P_m)]
=\sum_{j=1}^m [\theta (P_1)...\theta (P_{j-1})\delta (P_j)(grad (
P_j))\theta (P_{j+1})...\theta (P_m)]$, \\ where , $grad ~ g(z) :=
(\partial g(z)/\partial z_1,...,\partial g(z)/\partial z_n)$ (see
Formulas III.1.3$(1,7,7',9)$ and III.1.9$(6)$ \cite{gelshil}).
\par Let for the domain $U$ in the Euclidean space ${\bf R^n}$
the set of internal surfaces $cl_U [ Int_{{\bf R}^n} (U) \cap
\bigcup_{j\ne k} (\partial U_j\cap \partial U_k)]$ in $U$ on which a
function $g: U\to {\cal A}_v$ or its derivatives may be
discontinuous is presented as the disjoint union of surfaces $\Gamma
_j$, where each surface $\Gamma ^j$ is the boundary of the
sub-domain \par $(4)$ ${\cal P}^j := \{ P_{j,1}(z)\ge 0,...,
P_{j,m_j}(z)\ge 0 \} $, $\Gamma ^j =\partial {\cal
P}^j=\bigcup_{k=1}^{m_j} \partial ^k {\cal P}^j$, \\ $m_j\in {\bf
N}$ for each $j$, $cl_X (V)$ denotes the closure of a subset $V$ in
a topological space $X$, $Int_X(V)$ denotes the interior of $V$ in
$X$. By its construction $\{ {\cal P}^j: ~ j \} $ is the covering of
$U$ which is the refinement of the covering $\{ U_k: ~ k \} $, i.e.
for each ${\cal P}^j$ a number $k$ exists so that ${\cal P}^j\subset
U_k$ and $\partial {\cal P}^j\subset
\partial U_k$ and $\bigcup_j {\cal P}^j = \bigcup_k U_k=U$. We put
\par $(5)$ $h_j(z(x))=h(x)|_{x\in \Gamma ^j} := $\par $\lim_{y_{j,1}\downarrow
0,...,y_{j,n}\downarrow 0} g(z(x+y)) - \lim_{y_{j,1}\downarrow
0,...,y_{j,n}\downarrow 0} g(z(x-y))$ \\ in accordance with the
supposition made above that $g$ can have only discontinuous of the
first kind, i.e. the latter two limits exist on each $\Gamma ^j$,
where $x$ and $y$ are written in coordinates in ${\cal P}^j$,
$z=z(x)$ denotes the same point in the global coordinates $z$ of the
Euclidean space ${\bf R}^n$. We take new continuous function
\par $(6)$ $g^1(z) = g(z) - \sum_j h_j(z)\theta
(P_{j,1}(z),...,P_{j,m_j}(z))$.
\\ Let the partial derivatives and the gradient of the function
$g^1$ be calculated piecewise one each $U_k$. Since $g^1$ is the
continuous function, it is the regular generalized function by the
definition, consequently, its partial derivatives exist as the
generalized functions. If $g^1|_{U_j}\in C^1(U_j,{\cal A}_v)$, then
$\partial g^1(z)/\partial z_k$ is the continuous function on $U_j$,
i.e. in this case $\partial g^1(z)\chi_{U_j}(z)/\partial z_k$ is the
regular generalized function on $U_j$ for each $k$, where $\chi
_G(z)$ denotes the characteristic function of a subset $G$ in ${\cal
A}_v$, $\chi _G(z)=1$ for each $z\in G$, while $\chi (z)=0$ for
$z\in {\cal A}_v\setminus G$. Therefore, $g^1(z) =\sum_j
g^1(z)\chi_{U_j\setminus \bigcup_{k<j}U_k}(z)$, where $U_0=\emptyset
$, $j, k\in {\bf N}$. \par On the other hand, the function $g(z)$ is
locally continuous on $U$, consequently, it defines the regular
generalized function on the space ${\cal D}(U,{\cal A}_v)$ of test
functions as
\par $[g,\omega ) := \int_U g(z)\omega (z)\lambda (dz)$, \\ where $\lambda $
is the Lebesgue measure on $U$ induced by the Lebesgue measure on
the real shadow ${\bf R}^{2^v}$ of the Cayley-Dickson algebra ${\cal
A}_v$, $ ~ \omega \in {\cal D}(U,{\cal A}_v)$. Thus partial
derivatives of $g$ exist as generalized functions. \par  In
accordance with Formulas $(2,3,6)$ we infer that the gradient of the
function $g(z)$ on the domain $U$ is the following:
\par $(7)$ $grad ~ g(z) = grad ~ g^1(z) + \sum_j h_j(z) grad ~
\theta (P_{j,1},...,P_{j,m_j})$.
\par Thus Formulas $(3,7)$ permit calculations of coefficients
of the partial differential operator $Q$ given by Formula 14$(7)$.

\par {\bf 16. Line generalized functions.}
\par Let $U$ be a domain
satisfying Conditions 1$(D1,D2)$ and 15$(D1,i-vii)$. We embed the
Euclidean space ${\bf R}^n$ into the Cayley-Dickson algebra ${\cal
A}_v$, $2^{v-1}\le n\le 2^v-1$, as the $\bf R$ affine sub-space
putting ${\bf R}^n\ni x=(x_1,...,x_n)\mapsto
x_1i_{j_1}+...+x_ni_{j_n}+x^0\in {\cal A}_v$, where $j_k\ne j_l$ for
each $k\ne l$, $x^0$ is a marked Cayley-Dickson number, for example,
$j_k=k$ for each $k$, $x^0=0$. Moreover, each $z_j$ can be written
in the $z$-representation in accordance with Formulas 1$(1-3)$.
\par We denote by ${\bf P}={\bf P}(U)$ the family of all rectifiable
paths $\gamma : [a_{\gamma },b_{\gamma }]\to U$ supplied with the
metric
\par $(1)$ $\rho (\gamma ,\omega ) := |\gamma (a)-\omega (a_{\omega })|+
\inf_{\phi } V_a^b(\gamma (t))
- \omega (\phi (t))$ \\
where the infimum is taken by all diffeomorphisms $\phi : [a_{\gamma
},b_{\gamma }]\to [a_{\omega },b_{\omega }]$ so that $\phi
(a_{\gamma })=a_{\omega }$, $a=a_{\gamma }<b_{\gamma }=b$ (see \S
3).
\par Let us introduce a continuous mapping
$g : {\cal B}(U,{\cal A}_v) \times {\bf P}(U) \times {\cal
V}(U,{\cal A}_v)\to Y$ such that its values are denoted by
$[g;\omega ,\gamma ;\nu ]$, where $Y$ is a module over the
Cayley-Dickson algebra ${\cal A}_v$, $ ~ \omega \in {\cal B}(U,{\cal
A}_v)$, $ ~ \gamma \in {\bf P}(U)$, ${\cal V}(U,{\cal A}_v)$ denotes
the family of all functions on $U$ with values in the Cayley-Dickson
algebra of bounded variation (see \S 3), $\nu \in {\cal V}(U,{\cal
A}_v)$. For the identity mapping $\nu (z)= id(z)=z$ values of this
functional will be denoted shortly by $[g;\omega ,\gamma ]$. Suppose
that this mapping $g$ satisfies the following properties $(G1-G5)$:
\par $(G1)$ $[g;\omega ,\gamma ;\nu ]$ is bi- $\bf R$ homogeneous and ${\cal
A}_v$ additive by a test function $\omega $ and by a function of
bounded variation $\nu $;
\par $(G2)$ $[g;\omega ,\gamma ;\nu ] = [g;\omega ,\gamma ^1;\nu ] +
[g;\omega ,\gamma ^2;\nu ]$ for each $\gamma , \gamma ^1$ and
$\gamma ^2 \in {\bf P}(U)$ such that $\gamma (t)=\gamma ^1(t)$ for
all $t\in [a_{\gamma ^1},b_{\gamma ^1}]$ and $\gamma (t) = \gamma
^2(t)$ for any $t\in [a_{\gamma ^2},b_{\gamma ^2}]$ and $a_{\gamma
^1}=a_{\gamma }$ and $a_{\gamma ^2}=b_{\gamma ^1}$ and $b_{\gamma
}=b_{\gamma ^2}$.
\par $(G3)$ If a rectifiable curve $\gamma $ does not intersect
a support of a test function $\omega $ or a function of bounded
variation $\nu $, $\gamma ([a,b]\cap (supp (\omega )\cap supp (\nu
))=\emptyset $, then $[g;\omega ,\gamma ;\nu ]=0$, where $supp
(\omega ) := cl \{ z\in U: ~ \omega (z)\ne 0 \} $.
\par Further we put
\par $(G4)$ $[\partial ^{|m|}g(z)/\partial z_0^{m_0}...
\partial z_{2^v-1}^{m_{2^v-1}};\omega ,\gamma ] := (-1)^{|m|}
[g;\partial ^{|m|}\omega (z)/\partial z_0^{m_0}...
\partial z_{2^v-1}^{m_{2^v-1}}, \gamma ]$ \\ for each
$m=(m_0,...,m_{2^v-1})$, $ ~ m_j$ is a non-negative integer $0\le
m_j\in {\bf Z}$ for each $j$, $ ~ |m| := m_0+...+m_{2^v-1}$. In the
case of a super-differentiable function $\omega $ and a generalized
function $g$ we also put
\par $(G5)$ $[(d^kg(z)/dz^k).(h_1,...,h_k);\omega ,\gamma ] :=
(-1)^k [g;(d^k\omega (z)/dz^k).(h_1,...,h_k),\gamma ]$ \\
for any natural number $k\in {\bf N}$ and Cayley-Dickson numbers
$h_1,..,h_k\in {\cal A}_v$. \par Then $g$ is called the $Y$ valued
line generalized function on ${\cal B}(U,{\cal A}_v) \times {\bf
P}(U)\times {\cal V}(U,{\cal A}_v)$. Analogously it can be defined
on ${\cal D}(U,{\cal A}_v) \times {\bf P}(U)\times {\cal V}(U,{\cal
A}_v)$ (see also \S 9). If $Y={\cal A}_v$ we call it simply the line
generalized function, while for $Y=L_q({\cal A}_v^k,{\cal A}_v^l)$
we call it the line generalized operator valued function, $k, l\ge
1$, omitting "on ${\cal B}(U,{\cal A}_v) \times {\bf P}(U)\times
{\cal V}(U,{\cal A}_v)$" or "line" for short, when it is specified.
Their spaces we denote by $L_q({\cal B}(U,{\cal A}_v) \times {\bf
P}(U)\times {\cal V}(U,{\cal A}_v);Y)$.
 Thus if $g$ is a generalized function,
then Formula $(G5)$ defines the operator valued generalized function
$d^kg(z)/dz^k$ with $k\in {\bf N}$ and $l=1$.
\par If $g$ is a continuous function on $U$ (see \S 3),
then the formula
\par $(G6)$ $[g;\omega ,\gamma ;\nu ] =
\int_{\gamma } g(y) \omega (y)d\nu (y)$
\\ defines the generalized function. If ${\hat f}(z)$ is a
continuous $L_q({\cal A}_v,{\cal A}_v)$ valued function on $U$, then
it defines the generalized operator valued function with
$Y=L_q({\cal A}_v,{\cal A}_v)$ such that
\par $(G7)$ $[{\hat f};\omega ,\gamma ;\nu ] = \int_{\gamma } \{ {\hat
f}(z).\omega (z) \} d\nu (z)$. \\ Particularly, for $\nu =id$ we
certainly have $d\nu (z)=dz$.
\par We consider on $L_q({\cal B}(U,{\cal A}_v) \times {\bf
P}(U)\times {\cal V}(U,{\cal A}_v);Y)$ the strong topology:
\par $(G8)$ $\lim_l f^l =f$ means that for each marked test function $\omega \in
{\cal B}(U,{\cal A}_v)$ and rectifiable path $\gamma \in {\bf P}(U)$
and function of bounded variation $\nu \in {\cal V}(U,{\cal A}_v)$
the limit relative to the norm in $Y$ exists
\par $\lim_l [f^l;\omega ,\gamma ;\nu ] = [f;\omega ,\gamma ;\nu ]$.

\par {\bf 17. Line integration of generalized functions.}
\par Let $C^m_{ab}(V,{\cal A}_v)$ denote the $\bf R$ linear
space and right and left ${\cal A}_v$ module of all functions
$\gamma : V\to {\cal A}_v$ such that $\gamma (z)$ and each its
derivative $\partial ^{|k|}g(z)/\partial z_1^{m_1}...
\partial z_{n}^{m_n}$ for $1\le |k|\le m$ is absolutely
continuous on $V$ (see \S \S 3 and 16). This definition means that
$C^{m+1}(V,{\cal A}_v)\subset C^m_{ab}(V,{\cal A}_v)$, where
$C^m(V,{\cal A}_v)$ denotes the family of all $m$ times continuously
differentiable functions on a domain $V$ either open or canonical
closed in ${\bf R}^n$, which may be a a real shadow of $U$ as well.
\par {\bf 17.1. Lemma.} {\it Let  $\gamma \in C^m_{ab}([a,b],{\cal
A}_v) \cap {\bf P}(U)$ and $\omega \in {\cal B}(U,{\cal A}_v)$ and
$\nu \in C^0_{ab}(U,{\cal A}_v)$ for $m=0$ or $\nu =id$ for $m\ge
1$, where $0 \le m\in {\bf Z}$, then a line generalized function
$[g;\omega ,\gamma |_{[a,x]};\nu ]$ is continuous for $m=0$ or of
class $C^m$ by the parameter $x\in [a,b]$ for $m\ge 1$.}
\par {\bf Proof.} For absolutely continuous functions $\gamma (t)$
and $\nu $ (i.e. when $m=0$) the continuity by the parameter $x$
follows from the definition of the line generalized function, since
$\lim_{\Delta x\to 0}\rho (\gamma |_{[a,x]}, \gamma |_{[a,x+\Delta
x]})=0$ and \par $\lim_{\Delta x\to 0}\rho (\nu \circ \gamma
|_{[a,x]}, \nu \circ \gamma |_{[a,x+\Delta x]})=0$.
\par Consider now the case $m\ge 1$ and $\nu =id$.
In view of properties 16$(G1,G2)$ for any $\Delta x\ne 0$ so
that $x\in (a,b] := \{ t\in {\bf R}: ~ a<t\le b \} $ and $x+\Delta
x\in (a,b) := \{ t\in {\bf R}: ~ a<t<b \} $ the difference quotient
satisfies the equalities:
\par $(1)$ $ \{ [g; \omega , \gamma |_{[a,x+\Delta x]}] -
[g; \omega , \gamma |_{[a,x]} ] \}/ \Delta x = [g; \omega /\Delta x,
\gamma \circ \phi |_{[a,x]}] - [g; \omega /\Delta x,
\gamma |_{[a,x]}] $, \\
where $\phi : [a,x]\to [a,x+\Delta x]$ is a diffeomorphism of
$[a,x]$ onto $[a,x+\Delta x]$ with $\phi (a)=a$. Therefore, $\Delta
\omega := \omega (z+\Delta z) - \omega (z)$ for $z=\gamma (t)$ and
$z+\Delta z =\gamma (\phi (t))$, $t\in [a,x]$ in the considered
case.  Using Conditions $(G1,G3)$ one can mention that if $\omega
=\omega ^1$ on an open neighborhood $V$ of $\gamma $ in $U$, then
\par $(2)$ $[g;\omega ,\gamma ] = [g;\omega ^1,\gamma ]$, \\
since $\omega - \omega ^1=0$ on $V$ and $\gamma \cap supp (\omega
-\omega ^1)=0$. \par From Conditions 16$(G1,G4)$ and Formula $(2)$
we deduce that
\par $(3)$ $ \lim_{\Delta x\to 0}
\{ [g; \omega , \gamma |_{[a,x+\Delta x]}] - [g; \omega , \gamma
|_{[a,x]}] \}/ \Delta x $ \par $= \sum_{j=0}^{2^v-1} [g; (\partial
\omega (z)/\partial z_j), (d\gamma _j(t)/dt) \gamma |_{[a,x]}]$, \\
where ${z_j}' = d\gamma _j(t)/dt $ for $ ~ z=\gamma (t)$, $ ~ t\in
[a,b]$, since each partial derivative of the test function $\omega $
is again the test function. From the first part of the proof we get
that $[g; \omega , \gamma |_{[a,x]}]$ is of class $C^1$ by the
parameter $x$, since the product $(d\gamma _j(t)/dt) \gamma (t)$ of
absolutely continuous functions $(d\gamma _j(t)/dt)$ and $\gamma
(t)$ is absolutely continuous for each $j$. Repeating this proof by
induction for $k=1,...,m$ one gets the statement of this lemma for
$\gamma \in C^m_{ab}([a,b],{\cal A}_v) \cap {\bf P}(U)$.

\par {\bf 17.2. Lemma.} {\it If $\gamma $ is a rectifiable path,
then a line generalized function $[g;\omega ,\gamma |_{[a,x]}]$ is
of bounded variation by the parameter $x\in [a,b]$.}
\par {\bf Proof.}  Let $\gamma \in {\bf P}(U)$ be a rectifiable path
in $U$, $\gamma : [a,b]\to U$. We can present $\gamma $ in the form
\par $(1)$ $\gamma (t)=\sum_{j=0}^{2^v-1} \gamma _j(t)i_j$, \\ where each
function $\gamma _j(t)$ is real-valued. Therefore, $\gamma _j(t)$ is
continuous and of bounded variation for each $j$, since $\gamma $ is
such. Thus the function $\omega (\gamma (t))$ is of bounded
variation $V_a^b\omega (\gamma ))<\infty $, since $\omega $ is
infinite differentiable and $\gamma ([a,b])$ is compact.

\par  On the other hand, each function $f: [a,b]\to
{\bf R}$ of bounded variation can be written as the difference
$f=f^1-f^2$ of two monotone non-decreasing functions $f^1$ and $f^2$
of bounded variations: $f^1(t):= V_a^tf$ and $f^2(t)=f^1(t)-f(t)$
for each $t\in [a,b]$ (see \cite{fihteng,kolmfom}). This means that
$f^k = g^k+h^k$, where a function $g^k$ is continuous monotone and
of bounded variation, while $h^k$ is a monotone step function, where
$k=1, 2$. When the function $f$ is continuous one gets $f=g^1-g^2$.
For a monotone non-decreasing function $p$ one has
$V_a^tp=p(t)-p(a)$. \par In view of Property 17$(G1)$ we infer that
\par $(2)$ $[g;\omega ,\gamma |_{[a,x]}] = \sum_{j=0}^{2^v-1}
[g_j;\omega ,\gamma |_{[a,x]}] i_j $, \\
where the function $[g_j;\omega ,\gamma |_{[a,x]}]$ by $x$ is
real-valued for any $\omega \in {\cal B}(U,{\cal A}_v)$ and $\gamma
\in {\bf P}(U)$ for all $j=0,...,2^v-1$.

\par The metric space ${\bf P}({\bar U})$ is complete, where
${\bar U}=cl (U)$. Indeed, let $g^n$ be a sequence of rectifiable
paths in $\bar U$ fundamental relative to the metric $\rho $ given
by Formula 16$(1)$. Using diffeomorphism preserving orientations of
segments we can consider without loss of generality that each path
$g^n$ is defined on the unit segment $[0,1]$, $a=0$, $b=1$. It is
lightly to mention that \par $(3)$ $|g(a)-f(a))| + V_a^b (g-f)\ge
\sup_{t\in [a,b]} |g(t)-f(t)|$ \\ for any two functions of bounded
variation, $f, ~ g : [a,b]\to {\bar U}$. For each $\epsilon
>0$ a natural number $n_0=n_0(\epsilon )$ exists so that $\rho
(g^n,g^m)<\epsilon /2$ for all $n, m \ge n_0$. That is $\phi ^n:
[0,1]\to [0,1]$ diffeomorphisms exist such that \par
$|g^n(a)-g^m(a))| + V_a^b (g^n\circ \phi ^n - g^m\circ \phi ^m)
<\epsilon $ for all $n, m\ge n_0$, since $\phi ^m\circ (\phi
^n)^{-1}$ is also the diffeomorphism preserving the orientation of
the segment. Using induction by $\epsilon =1/l$ with $l\in {\bf N}$
one chooses a sequence of diffeomorphisms $\phi ^n$ such that for
each $l\in {\bf N}$ a natural number $n_0=n_0(l)$ exists such that
\par $|g^n(a)-g^m(a))| + V_a^b (g^n\circ \phi ^n - g^m\circ \phi ^m)
<1/l$ for all $n, m\ge n_0(l)$, consequently,
\par $\sup_{t\in [a,b]} |g^n(\phi ^n(t))-g^m(\phi ^m(t))|<1/l$
for all $n, m\ge n_0(l)$. \\ Thus the sequence $g^n\circ \phi ^n$ is
fundamental in $C^0([a,b],{\bar U})$. The latter metric space is
complete relative to the metric \par $d(f,g) := \sup_{t\in [a,b]}
|f(t) - g(t)|$, \\ since from the completeness of the Cayley-Dickson
algebra ${\cal A}_v$ considered as the normed space over the real
field the completeness of the closed subset $\bar U$ follows (see
also Chapter 8 in \cite{eng}). Therefore, the sequence $g^n\circ
\phi ^n$ converges to a continuous function $f: [a,b]\to {\bar U}$.
On the other hand, $\lim_{m\to \infty } \rho (g^n\circ \phi ^n,
g^m\circ \phi ^m)= \rho (g^n\circ \phi ^n, f) \le 1/l$ for each
$n>n_0(l)$, $l\in {\bf N}$. The function $g^n\circ \phi ^n$ is of
bounded variation, consequently, the function $f$ is also of bounded
variation. That is $f\in {\bf P}({\bar U})$. Thus ${\bf P}({\bar
U})$ is complete.

\par Take any sequence $\gamma ^n$ of $C^2_{ab}([a,b],{\cal A}_v)$ paths
in $U$ converging to $\gamma $ relative to the metric $\rho $ on
${\bf P}({\bar U})$ and the latter metric space is complete as it
was demonstrated above. In view of Formula 17.1$(3)$ and Property
16$(G3)$ the sequence $[g;\omega ,\gamma ^n|_{[a,x]}]$ is
fundamental in ${\bf P}({\bar U})$.  On the other hand, the
generalized function $g$ is continuous on ${\cal B}(U,{\cal
A}_r)\times {\bf P}({\bar U})$, consequently, the sequence
$[g;\omega ,\gamma ^n|_{[a,x]}]$ converges in ${\cal B}(U,{\cal
A}_r)\times {\bf P}({\bar U})$ to $[g;\omega ,\gamma ^n|_{[a,x]}]$
for each $a<x\le b$, hence $[g;\omega ,\gamma
|_{[a,x]}]=\lim_n[g;\omega ,\gamma ^n|_{[a,x]}]$ in ${\bf P}({\bar
U})$. By the conditions of this lemma $[g;\omega ,\gamma
|_{[a,x]}]\in {\bf P}(U)$, since $\gamma ([a,b])\subset U$. Thus the
function $[g;\omega ,\gamma |_{[a,x]}]$ by $x\in [a,b]$ is of
bounded variation:
\par $V_a^b[g;\omega ,\gamma |_{[a,x]}]<\infty $.

\par {\bf 18. Definition.} Let $f$ and $\eta $ be two line generalized
functions on ${\cal B}(U,{\cal A}_v) \times {\bf P}(U)\times {\cal
V}(U,{\cal A}_v)$. We define a line functional with values denoted
by
\par $ [ \int_{\gamma } fd\eta , \omega ^1\otimes \omega ) :=
[{\hat f};\omega ^1,\gamma ;[\eta ;\omega ,\kappa ]]|_{\kappa
=\gamma } = [{\hat f};\omega ^1,* ;[\eta ;\omega ,*]](\gamma )$, \\
where $\gamma \in {\bf P}(U)$ is a rectifiable path in $U$, $\omega
, \omega ^1\in {\cal B}(U,{\cal A}_v)$ are any test functions. The
functional $\int_{\gamma } fd\eta $ is called the non-commutative
line integral over the Cayley-Dickson algebra ${\cal A}_v$ of line
generalized functions $f$ by $\eta $. Quite analogously such
integral is defined for line generalized functions $f$ and $\eta $
on ${\cal D}(U,{\cal A}_v) \times {\bf P}(U)\times {\cal V}(U,{\cal
A}_v)$.
\par {\bf 19. Theorem.} {\it Let $F$ and $\Xi $ be two generalized
functions on $U$, $~ F, ~\Xi \in {\cal B}'(U,{\cal A}_v)$ or $ ~ F,
~ \Xi \in {\cal D}'(U,{\cal A}_v)$, then the the non-commutative
line integral over the Cayley-Dickson algebra ${\cal A}_v$ of line
generalized functions $f$ by $\xi $ exists, where $f$ is induced by
$F$ and $\xi $ by $\Xi $.}
\par {\bf Proof.} At first it easy to mention that Definition 18
is justified by Definition 16 and Lemma 17.2, since the function
$[\eta ;\omega ,\kappa |_{[a,x]}]$ is of bounded variation by the
variable $x$ for each rectifiable path $\kappa \in {\bf P}(U)$ and
any test function $\omega $ (see Properties 16$(G1-G3)$), while the
operator ${\hat f}$ always exists in the class of generalized line
operators, ${\hat f}=dg/dz$, $(dg(z)/dz).1=f(z)$ (see Property
16$(G5)$).
\par Each generalized function $f\in {\cal B}(U,{\cal A}_v)$ can be
written in the form:
\par $(1)$ $[f,\omega ) = \sum_{j,k=0}^{2^v-1} [f_{j,k},\omega _k)i_j$,
\\ where each $f_{j,k}$ is a real valued generalized function,
$f_{j,k}\in {\cal B}'(U,{\bf R})$, $\omega = \sum_k \omega _ki_k$,
$\omega _k\in {\cal B}(U,{\bf R})$ is a real valued test function,
$[f_{j,k},\omega _k) = [f_j,\omega _ki_k)$, $[f,\omega
)=\sum_j[f_j,\omega )i_j$, $[f_j,\omega )\in {\bf R}$ for each
$j=0,...,2^v-1$ and $\omega \in {\cal B}(U,{\cal A}_v)$,
$i_0,...,i_{2^v-1}$ is the standard base of generators of the
Cayley-Dickson algebra ${\cal A}_v$. It is well-known that in the
space ${\cal B}'(U,{\bf R})$ of generalized functions the space
${\cal B}(U,{\bf R})$ of test functions is everywhere dense (see
\cite{gelshil} and \S 9 above). In view of the decomposition given
by Formula $(1)$ we get that ${\cal B}(U,{\cal A}_v)$ is everywhere
dense in ${\cal B}'(U,{\cal A}_v)$. Thus sequences of test functions
$F^l$ and $\Xi ^l$ exist converging to $F$ and $\Xi $
correspondingly. \par Without loss of generality we can embed $U$
into ${\cal A}_v$ taking its $\epsilon $-enlargement (open
neighborhood) in case of necessity. So it is sufficient to treat the
case of a domain $U$ in ${\cal A}_v$. In view of the analog of the
Stone-Weierstrass theorem (see \cite{ludfov,ludoyst}) in
$C^0(Q,{\cal A}_v)$ for each compact canonical closed subset in
${\cal A}_v$ the family of all super-differentiable on $Q$ functions
is dense, consequently, the space ${\cal H}(U,{\cal A}_v)$ of all
super-differentiable functions on $U$ is everywhere dense in ${\cal
D}(U,{\cal A}_v)$. For each rectifiable path $\gamma $ in the domain
$U$ a compact canonical closed domain $Q$ exists $Q\subset U$ so
that $\gamma ([a,b])\subset Q$. Therefore, it is sufficient to
consider test functions with compact supports in $Q$. Thus we take
super-differentiable functions $F^n$ and $\Xi ^n$. \par Let $\gamma
^l$ be a sequence of rectifiable paths continuously differentiable,
$\gamma ^l\in C^1([a,b],{\cal A}_v)$, converging to $\gamma $ in
${\bf P}(U)$ relative to the metric $\rho $.
\par Then for any super-differentiable functions $p$ and $q$ we have \par
$(2)$ $\int_{\gamma ^l} p(z) dq(z) = \int_a^b (d\zeta (z)/dz).[(dq
(z)/dz). d\gamma ^l(t)]|_{z=\gamma ^l(t)} $\par $=
\int_a^b\sum_{k=0}^{2^v-1} (\partial \zeta (z)/\partial
z_k)[\sum_{j=0}^{2^v-1}(\partial q_k(z)/\partial z_j)d\gamma
_j^l(t)]$, \\ since each super-differentiable function is Fr\'echet
differentiable, $d\gamma _j^l(t)={\gamma _j^l}'(t)dt$, where
$(d\zeta (z)/dz).1 = p(z)$ and for the corresponding phrases of them
for each $z\in U$. On the other hand, the functional
\par $(3)$ $\int_a^b\sum_{k=0}^{2^v-1} (\partial \zeta (z)/\partial
z_k)[\sum_{j=0}^{2^v-1}(\partial q_k(z)/\partial z_j)d\gamma
_j^l(t)]$ is continuous on ${\cal B} (U,{\cal A}_v)^2\times {\bf
P}(U)$, i.e. for $\zeta , p\in {\cal B} (U,{\cal A}_v)$ and $\gamma
\in {\bf P}(U)$ as well. \par For a rectifiable path $\gamma $ in
$U$ it is possible to take a sequence of open $\epsilon $
neighborhoods $\Gamma ^{\epsilon } := \bigcup_{z\in \gamma ([a,b])}
\breve{B} ({\cal A}_v,z,\epsilon )$, $\epsilon =\epsilon (l) =1/l$,
where $\breve{B} ({\cal A}_v,z,\epsilon ) := \{ y: ~ y\in {\cal
A}_v; ~ |y-z|<\epsilon \} $. Therefore, for each function $\nu $ of
bounded variation on $U$ and each rectifiable path $\gamma $ in $U$
a sequence of test functions $\theta ^l$ with supports contained in
$\Gamma ^{1/l}$ exists such that
\par $\lim_l \int_U [(d\zeta (z)/dz).\theta ^l(z)]\lambda (dz) =
\int_{\gamma } p(z) d\nu (z)$ \\
for each super-differentiable test functions $p, \zeta \in {\cal
H}(U,{\cal A}_v)$ with $(d\zeta (z)/dz).1 = p(z)$ on $U$, where
$\lambda $ denotes the Lebesgue measure on $U$ induced by the
Lebesgue measure on the real shadow ${\bf R}^{2^v}$ of the
Cayley-Dickson algebra ${\cal A}_v$, where ${\cal H}(U,{\cal A}_v)$
denotes the family of all super-differentiable functions on the
domain $U$ with values in the Cayley-Dickson algebra ${\cal A}_v$.
\par Using the latter property and in accordance with Formulas
$(1-3)$ and 16$(G6,G7)$ we put:
\par $(4)$ $[\xi ;\omega ,\gamma ] :=\lim_l [\Xi ^l;\omega ,\gamma ]
= \lim_l \int_{\gamma } \Xi ^l(y) \omega (y)dy$ and
\par $(5)$ $[{\hat f};\omega ^1,\gamma ;\nu ] =\lim_l [dG^l/dz;
\omega ^1, \gamma ;\nu ] = \lim_l \int_{\gamma } \{
(dG^l(z)/dz).\omega ^1(z) \} d\nu (z)$ \\ for any $\nu \in {\cal
V}(U,{\cal A}_v)$, where $(dG^l/dz).1=F^l(z)$ on $U$.
\par Therefore $\Xi ^l$ converges to $\xi $  and
$dG^l/dz$ converges to ${\hat f}$, where $[\xi ;\omega ,*](\kappa
|_{[a,x]}) = [\xi ;\omega ,\kappa |_{[a,x]}]$ for each $\kappa \in
{\bf P}(U)$, $a<x\le b$ (see Lemma 17.2). Therefore, from Formulas
$(2-5)$ and Lemmas 17.1 and 17.2 we infer that
\par $(6)$ $[ \int_{\gamma } fd\eta , \omega ^1\otimes \omega )=
\lim_l [dG^l/dz; \omega ^1, * ; [\Xi ^l;\omega ,* ]](\gamma ^l)
$\par $ = \lim_l \int_{\gamma
^l} [dG^l/dz; \omega ^1, * ; d[\Xi ^l;\omega ,* ](z)] $, \\
where $z=\gamma ^l(t)$, $a\le t\le b$.
\par {\bf 19.1. Corollary.} {\it If $F: U\to {\cal A}_v$
is a continuous function on $U$ and $\Xi $ is a generalized function
on $U$, then  the non-commutative line integral over the
Cayley-Dickson algebra ${\cal A}_v$ of line generalized functions
$f$ by $\xi $ \par $(1)$ $[ \int_{\gamma } fd\xi , \omega ^1\otimes
\omega )$ \\ exists, where $f$ is induced by $F$ and $\xi $ by $\Xi
$.}
\par {\bf Proof.} This follows from Theorem 19 and
the fact that each continuous function $F$ on $U$ gives the
corresponding regular line operator valued generalized function on
the space of test functions $\omega ^1$ in ${\cal B}(U,{\cal A}_v)$
or ${\cal D}(U,{\cal A}_v)$:
\par $[{\hat F};\omega ^1,\gamma ] = \int_{\gamma }
({\hat F}(z).\omega ^1(z))dz$.
\\ In this case one can take the marked function $\omega ^1=\chi _V$,
where $V$ is a compact canonical closed sub-domain in $U$, since
$\gamma ([a,b])$ is compact for each rectifiable path $\gamma $ in
$U$ so that $\gamma ([a,b])\subset V$ for the corresponding compact
sub-domain $V$. This gives ${\hat F}.\chi _V(z)=F(z)$ for each $z\in
V$ and ${\hat F}.\chi _V(z)=0$ for each $z\in U\setminus V$.
\par {\bf 19.2. Corollary.} {\it If $F\in {\cal B}'(U,{\cal A}_v)$
or $F\in {\cal D}'(U,{\cal A}_v)$ is a generalized function on $U$
and $\Xi $ is a function of bounded variation on $U$, then the
non-commutative line integral over the Cayley-Dickson algebra ${\cal
A}_v$ of line generalized functions $f$ by $\xi $ \par $(1)$ $[
\int_{\gamma } fd\xi , \omega ^1\otimes \omega )$ \\ exists, where
$f$ is induced by $F$ and $\xi $ by $\Xi $.}
\par {\bf Proof.} In this case we put \par $[\xi ;\omega ,\kappa ]
:= \int_{\kappa } \omega (z)d\Xi (z)$ \\ for each test function
$\omega $ and each rectifiable path $\kappa $ in $U$. It is
sufficient to take marked test function $\omega (z)=1$ for each
$z\in U$ that gives $d[\xi ;1,*]=d\Xi $. Thus this corollary follows
from Theorem 19.
\par {\bf 19.3. Corollary.} {\it If $F$ is a continuous function on
$U$ and $\Xi $ is a function of bounded variation on $U$, then the
non-commutative line integral over the Cayley-Dickson algebra ${\cal
A}_v$ of line generalized functions $f$ by $\xi $ \par $(1)$ $[
\int_{\gamma } fd\xi , \omega ^1\otimes \omega )$ \\ exists, where
$f$ is induced by $F$ and $\xi $ by $\Xi $. Moreover, this integral
coincides with the non-commutative line integral from \S 3 for the
unit test functions $\omega (z)=\omega ^1(z)=1$ for each $z\in U$:
\par \par $(2)$ $[\int_{\gamma } fd\xi , 1\otimes 1) =
\int_{\gamma }fd\xi $.}
\par {\bf Proof.} This follows from the combination of two
preceding corollaries, since for a rectifiable path $\gamma $ its
image in $U$ is contained in a compact sub-domain $V$ in $U$, i.e.
$\gamma ([a,b])\subset V$.
\par {\bf 19.4. Convolution formula for solutions of partial
differential equations.}
\par Using
convolutions of generalized functions  a solution of the equation
\par $(C1)$ $(\Upsilon ^s +\beta )f = g $ in ${\cal B}({\bf R}^n,Y)$ or in
the space ${\cal B}'({\bf R}^n,Y)_l$ is: \par $(C2)$  $f = {\cal
E}_{\Upsilon ^s+\beta }*g$, \\ where ${\cal E}_{\Upsilon ^s +\beta
}$ denotes a fundamental solution of the equation \par $(C3)$
$(\Upsilon ^s +\beta ){\cal E}_{\Upsilon +\beta }=\delta $,
\\ $(\delta ,\phi )=\phi (0)$ (see \S 9). The fundamental solution of
the equation
\par $(C4)$ $A_0 {\cal V} = \delta $ with
$A_0 = (\Upsilon ^s +\beta ) (\Upsilon ^{s_1}_1+\beta _1)$ \\ can be
written as the convolution
$$(C5)\quad {\cal V} =: {\cal V}_{A_0} = {\cal E}_{\Upsilon ^s +\beta
} * {\cal E}_{\Upsilon ^{s_1}_1+\beta _1}.$$ In view of Formulas
4$(7-9)$ each generalized function ${\cal E}_{\Upsilon ^s +\beta }$
can also be found from the elliptic partial differential equation
\par $(C6)$ $\Xi _{\beta } \Psi _{\Upsilon ^s+\beta } =\delta $ by
the formula:
\par $(C7)$ ${\cal E}_{\Upsilon ^s +\beta } = [(\Upsilon ^s +\beta
)^*] \Psi _{\Upsilon ^s+\beta } $,  where \par $(C8)$ $\Xi _{\beta }
:=(\Upsilon ^s +\beta )(\Upsilon ^s +\beta )^*$
\\ (see \S 33 \cite{lumdltcdla}).
\par {\bf 20. Poly-functionals.} Let ${\bf a}_k: {\cal B}(U,{\cal
A}_r)^k\to {\cal A}_r$ or ${\bf a}_k: {\cal D}(U,{\cal A}_r)^k\to
{\cal A}_r$ be a continuous mapping satisfying the following three
conditions:
\par $(P1)$ $[{\bf a}_k, \omega ^1\otimes ... \otimes \omega ^k)$
is $\bf R$ homogeneous \par $[{\bf a}_k, \omega ^1\otimes ...
\otimes (\omega ^lt)\otimes... \otimes \omega ^k)=[{\bf a}_k,\omega
^1\otimes ... \otimes \omega ^l\otimes... \otimes \omega ^k)t=[{\bf
a}_kt, \omega ^1\otimes ...\otimes \omega ^k)$ \\ for each $t\in
{\bf R}$ and ${\cal A}_r$ additive \par $[{\bf a}_k, \omega
^1\otimes ... \otimes (\omega ^l+\kappa ^l)\otimes... \otimes \omega
^k)=[{\bf a}_k,\omega ^1\otimes ... \otimes \omega ^l\otimes...
\otimes \omega ^k) + [{\bf a}_k, \omega ^1\otimes ... \otimes \kappa
^l\otimes... \otimes \omega ^k)$ \\ by any ${\cal A}_r$ valued test
functions $\omega ^l$ and $\kappa ^l$, when other are marked,
$l=1,...,k$, i.e. it is $k$ ${\bf R}$ linear and $k$ ${\cal A}_r$
additive, where $[{\bf a}_k, \omega ^1\otimes ... \otimes \omega
^k)$ denotes a value of ${\bf a}_k$ on given test ${\cal A}_r$
valued functions $\omega ^1,...,\omega ^k$;
\par $(P2)$ $[{\bf a}_k\alpha ,\omega ^1\otimes ... \otimes (\omega ^l\beta )\otimes...
\otimes \omega ^k)=([{\bf a}_k,\omega ^1\otimes ... \otimes \omega
^l\otimes... \otimes \omega ^k)\alpha )\beta = [({\bf a}_k\alpha
)\beta ,\omega ^1\otimes ... \otimes \omega ^l\otimes... \otimes
\omega ^k)$ for all real-valued test functions and $\alpha , \beta
\in {\cal A}_r$;
\par $(P3)$ $[{\bf a}_k,\omega ^{\sigma (1)}\otimes ...
\otimes \omega ^{\sigma (k)})=[{\bf a}_k,\omega ^1\otimes ...
\otimes \omega ^k) $ for all real-valued test functions and each
transposition $\sigma $, i.e. bijective surjective mapping $\sigma :
\{ 1,...,k \} \to \{ 1,...,k \} $. \par Then ${\bf a}_k$ will be
called the symmetric $k$ ${\bf R}$ linear $k$ ${\cal A}_r$ additive
continuous functional, $1\le k\in {\bf Z}$. The family of all such
symmetric functionals is denoted by ${{\cal B}'}_{k,s}(U,{\cal
A}_v)$ or ${{\cal D}'}_{k,s}(U,{\cal A}_r)$ correspondingly. A
functional satisfying Conditions $(P1,P2)$ is called a continuous
$k$-functional over ${\cal A}_r$ and their family is denoted by
${{\cal B}'}_k(U,{\cal A}_r)$ or ${{\cal D}'}_k(U,{\cal A}_r)$. When
a situation is outlined we may omit for short "continuous" or "$k$
${\bf R}$ linear $k$ ${\cal A}_v$ additive".
\par The sum of two $k$-functionals over the Cayley-Dickson algebra
${\cal A}_r$ is prescribed by the equality:
\par $(P4)$ $[{\bf a}_k+{\bf b}_k, \omega ^1\otimes ... \otimes
\omega ^k) = [{\bf a}_k, \omega ^1\otimes ... \otimes \omega ^k) +
[{\bf b}_k, \omega ^1\otimes ... \otimes \omega ^k)$ \\
for each test functions. Using Formula $(P4)$ each $k$-functional
can be written as
\par $(1)$ $[{\bf a}_k,\omega ^1\otimes ... \otimes \omega ^k)
= [{\bf a}_{k,0}i_0,\omega ^1\otimes ... \otimes \omega ^k)+...+
[{\bf a}_{k,2^r-1}i_{2^r-1},\omega ^1\otimes ... \otimes \omega ^k)$,\\
where $[{\bf a}_{k,j},\omega ^1\otimes ... \otimes \omega ^k)\in
{\bf R}$ is real for all real-valued test functions $\omega
^1,...,\omega ^k$ and each $j$; $~ i_0$,...,$i_{2^r-1}$ denote the
standard generators of the Cayley-Dickson algebra ${\cal A}_r$.
\par The direct product ${\bf a}_k\otimes {\bf b}_p$
of two functionals ${\bf a}_k$ and ${\bf b}_p$ for the same space of
test functions is a $k+p$-functional over ${\cal A}_r$ given by the
following three conditions:
\par $(P5)$ $[{\bf a}_k\otimes {\bf b}_p, \omega
^1\otimes ... \otimes \omega ^{k+p})= [{\bf a}_k, \omega ^1\otimes
... \otimes \omega ^k) [{\bf b}_p, \omega ^{k+1}\otimes ... \otimes
\omega ^{k+p})$ \\ for any real-valued test functions $\omega
^1,...,\omega ^{k+p}$;
\par $(P6)$ if $[{\bf b}_p, \omega
^{k+1}\otimes ... \otimes \omega ^{k+p})\in {\bf R}$ is real for any
real-valued test functions, then \par $[({\bf a}_kN_1)\otimes ({\bf
b}_pN_2), \omega ^1\otimes ... \otimes \omega ^{k+p})= ([{\bf
a}_k\otimes {\bf b}_p, \omega ^1\otimes ... \otimes \omega
^{k+p})N_1)N_2$ \\ for any real-valued test functions $\omega ^1,
...,\omega ^{k+p}$ and Cayley-Dickson numbers $N_1, N_2\in {\cal
A}_r$;
\par $(P7)$ if $[{\bf a}_k, \omega
^1\otimes ... \otimes \omega ^k)\in {\bf R}$ and $[{\bf b}_p, \omega
^{k+1}\otimes ... \otimes \omega ^{k+p})\in {\bf R}$ are real for
any real-valued test functions, then \par $[{\bf a}_k\otimes {\bf
b}_p, \omega ^1\otimes ... \otimes (\omega ^lN_1)\otimes ... \otimes
\omega ^{k+p})= [{\bf a}_k\otimes {\bf b}_p, \omega ^1\otimes ...
\otimes \omega ^{k+p})N_1$ \\ for any real-valued test functions
$\omega ^1, ...,\omega ^{k+p}$ and each Cayley-Dickson number
$N_1\in {\cal A}_r$ for each $l=1,...,k+p$.

\par Therefore, we can now consider a partial
differential operator of order $u$ acting on a generalized function
$f\in {\cal B}'(U,{\cal A}_r)$ or $f\in {\cal D}'(U,{\cal A}_r)$ and
with generalized coefficients either ${\bf a}_{\alpha }\in {\cal
B}'_{|\alpha |}(U,{\cal A}_r)$ or all ${\bf a}_{\alpha } \in {{\cal
D}'}_{|\alpha |}(U,{\cal A}_r)$ correspondingly:
$$(1)\quad Af(x)= \sum_{|\alpha |\le u} (\partial
^{\alpha } f(x))\otimes [({\bf a}_{\alpha }(x))\otimes 1^{\otimes
(u-|\alpha |)}],$$ where $\partial ^{\alpha } f=\partial ^{|\alpha
|}f(x)/\partial x_0^{\alpha _0}...\partial x_n^{\alpha _n}$,
$x=x_0i_0+...x_ni_n$, $x_j\in {\bf R}$ for each $j$, $1\le n=2^r-1$,
$\alpha = (\alpha _0,...,\alpha _n)$, $|\alpha |=\alpha
_0+...+\alpha _n$, $0\le \alpha _j\in {\bf Z}$, $[1,\omega
):=\int_U\omega (y)\lambda (dy)$, $\lambda $ denotes the Lebesgue
measure on $U$, for convenience $1^{\otimes 0}$ means the
multiplication on the unit $1\in {\bf R}$. The partial differential
equation
\par $(2)$ $Af=g$ in terms of generalized functions has a solution $f$
means by the definition that
\par $(3)$ $[Af,\omega ^{\otimes (u+1)}) = [g, \omega ^{\otimes
(u+1)})$ \\ for each real-valued test function $\omega $ on $U$,
where $\omega ^{\otimes k}=\omega \otimes ... \otimes \omega $
denotes the $k$ times direct product of a test functions $\omega $.
\par {\bf 21. Theorem.} {\it Let $A=A_u$ be a partial differential operator
with generalized over the Cayley-Dickson algebra ${\cal A}_r$
coefficients of an even order $u=2s$ on $U$ such that each ${\bf
a}_{\alpha }$ is symmetric for $|\alpha |=u$ and $A$ has the form
\par $(4)$ $Af = (B_{u,1}f)c_{u,1} +...+ (B_{u,k}f)c_{u,k}$, where
each
\par $(5)$ $B_{u,p}=B_{u,p,0}+Q_{u-1,p}$ \\ is a partial
differential operator by variables
$x_{m_{u,1}+...+m_{u,p-1}+1}$,...,$x_{m_{u,1}+...+m_{u,p}}$ and of
the order $u$, $m_{u,0}=0$, $c_{u,k}(x)\in {\cal A}_r$ for each $k$,
its principal part \par $(6)$ $B_{u,p,0}f= \sum_{|\alpha |=s}
(\partial ^{2\alpha }f)\otimes {\bf a}_{p,2\alpha }(x)$ \\ is
elliptic, i.e. $\sum_{|\alpha |=s} y^{2\alpha } [{\bf a}_{p,2\alpha
}, \omega ^{\otimes 2s})\ge 0$ for all
$y_{k(1)}$,...,$y_{k(m_{u,p})}$ in ${\bf R}$ with
$k(1)=m_{u,1}+...+m_{u,p-1}+1$,...,$k(m_{u,p})=m_{u,1}+...+m_{u,p}$,
$y^{\beta } := y_{k(1)}^{\beta _{k_1}}...y_{k(m_{u,p})}^{\beta
_{k(m_{u,p})}}$ and $[{\bf a}_{p,2\alpha }, \omega ^{\otimes 2s})\ge
0$ for each real test function $\omega $, either $0\le r\le 3$ and
$f$ is with values in ${\cal A}_r$, or $r\ge 4$ and $f$ is
real-valued on real-valued test functions. Then three partial
differential operators $\Upsilon ^s$ and $\Upsilon _1^s$ and $Q$ of
orders $s$ and $p$ with $p\le u-1$ with generalized on $U$
coefficients with values in ${\cal A}_v$ exist such that
\par $(7)$ $[Af, \omega ^{\otimes (u+1)})= [\Upsilon ^s
(\Upsilon _1^sf) +Qf, \omega ^{\otimes (u+1)})$ for each real-valued
test function $\omega $ on $U$.}
\par {\bf Proof.} If $a_{2s}$ is a symmetric functional and
$[{\bf c}_s, \omega ^{\otimes s}) = [{\bf a}_{2s}, \omega ^{\otimes
2s})^{1/2}$ for each real-valued test function $\omega $, then by
Formulas 20$(P1,P2)$ this functional ${\bf c}_s$ has an extension up
to a continuous $s$-functional over the Cayley-Dickson algebra
${\cal A}_r$. This is sufficient for Formula $(7)$, where only
real-valued test functions $\omega $ are taken.
\par Consider a continuous $p$-functional ${\bf c}_p$ over ${\cal A}_v$,
$p\in {\bf N}$. Supply the domain $U$ with the metric induced from
the corresponding Euclidean space or the Cayley-Dickson algebra in
which $U$ is embedded. It is possible to take a sequence of
non-negative test functions $\mbox{}_l\omega $ on $U$ with a support
$supp (\mbox{}_l\omega )$ contained in the ball $B(U,z, 1/l)$ with
center $z$ and radius $1/l$ and $\mbox{}_l\omega $ positive on some
open neighborhood of a marked point $z$ in $U$ so that $\int_U
\mbox{}_l\omega (z)\lambda (dz)=1$ for each $l\in {\bf N}$. If the
$p$-functional ${\bf c}_p$ is regular and realized by a continuous
${\cal A}_v$ valued function $g$ on $U^p$, then $\lim_l [{\bf
c}_p,\omega ^{\otimes p})= g(z,...,z)$. Thus the partial
differential equation 20$(2)$ for regular functionals and their
derivatives implies the classical partial differential equation
2$(1)$.
\par Therefore, the statement of this theorem
follows from Theorem 10, and \S \S 14, 15 and 20, since the spaces
of test functions are dense in the spaces of generalized functions
(see \S 19).
\par {\bf 22. Corollary.} {\it If $Af=\sum_{j,k}(\partial ^2f(z)/\partial
z_k\partial z_j) \otimes a_{j,k}(z) + \sum_j (\partial f(z)/\partial
z_j)\otimes b_j(z)\otimes 1 + f(z)\otimes \eta (z)\otimes 1$ is a
second order partial differential operator with generalized
coefficients in ${\cal B}'(U,{\cal A}_r)$ or ${\cal D}'(U,{\cal
A}_r)$, where each $a_{j,k}$ is symmetric, $f$ and ${\cal A}_r$ are
as in \S 20, then three partial differential operators $\Upsilon
+\beta $, $\Upsilon _1+\beta _1$ and $Q$ of the first order with
generalized coefficients with values in ${\cal A}_v$ for suitable
$v\ge r$ of the same class exist such that
\par $(1)$ $[Af, \omega ^{\otimes 3}) =[(\Upsilon +\beta )((\Upsilon _1 +\beta _1)f +Qf),
\omega ^{\otimes 3})$ for each real-valued test function $\omega $
on $U$.}
\par {\bf Proof.} This follows from Theorem 21 and Corollary 12
and \S \S 2 and 8.
\par {\bf 23. Anti-derivatives of first order partial
differential operators with generalized coefficients.} \par {\bf
Theorem.} {\it Let $\Upsilon $ be a first order partial differential
operator given by the formula
\par $(1)$ $\Upsilon f= \sum_{j=0}^n
(\partial f/\partial z_j) \otimes [i_j^* {\psi }_j(z)]$ or \\
\par $(2)$ $\Upsilon f= \sum_{j=0}^n  (\partial f/\partial z_j)\otimes
{\phi }_j^*(z)$,
\\ where $supp ({\psi }_j(z))=U$ or $supp ({\phi }_j(z))=U$ for each
$j$ respectively, $f$ and ${\psi }_j(z)$ or ${\psi }_j(z)$ are
${\cal A}_v$-valued generalized functions in ${\cal B}'(U,{\cal
A}_r)$ or ${\cal D}'(U,{\cal A}_r)$ on the domain $U$ satisfying
Conditions 1$(D1,D2)$, $ ~ alg_{\bf R} \{ [\phi _j, \omega ), [\phi
_k,\omega ), [\phi _l,\omega ) \} $ is alternative for all $0\le j,
k, l \le 2^v-1$ and $alg_{\bf R} \{ [\phi _0,\omega ),...,[\phi
_{2^v-1},\omega ) \} \subset {\cal A}_v$ for each real-valued test
function $\omega $ on $U$. Then its anti-derivative operator ${\cal
I}_{\Upsilon }$ exists such that $\Upsilon {\cal I}_{\Upsilon }f=f$
for each continuous generalized function $f: U\to {\cal A}_v$ and it
has an expression through line integrals of generalized functions.}
\par {\bf Proof.} When an operator with generalized coefficients
is given by Formula $(1)$, we shall take unknown generalized
functions $\nu _j(z)\in {\cal A}_v$ as solutions of the system of
partial differential equations by real variables $z_k$:
\par $(3)$ $[(\partial \nu _j(z)/\partial z_j)\otimes \psi _j(z),
\omega ^{\otimes 2}) = [1, \omega ^{\otimes 2})$ for all $1\le j\le
n$;
\par $(4)$ $[\psi _k(z) \otimes (\partial \nu _j(z)/\partial z_k),
\omega ^{\otimes 2}) = [\psi _j(z)\otimes (\partial \nu
_k(z)/\partial z_j), \omega ^{\otimes 2}) $ for all $1\le j<k\le n$
and and real-valued test functions $\omega $ on $U$.
\par If the operator is given by Formula $(2)$ we consider the
system of partial differential equations:
\par $(5)$ $[((dg(z)/dz).[\partial \nu _j(z)/\partial z_k])\otimes \phi
_k^*(z) + ((dg(z)/dz).[\partial \nu _k(z)/\partial z_j])\otimes \phi
_j^*(z), \omega ^{\otimes 2}) =0$ for all $0\le j<k\le n$;
\par $(6)$ $\partial \nu _j(z)/\partial z_j = \phi _j(z)$ for all
$j=0,....,n$; \par $(7)$ $[((dg(z)/dz).\phi _j(z))\otimes \phi
_j^*(z), \omega ^{\otimes 2}) =[f(z)\otimes 1,\omega ^{\otimes 2}) $
for each $j=0,...,n$ and every real-valued test function $\omega $.
\par Certainly the system of differential equations given by Formulas
$(3,4)$ or $(5-7)$ have solutions in the spaces of test functions
${\cal B}(U,{\cal A}_r)$ or ${\cal D}(U,{\cal A}_r)$, when all
functions $\psi _j$ or $\phi _j$ are in the same space respectively.
Applying \S \S 4 or 5 we find generalized functions $\nu _j$
resolving these system of partial differential equations
correspondingly, when all functions $\psi _j$ or $\phi _j$ are
generalized functions, since the spaces of test functions are dense
in the spaces of generalized functions (see \S 19). Substituting
line integrals $\int_{\gamma } q(y)d\nu _j(y)$ from \S \S 4 and 5 on
line integrals $[\int_{\gamma } q(y)d\nu _j(y), \omega ^1\otimes
\omega )$ from \S 19 one gets the statement of this theorem, since
test functions $\omega ^1$ and $\omega $ in the line integrals of
generalized functions can also be taken real-valued and the real
field is the center of the Cayley-Dickson algebra ${\cal A}_v$.
Therefore, we infer that \par $(8)$ $\partial [\int_{\gamma ^{\alpha
}|_{<a_{\alpha },t_z]}} f(y)d\nu _j(y), \omega \otimes \omega )
/\partial z_k =
[{\hat f} (z).[d\nu _j(z)/dz_k], \omega \otimes \omega ) $ \\
for each real-valued test function $\omega $ and $z\in U$, where
$\gamma ^{\alpha }(t_z)=z$, $ ~ t_z\in <a_{\alpha }, b_{\alpha }>$,
$ ~ \alpha \in \Lambda $. Equality $(8)$, Theorem 19 and Corollaries
19.1-19.3 and Conditions 20$(P1-P7)$ give the formula for an
anti-derivative operator:
$$(9)\quad [{\cal I}_{\Upsilon }f,\omega \otimes \omega ) =
[\mbox{}_{\Upsilon }\int f(z)dz,\omega \otimes \omega ) = (n+1)^{-1}
\sum_{j=0}^n \{ [\int_{\gamma ^{\alpha }|_{[a_{\alpha },t]}}
q(y)d\nu _j(y), \omega \otimes \omega )$$ for each real-valued
test-function $\omega $, where $\alpha \in \Lambda $, $a_{\alpha
}\le t\le b_{\alpha }$, $t=t_z$, $z=\gamma (t)$, consequently,
\par $(10)$ $[\Upsilon \mbox{}_{\Upsilon }\int f(y)dy, \omega ^{\otimes
3})= [f\otimes 1\otimes 1, \omega ^{\otimes 3})$.
\par {\bf 23.1. Note.} Certainly, the case of the partial differential operator
\par $(1)$ $\Upsilon f= \sum_{j=0}^n  (\partial f/\partial z_{k(j)})\otimes
{\phi }_{k(j)}^*(z)$, \\ where $0\le k(0)<k(1)<...<k(n)\le 2^v-1$
reduces to the considered in \S 23 case by a suitable change of
variables $z\mapsto y$ so that $z_{k(j)}=y_j$.
\par {\bf 24. Example.} We consider a consequence of Formulas
15$(2-6)$. If $q(t)$ is a differentiable function on the real field
$\bf R$ having simple zeros $q(t_j)=0$ (i.e. zeros of the first
order), then \par $(1)$ $\delta (q(t)) = \sum_j \frac{1}{q(t_j)}
\delta (t-t_j)$, \\ where the sum is accomplished by all zeros $t_j$
of the equation $q(t)=0$ (see Formula II.2.6$(IV)$ \cite{gelshil}).
 Therefore, if $\gamma (\tau )$ is a $C^1$ path in $U$
intersecting the surface $\partial U_s\cap \partial U_p$ at the
marked point $x$ of index $l=l_{s,p}(x)$, $\gamma (\tau _0)=x$,
$0<\tau _0<1$, such that $\gamma (\tau )\in U_s$ for each $\tau
<\tau _0$ and $\gamma (\tau )\in U_p$ for each $\tau >\tau _0$ then
\par $(2)$ $dg(\gamma (\tau ))/d\tau = dg^1(\gamma (\tau ))/d\tau +
h\sum_{j=1}^l \delta (P_j)dP_j(\gamma (\tau ))/d\tau $, \\
where $\theta (t)=0$ for $t<0$ and $\theta (t)=1$ for $t\ge 0$,
$g^1(\gamma (\tau )) = g(\gamma (\tau )) - h \theta (\tau -\tau
_0)$, \par $(2.1)$ $h= \lim_{\tau \downarrow \tau _0}g(\gamma (\tau
)) - \lim_{\tau \uparrow \tau _0}g(\gamma (\tau ))$, \par $(2.2)$
$dP_j(\gamma (\tau ))/d\tau = \sum_{k=1}^n (\partial P_j(z)/\partial
z_k)(\partial \gamma _k(\tau )/d\tau )|_{z=\gamma (\tau )} $  (see
Example I.2.2.2 \cite{gelshil}). Particularly, if a point $x$ is of
index $1$, then Formula $(2)$ simplifies:
\par $(3)$ $dg(\gamma (\tau ))/d\tau = dg^1(\gamma (\tau ))/d\tau +
h \delta (P)[dP(\gamma (\tau ))/d\tau ]$. Particularly, these
formulas can be applied to $d\nu _j$.
\par Let a partial differential differential operator $Q$ be given
by Formula 14$(7)$ and functions $\nu _k$ are found (see Theorems 5
and 21 above). We put in accordance with Formula 15$(6)$
\par $(4)$ $\nu _k^1(z) = \nu _k(z) - \sum_{s,p} h_{k;s,p}(z)\theta
(P_{s,p;1}(z),...,P_{s,p;m_j}(z))$, \\ where
\par $(5)$ $h_{k;s,p}(z(x))=h_k(x)|_{x\in \Gamma ^{s,p}} := $\par
$\lim_{y_{1;s,p}\downarrow 0,...,y_{n;s,p}\downarrow 0} \nu
_k(z(x+y)) - \lim_{y_{1;s,p}\downarrow 0,...,y_{n;s,p}\downarrow 0}
\nu _k(z(x-y))$, \\ where the sum is by $s$ and $p$ with $\partial
U_s\cap \partial U_p\ne \emptyset $.
\\ Let a domain $W$ be a canonical closed compact set in the
Euclidean space ${\bf R}^{n+1}$ embedded into ${\cal A}_v$ and
contained in a canonical closed compact domain $U$ so that $W= \{
z\in U: ~ z_j=0 ~ \forall n<j\le 2^v-1 \} $. Thus $\Upsilon $ from
test and generalized functions on $W$ is extended on test and
generalized functions on $U$. We can put $\nu _j=0$ for $n<j\le
2^v-1$, when $n<2^v-1$. Then for the rectifiable path $\gamma $ (see
above) we get
$$(6)\quad (n+1)^{-1} \sum_{k=0}^n [\int_{\gamma |_{[a,t]}} q(y)d\nu
_k(y), \omega ^1\times \omega )$$
$$ = (n+1)^{-1} \sum_{k=0}^n \{ [\int_{\gamma } q(y)d\nu ^1_k(y),
\omega ^1\times \omega )$$  $$ + [{\hat q};\omega ^1,*; [\sum_{s,p}
h_{k;s,p}(z)\sum_{j=1}^l \delta (P_j)\sum_{m=1}^n (\partial
P_j(z)/\partial z_m)(\partial \gamma _m(\tau )/d\tau )|_{z=\gamma
(\tau )};\omega ,*]](\gamma ) \} ,$$ where $\gamma \in \{ \gamma
^{\alpha }: ~ \alpha \in \Lambda \} $ is taken from the foliation
$C^1$ family of paths (see \S 6.1 above and also Theorem 2.13
\cite{luviralgoclanl}), $z=\gamma (t_z)$, $t_z\in [a_{\alpha },
b_{\alpha }]$, $[a,t]=[a_{\alpha },t_z]$, $l=l_{s,p}(z)$ denotes an
index of a point $z$ in the intersection of boundaries $\Gamma
^{s,p} :=
\partial U_s\cap
\partial U_p\ne \emptyset $, $ ~ \omega ^1$ and $\omega $ are
real-valued test functions. Since $\omega ^1$ is real-valued, we get
${\hat f}(z).\omega ^1=f(z)\omega ^1(z)$ and
\par $(7)$ $[{\hat q};\omega ^1,*; [\sum_{s,p} h_{k;s,p}(z)\sum_{j=1}^l
\delta (P_j)\sum_{m=1}^n (\partial P_j(z)/\partial z_m)(\partial
\gamma _m(\tau )/d\tau )|_{z=\gamma (\tau )};\omega ,*]](\gamma
)$\par $ = \sum_{s,p} [q(z)\omega ^1(z); [\theta (\tau -\tau _{s,p})
h_{k;s,p}(z), \omega ))|_{z=\gamma (\tau )}$, \\ where $\tau _{s,p}$
corresponds to the intersection point $\gamma (\tau _{s,p})$ of
$\gamma $ with $\Gamma ^{s,p}\ne \emptyset $. Here the expression
$[q,\omega )|_{z=\gamma (\tau )} := \lim_j [q\circ \kappa ^j, \omega
\circ \kappa ^j)$ denotes the restriction of the generalized
function from $U$ onto $\gamma ([a,b])$, $ ~ \kappa ^j\in {\cal
D}(U,{\cal A}_v)$ is a sequence of test functions and $\kappa
^j(\phi ([a,b]))\subset U$ for each $j\in {\bf N}$, $~ \phi \in
{\cal D}([a,b],{\cal A}_v)$, $ ~ \bigcap_{j=1}^{\infty } supp
(\kappa ^j)= \phi ([a,b])$, $\lim_j \kappa ^j\circ \phi =\gamma $ in
${\bf P}(U)$. Therefore, the derivative of the operator $[(n+1)^{-1}
\sum_{k=0}^n [\int_{\gamma |_{[a,t]}} q(y)d\nu _k(y), \omega
^1\times \omega )$ by the parameter $\tau \in [a,b]$ for the real
test functions $\omega ^1$ and $\omega $ is the following:
\par $(8)$ $\partial (n+1)^{-1} \sum_{k=0}^n [\int_{\gamma |_{[a,t]}} q(y)d\nu
_k(y), \omega ^1\times \omega )/\partial \tau =$\par $ [(n+1)^{-1}
\sum_{k=0}^n \{ 1\otimes {\hat q}^1(z).(d\nu ^1_k(\gamma (\tau
))/d\tau ) + \sum_{s,p} (h^{g'}_{s,p}(z).h_{k;s,p}(z) + {\hat
q}^1(z).h_{k;s,p}(z) + h^{g'}_{s,p}(z).(d\nu ^1_k(\gamma (\tau
))/d\tau ))|_{z=\gamma (\tau )})\otimes \delta (\tau -\tau _{s,p}),
\omega ^1\otimes \omega )$, \\
where $dg(z)/dz={\hat q}(z)$ on $U$ in the class of generalized
operator $L_q({\cal A}_v,{\cal A}_v)$ valued functions,
$(dg(z)/dz).1=q(z)$ on $U$, $~ h^{g'}_{s,p}(z)=h(z)$ is given by
Formula $(2.1)$ for the derivative operator $dg(z)/dz=g'$ instead of
$g$ on each $\Gamma ^{s,p}\ne \emptyset $, ${\hat q}^1$ is given by
Formula 15$(6)$ for the function ${\hat q}(z)$ with values in
$L_q({\cal A}_v,{\cal A}_v)$ instead of $g(z)$. The terms like
${\hat q}^1(z).(d\nu ^1_k(\gamma (\tau ))/d\tau )$ correspond to the
action of the operator valued generalized function ${\hat q}^1(z)$
on the generalized function $(d\nu ^1_k(\gamma (\tau ))/d\tau )$
which gives a generalized function. \par Using Formulas $(6-8)$ for
$n$ constant on $U$ and $\psi _j(z)$ or $\phi _j(z)$ respectively
non-zero for each $z\in U$ and all $j=0,...,n$ we infer that for a
continuous or generalized function $f$ \par $(9)$ $\Upsilon {\cal
I}_{\Upsilon }f(z)=f(z)$, where
\par $$(10)\quad \mbox{}_{\Upsilon }\int f(z)dz := \{ (n+1)^{-1}
\sum_{j=0}^n \int_{\gamma ^{\alpha }|_{[a_{\alpha },t]}} q(z)d\nu ^1
_j(z), \quad \alpha \in \Lambda , ~ a_{\alpha }\le t \le b_{\alpha }
\} ,$$ where $q=(dg/dz).1$ and $g$ is given by the Equation 5$(9)$,
since $f^1=f$ and $h^{g'}_{s,p}=0$ in the class of generalized
functions $f$ and in the class of continuous functions $f$, also
$h_{k;s,p}=0$ for $\nu _k=\nu ^1_k$ on $U$.
\par Formulas $(9,10)$ show what sort of boundary conditions is
sufficient to specify a unique solution for a given domain $U$ with
sub-domains $U_s$. If $U$ is $C^1$ diffeomorphic to the half-space
$H_p := \{ z\in {\cal A}_v: ~  z_0p_0+...+z_{2^v-1}p_{2^v-1}\ge 0 \}
$, where $p=p_0i_0+...+p_{2^v-1}i_{2^v-1}$ is a marked
Cayley-Dickson number, $p_0,...,p_{2^v-1}\in {\bf R}$, and
sub-domains $U_s$ are not prescribed, then it is sufficient to give
the boundary condition $F|_{\partial U}=G$ when a solution is in the
class of continuous or generalized functions with the corresponding
$f$ and $\psi _j$ or $\phi _j$. Indeed, if the functions $\nu _k$
along $\gamma ^{\alpha }$ are defined up to constants $\mu _k$, the
differentials are the same $d(\nu _k+\mu _k)(z)|_{z=\gamma ^{\alpha
}(\tau )}=d\nu _k|_{z=\gamma ^{\alpha }(\tau )}$ in the
anti-derivative operator, when $d\mu _k|_{z=\gamma ^{\alpha }(\tau
)}=0$ for each $\alpha \in \Lambda $ and $\tau \in [a_{\alpha
},b_{\alpha }]$. \par The operator ${\cal I}_{\Upsilon }$ may be
applied also piecewise on each $U_s$. If a solution $F$ is locally
continuous on $U$ and continuous on each sub-domain $U_s$, then
boundary conditions $F|_{\partial U_s}=G^s$ for all $s=1,...,m$ may
be necessary to specify a solution $F$.  Without boundary conditions
the anti-derivative operator applied to $f$ gives the general
solution ${\cal I}_{\Upsilon }f$ of the differential equation
$\Upsilon F =f$.
\par If each $\nu
_k$ is continuously differentiable, which is possible, when each
function $\psi _k$ or $\phi _k$ is continuous, and $f$ is continuous
on $U$, then a solution $F=\{ \mbox{}_{\Upsilon }\int_{\gamma
^{\alpha }|_{[a_{\alpha },\tau ]}} f(y)dy: ~ \alpha \in \Lambda , ~
a_{\alpha } \le \tau \le b_{\alpha } \} $ is continuously
differentiable by each $z_k$, $z\in U$, $~ z=\gamma ^{\alpha
}(t_z)$. \par One can also mention that the sequence
$\mbox{}_m\omega (z_0,...,z_n) = (2\pi m)^{-(n+1)/2}\exp \{ -
(z_0^2+...+z_n^2)/(2m) \} $ converges to the delta-function on ${\bf
R}^{n+1}$ embedded into ${\cal A}_v$,  $m\in {\bf N}$.  Then the
sequence $\mbox{}_m\theta (z_0) = \int_{-m}^{z_0} (2\pi m)^{-1/2}
\exp \{ - t^2/(2m) \} dt $ converges to the $\theta $ function,
while each function $\mbox{}_m\theta (z_0)$ is analytic by $z_0$,
since $m\in {\bf N}$, $\exp \{ - t^2/(2m) \}$ is the analytic
function with the infinite radius of convergence of its power
series, while $\lim_{m\to + \infty } \int_{- \infty }^{-m}(2\pi
m)^{-1/2} \exp \{ - t^2/(2m) \} dt =0$. Then each $\mbox{}_m\omega
(z_0,...,z_n)$ and $\prod_{j=0}^n\mbox{}_m\theta (z_j)$ can be
written in the $z$-representation over ${\cal A}_v$ as the analytic
function with the help of Formulas 1$(1-3)$, where $n\le 2^v-1$,
$z_j\in {\bf R}$ for each $j$, $z=z_0i_0+...+z_{2^v-1}i_{2^v-1}$.
Thus $\lim_{m\to \infty } \mbox{}_m\omega (z_0,...,z_n) = \delta
(z_0,...,z_n)$ and $\lim_{m\to \infty } \prod_{j=0}^n\mbox{}_m\theta
(z_j)= \prod_{j=0}^n\theta (z_j)$, $ ~ d\theta (z_0)/dz_0=\delta
(z_0)$.
\par {\bf 25. Boundary conditions.}
\par If $U$ is a domain as in \S 15, we put ${\cal B}(\partial U,Y) =
\{ f|_{\partial U}: ~ f\in {\cal B}(U,Y) \} $ and  ${\cal
D}(\partial U,Y) = \{ f|_{\partial U}: ~ f\in {\cal D}(U,Y) \} $
when a boundary $\partial U$ is non-void so that the topologically
adjoint linear over ${\bf R}$ spaces, left and right ${\cal A}_r$
modules, of generalized functions are ${\cal B}'(\partial U,Y)$ and
${\cal D}'(\partial U,Y)$. \par Let us consider a generalized
function $f\in {\cal B}'(\partial U,Y)$ or ${\cal D}'(\partial U,Y)$
and a test function $h\in {\cal B}(\partial U,Y)$ or ${\cal
D}(\partial U,Y)$ respectively. One can take $g\in {\cal B}(U,Y)$ or
${\cal D}(U,Y)$ and a sequence $q^m\in {\cal B}(U,Y)$ or ${\cal
D}(U,Y)$ with supports non intersecting with the boundary $supp
(q^m)\cap \partial U=\emptyset $ such that $(g-q^m)$ tends to zero
in ${\cal B}(V,Y)$ or ${\cal D}(V,Y)$ for each compact subset $V$ in
the interior $Int (U)$, when $m$ tends to the infinity, while
$\lim_m (g-q^m) =h$ in ${\cal B}(\partial U,Y)$ or ${\cal
D}(\partial U,Y)$ respectively. Here as usually the interior $Int
(U)$ is taken in the corresponding topological space ${\bf R}^n$ or
${\cal A}_r$. Each generalized function is a limit of test
functions, consequently, a generalized function $\xi \in {\cal
B}'(U,Y)$ or ${\cal D}'(U,Y)$ exists so that \par $(B1)$ $\lim_m
[\xi , (g-q^m)) = [f,h)$. \\ Vise versa if $\xi \in {\cal B}'(U,Y)$
or ${\cal D}'(U,Y)$ is a generalized function on $U$, then Formula
$(B1)$ defines  a generalized function $f\in {\cal B}'(\partial
U,Y)$ or ${\cal D}'(\partial U,Y)$, which we consider as  the
restriction of $\xi $ on ${\cal B}(\partial U,Y)$ or ${\cal
D}(\partial U,Y)$ correspondingly. In view of the definition of
convergence of test and generalized functions Formula $(B1)$ defines
the unique restriction $f$ for the given generalized function $\xi
$.
\par A subsequent use of decomposition of operators into compositions
of first order partial differential operators and their
anti-derivation operators permits to integrate partial differential
equations with locally continuous or generalized coefficients.
\par The results and definitions of previous sections show that for the
differential equation
\par $(1)$ $Af=g$, \\ where a partial differential operator is written
in accordance with Formulas 10$(1,2)$. When $\partial U$ is a
$C^1$-manifold without corner points of index greater than one, the
following boundary conditions may be used:
\par $(2)$ $f(t)|_{\partial U} = f_0(t'),$
$(\partial ^{|q|}f(t)/\partial s_1^{q_1}...\partial s_n^{q_n}
)|_{\partial U} = f_{(q)}(t')$ for $|q|\le \alpha -1$, where
$s=(s_1,...,s_n) \in {\bf R}^n$, $(q)=(q_1,...,q_n)$, $|q|
=q_1+...+q_n$, $0\le q_k\in {\bf Z}$ for each $k$, $t\in \partial U$
is denoted by $t'$, $f_0$, $f_{(q)}$ are given functions. Generally
these conditions may be excessive, so one uses some of them or their
linear combinations (see $(4)$ below). Frequently, the boundary
conditions
\par $(3)$ $f(t)|_{\partial U} = f_0(t'),$
$(\partial ^lf(t)/\partial \nu ^l)|_{\partial U} = f_l(t')$ for
$1\le l\le \alpha -1$ are also used, where $\nu $ denotes a real
variable along a unit external normal to the boundary $\partial U$
at a point $t'\in \partial U_0$. Using partial differentiation in
local coordinates on $\partial U$ and $(3)$ one can calculate in
principle all other boundary conditions in $(2)$ almost everywhere
on $\partial U$. \par It is possible to describe as an example a
particular class of domains and boundary conditions. Suppose that a
domain $U_1$ and its boundary $\partial U_1$ satisfy Conditions
$(D1,i-vii)$ and $g_1=g\chi _{U_1}$ is a regular or generalized
function on ${\bf R}^n$ with its support in $U_1$. Then any function
$g$ on ${\bf R}^n$ gives the function $g\chi _{U_2}=: g_2$ on ${\bf
R}^n$, where $U_2={\bf R}^n\setminus U_1$.  Take now new domain $U$
satisfying Conditions $(D1,i-vii)$ and $(D2-D5)$:
\par $(D2)$ $U\supset U_1$ and $\partial U\subset \partial U_1$;
\par $(D3)$ if a straight line $\xi $ containing a point $w_1$ (see
15$(vi)$) intersects $\partial U_1$ at two points $y_1$ and $y_2$,
then only one point either $y_1$ or $y_2$ belongs to $\partial U$,
where $w_1\in U_1$, $U-w_1$ and $U_1-w_1$ are convex; if $\xi $
intersects $\partial U_1$ only at one point, then it intersects
$\partial U$ at the same point;
\par $(D4)$  any straight line $\xi $ through the point $w_1$ either does not intersect
$\partial U$ or intersects the boundary $\partial U$ only at one
point. \par Take now $g$ with $supp (g)\subset U$, then $supp (g
\chi _{U_1})\subset U_1$. Therefore, any problem $(1)$ on $U_1$ can
be considered as the restriction of the problem $(1)$ defined on
$U$, satisfying $(D1-D4,i-vii)$. Any solution $f$ of $(1)$ on $U$
with the boundary conditions on $\partial U$ gives the solution as
the restriction $f|_{U_1}$ on $U_1$ with the boundary conditions on
$\partial U$. \par Henceforward, we suppose that the domain $U$
satisfies Conditions $(D1,D4,i-vii)$, which are rather mild and
natural.
\par Thus the sufficient boundary conditions are:
\par $(4)$ $(\partial ^{|\beta |} f(t^{(lj)})/\partial \tau _{\gamma (1)}^{\beta _1}...\partial
\tau _{\gamma (m)}^{\beta _m})|_{\partial U_{(lj)}} = \phi _{\beta
,(lj)}(t^{(lj)})$ \\ for $|\beta | = |q|$, where $m=|h(lj)|$,
$|j|\le \alpha $, $|(lj)|\ge 1$, ${\bf a}_j\ne 0$, $q_k=0$ for
$l_kj_k=0$, $m_k+q_k+h_k=j_k$, $~ h_k = sign (l_kj_k)$, $0\le q_k\le
j_k-1$ for $k>n-\kappa $; $\phi _{\beta ,(l)}(t^{(l)})$ are known
functions on $\partial U_{(l)}$, $t^{(l)}\in \partial U_{(l)}$. In
the half-space $t_n\ge 0$ only the partial derivatives by $t_n$
\par $(5)$ $\partial ^{\beta }f(t)/\partial t_n^{\beta }|_{t_n=0}$ \\ are
necessary for $\beta =|q|<\alpha $ and $q$ as above.
\par Depending on coefficients of the operator $A$ and the domain
$U$ some boundary conditions may be dropped, when the corresponding
terms vanish.
\par Conditions in $(5)$ are given on disjoint for different
$(l)$ sub-manifolds $\partial U_{(l)}$ in $\partial U$ and partial
derivatives are along orthogonal to them coordinates in ${\bf R}^n$,
so they are correctly posed. \par We recall, that a characteristic
surface of a partial differential operator given by Formula 10$(1)$
is a surface defined as a zero of $C^u$ differentiable function
$\phi (x_1,...,x_n)=0$ in the Euclidean space or in the
Cayley-Dickson algebra such that at each point $x$ of it the
condition is satisfied
\par $(CS)$ $\sum_{|j|=u} {\bf a}_j(t(x)) (\partial \phi /\partial x_1)^{j_1}...
(\partial \phi /\partial x_n)^{j_n} =0  $ \\
and at least one of the partial derivatives $(\partial \phi
/\partial x_k)\ne 0$ is non-zero. Generally a domain $U$ is
worthwhile to choose with its interior $Int (U)$ non-intersecting
with a characteristic surface $\phi (x_1,...,x_n)=0$ (see also
\cite{petrov,vladumf}).
\par {\bf 26. Solutions of second order partial differential equations
with the help of the line integration over the Cayley-Dickson
algebras.}
\par Mention that the operator $(\Upsilon +\beta )(z_0,...,z_n)$
acting on a function depending on variables $z_0,...,z_n$ only can
be written as \par $(1)$ $\Psi (z_0,...,z_{n+1})
(f(z)z_{n+1})|_{z_{n+1}=1} = \Upsilon (z_0,...,z_n) f(z)+f(z)\beta
(z)$\par $= [\sum_{j=0}^{n+1} (\partial (f(z)z_{n+1})/\partial
z_j)\phi _j^*(z)]|_{z_{n+1}=1},$
\\ where $\phi _{n+1}^*(z)=\beta (z)$, each function
$\phi _j(z)$ and $f(z)$ may depend on $z_0,...,z_n$ only, omitting
for short the direct product $\otimes $ in the case of generalized
coefficients in Formula $(1)$ and henceforth. The operator $\Psi
(z_0,...,z_{n+1}) (f(z)z_{n+1})$ may be reduced to the form
satisfying conditions of Theorems 5 or 23 using a suitable change of
variables.
 This procedure gives an anti-derivative operator \par $(1.1)$
${\cal I}_{\Upsilon +\beta }= {\cal I}_{\Psi }|_{z_{n+1}=1}$ such
that \par $(1.2)$ $(\Upsilon + \beta ){\cal I}_{\Upsilon +\beta
}f=f$ \\ for a continuous function or a generalized function $f$.
Therefore, we shall consider operators of the form $\Upsilon $ and
their compositions and sums.
\par  We take the partial differential equation
with piecewise continuous or generalized coefficients
\par $(2)$ $A= \Upsilon \Upsilon _1 f(z) + \Upsilon _2f(z) =
g$,\\ where
\par $(3)$ $\Upsilon _kf(z)=
[\sum_{j=0}^n (\partial f(z)/\partial z_j)\phi ^k_j(z)^*]$ \\ for
$k=1, 2$ or without this index $k$, $ ~ \phi ^k_j(z)=i_j\psi
^k_j(z)$ (see \S \S 4, 5 and 23). For solving it we write the
system:
\par $(4)$  $\Upsilon _1 f = g_1$, $\Upsilon g_1
= g - \Upsilon _2 f$. \par In accordance with Equation $(2)$ we
have: \par $(5)$ $\Upsilon g_1 + \Upsilon _2 (\Upsilon _1^{-1}g_1) =
g$, \\ where the inverse operator $\Upsilon _1^{-1}$ is the
anti-derivation operator ${\cal I}_{\Upsilon _1}$ described above in
Theorems 4, 5 and 23. If $\Upsilon _2\ne 0$ we suppose that either
$(i)$ $g_1$ or $[g_1,\omega )$ is real-valued or the algebra
$alg_{\bf R} \{ g_1(z), \phi ^1_j(z), \phi ^2_k(z) \} $ or $alg_{\bf
R} \{ [g_1, \omega ), [\phi ^1_j, \omega ), [\phi ^2_k,\omega ) \} $
for all $j, k=0,...,n$ is associative for each $z\in U$ in the case
of functions or for every real-valued test function $\omega $ in the
case of generalized functions correspondingly. Calculating the
expression
\par $(6)$ $\Upsilon g_1 + \Upsilon _2 ({\cal I}_{\Upsilon _1}g_1)
=(\Upsilon +\beta ^3)g_1=g$ \\
we get a new operator $(\Upsilon +\beta ^3) = \Psi $ in accordance
with Formulas 4$(5,6)$, 5$(6-8)$ and 23$(3-7)$ omitting $\otimes $
and $\otimes \omega $ and $[*,\omega ^{\otimes 3})$ for short in the
case of generalized coefficients, where \par $\beta ^3=
\sum_{j,k=0}^n (\partial \nu ^1_j/\partial z_k) (\phi ^2_k)^*$, \\
when a solution $ \{ \nu _j(z): ~ j \} $ is chosen real, i.e. each
function $\nu _j$ is real-valued or a real-valued generalized
function on real-valued test functions (see System 4$(10,11)$ and \S
23).
\par Generally without supposition $(i)$ we deduce that
\par $(6.1)$ $\Upsilon g_{1,l} + \Upsilon _2 ({\cal I}_{\Upsilon _1}g_{1,l})
= (\Upsilon +\beta ^3_l) g_{1,l} = i_l^*g^l$ \\
for each $l=0,...,n$, where $\beta ^3_l=  i_l^* [\sum_{j,k=0}^n(i_j
(\partial \nu ^1_j/\partial z_k))(\phi ^2_k)^*]$, \par $(6.2)$
$\sum_{l=0}^n g^l=g$,
\par $(6.3)$ $\sum_{l=0}^n g_{1,l} i_l=g_1$, \\ each
$g_{1,l}(z)$ or $[g_{1,l},\omega )$ is real-valued for each $z\in U$
or every real-valued test function $\omega $ respectively. Solving
the system $(6.1-3)$ with the help of known anti-derivative
operators ${\cal I}_{\Psi _l}$ one finds $g_1$, where $\Psi
_l=\Upsilon +\beta ^3_l$. Thus an anti-derivative operator ${\cal
J}:={\cal J}_{\Upsilon  + \Upsilon _2 {\cal I}_{\Upsilon _1}}$
exists so that \par $(6.4)$ $(\Upsilon  + \Upsilon _2 {\cal
I}_{\Upsilon _1}) {\cal J}_{\Upsilon + \Upsilon _2 {\cal
I}_{\Upsilon _1}}g=g$.
\\ In the particular case $(i)$ the equality ${\cal J}_{\Upsilon +
\Upsilon _2 {\cal I}_{\Upsilon _1}} = {\cal I}_{\Psi }$ is
satisfied.
\par
Therefore, in the case of either continuous coefficients of
operators and $g$ or generalized coefficients and $g$ the general
solution is:
\par $(7)$ $f=  {\cal I}_{\Upsilon _1} g_1 =  {\cal I}_{\Upsilon
_1}{\cal J}g,$ where
\par $(8)$ $g_1={\cal J}g.$
\par If $\Upsilon _2=0$ the formula simplifies to
\par $(9)$ $f= {\cal I}_{\Upsilon _1}{\cal I}_{\Upsilon }g.$
\par Examples of boundary conditions and domains permitting to specify
a unique solution are given in \S 25 above.

\par {\bf 27. Example.} Let us consider a function and its phrase
satisfying Condition 7$(P3)$. Therefore, we get
\par $(1)$  $\Upsilon _1^{-1}(g_1)(z_{a^1}(x)) =  - \phi _{f'} (Im (z_{a^1}(x)))
+ \int_{\alpha }^x g_1(z)dz + \phi _{g_1}(x')$ \\ in accordance with
Formula 7$(6)$, where $x'$ and $Im (z_{a^1}(x))$ are written in the
$x$ and $z$-representations respectively using Identities 1$(1-3)$.
In particular, for $\Upsilon ^1=\sigma _1$ we have the coefficients
$\psi ^1_l(z) = i_l ({\bf a}_k(z) w_j i_{2^rk}^*)$ for each $l=2^rk$
with $k\in \{ m_1+...+m_{j-1}+1,..., m_1+...+m_j \} $, while $\psi
^1_l(z)=0$ for all others $l$ for each $z$. A function $a^1$ is
given by Formula 7$(3)$ for $\psi ^1$ instead of $\psi $. Let the
first order operator $Q$ be written in its standard form:
$$(2)\quad Qf = \sum_{j=1}^m
\sum_{k=m_1+...+m_{j_1-1}+1}^{m_1+...+m_{j_1}} {\bf s}_{k}(z)
(\partial f/\partial z_{2^rk}) (u_ji_{2^rk}) ,$$  since
$i_{2^rk}w_j=w_j^*i_{2^rk}$ for each $w_j\in {\cal A}_r$ and $k\ge
1$, when $r\ge 1$, where $w_j\in {\cal A}_r$ and $u_j = u_j(z) \in
{\cal A}_v$ for each $j$, ${\bf s}_k(z)$ is the real-valued
(super-)differentiable function for each $k$. If $\xi =z_a(y)$, then
$(d\xi /dy).[(dy/d\xi ).h] = h$ for each Cayley-Dickson number $h\in
{\cal A}_v$. This implies that these two $\bf R$-linear ${\cal
A}_v$-additive operators are related by the equality $(dy/d\xi ).h =
(d\xi /dy)^{-1}.h$. On the other hand, $(dz_a(y)/dy).1= a(y) \in
{\cal A}_v$ and $y\in U\subset {\cal A}_{v}$ in the considered case.
The function $z_a(y)$ is defined up to the addendum $z_a(\alpha )$,
where $\alpha \in H_{\alpha _0}\cap U$, $H_{\alpha _0} := \{ z\in
{\cal A}_v: ~ Re (z) =Re (\alpha _0) \} $. We can choose $\phi
_a(y')$ so that $(dx/dy).(1/a) =1$ for each $y$ for which $a=a(y)\ne
0$ and inevitably we get $(dy/dx).1 = 1/a(y)$.
\par  In the particular case of $\sigma $, $\sigma _1$ and $Q$
accomplishing the differentiation with the help of the latter
identities we infer that: $$(3)\quad Q(\sigma _1^{-1}g_1)(x) = -
\sum_{j=1}^m \sum_{k=m_1+...+m_{j_1-1}+1}^{m_1+...+m_{j_1}} [(d\phi
_{f'} (x')/dx).i_{2^rk}$$  $$ + \{ {\hat g}_1(y) + (d\phi
_{g_1}(y')/dy) \} . [(dz_{a^1}(y)/dy)^{-1}.i_{2^rk}] ({\bf
s}_{k}(y)u_j(y) i_{2^rk}) ,$$  where ${\hat g}_1 = d\zeta _1(y)/dy$
for a (super-)differentiable function $\zeta _1$ such that $(d\zeta
_1(y)/dy).1= g_1(y')$, $\psi ^1_l(z) = i_l ({\bf a}_k(z) w_j
i_{2^rk}^*)$ for each $l=2^rk$ with $k\in \{ m_1+...+m_{j-1}+1,...,
m_1+...+m_j \} $, while $\psi ^1_l(z)=0$ for all others $l$ for each
$z$. Also $\psi _l(z) = i_l ({\bf a}_k(z) w_j^* i_{2^rk}^*)$ for
each $l=2^rk$ with $k\in \{ m_1+...+m_{j-1}+1,..., m_1+...+m_j \} $,
while $\psi _l(z)=0$ for all others $l$ and for each $z$. We
introduce the notation: \par $(4)\quad \theta _l(y) = i_l
[(dz_a(y)/dy)^{-1}.i_{2^rk} ] ( {\bf s}_{k}(y)u_j(y)i_{2^rk})$ for
$l=2^rk$ and $k\in \{ m_1+...+m_{j_1-1}+1,...,m_1+...+m_{j_1} \} $,
$\theta _l(y)=0$ for all others $l$ and for each $y$;
$$(5)\quad a(x) = - \sum_{j=1}^m w_j \sum_{k=m_1+...+m_{j-1}+1}^{m_1+...+m_j}
{\bf a}_k(x) = \sigma (x)\mbox{ and}$$  $(6)\quad \kappa _l(x) =
 i_l [s_k(x)u_j(x)i_{2^rk}]$ for $l=2^rk$ with $k \in \{
m_1+...+m_{j_1-1}+1,...,m_1+...+m_{j_1} \} $, $\kappa _l(x)=0$ for
every other $l$ and for each $x$.
\par In the general case $$(7)\quad \Upsilon _2(\Upsilon _1^{-1}g_1)(x)
= - \sum_{j=0}^{2^v-1} [(d\phi _{f'} (x')/dx).i_j$$  $$ + \{ {\hat
g}_1(y) + (d\phi _{g_1}(y')/dy) \} . [(dz_{a^1}(y)/dy)^{-1}.i_j]
(i_j^*\psi _j(x) ) ,$$ where ${\hat g}_1 = d\zeta _1(y)/dy$ for a
(super-)differentiable function $\zeta _1$ such that $(d\zeta
_1(y)/dy).1= g_1(y')$. We shall use the notation: \par $(8)\quad
\theta _j(y) = i_j [(dz_a(y)/dy)^{-1}.i_j] ( i_j^*\psi _j(y))$ and
for each $y$ and each $j$;
$$(9)\quad a(x) =  \Upsilon (x)\mbox{ and}$$  $(10)\quad \kappa _j(x) =
 i_j [i_j^*\psi _j(x)]$ for each $x$ and $j$.
Substituting $(7)$ into 26$(5)$ we deduce that:
\par $(11)$  $g_1(z_a(x)) + g_1(z_{\theta }(x)) = \zeta (x)$,
where $$(12)\quad \zeta (x) = - \phi _{g_1}(Im(z_a(x))) - \phi
_{g_1}(Im (z_{\theta }(x))) + \phi_{f'}(Im(z_{\kappa}(x))) +
\int_{\alpha }^x g(z)dz.$$  For $x = \alpha $ we certainly have
$z_a(\alpha ) = z_{\theta }(\alpha )$. Suppose that $a(x)\ne \theta
(x)$ identically. The dimension of the Cayley-Dickson algebra ${\cal
A}_v$ over the real field is not less, than four. Therefore, we can
choose a path $\gamma $ so that $\gamma $ is orthogonal to $\theta $
and $\kappa $ at each point on $\gamma $, that is $\gamma '(t)\perp
\theta ' (\gamma (t))$ and $\gamma '(t)\perp \kappa ' (\gamma (t))$
relative to the real-valued scalar product $(RS)$ for each $t\in
[0,1]$, where $\gamma '(t) := d\gamma (t)/dt$. Then $g_1(z_{\theta
}(x)) = g_1(\alpha )$ and $\phi _{f'}(Im (z_{\kappa}(x))) = \phi
_{f'}(\alpha ')$ for each $x=\gamma (t)$. Therefore, along such path
$\gamma $ one has
$$(13)\quad g_1(z_a(x)) + g_1(\alpha ) = \zeta (x) = - \phi _{g_1}(Im (z_a(x))
- \phi _{g_1}(\alpha ') + \phi _{f'}(\alpha ') + \int_{\alpha
}^xg(z)dz$$ for each $x=\gamma (t)$. Expressing $g_1(z)$ from
Equation $(11)$, substituting into 26$(4)$ and integrating one gets:
$$(14)\quad f(z_{a^1}(x)) = - \phi _f(Im (z_{a
^1}(x))) + \int_{\alpha }^x g_1(z)dz .$$  Particularly, if the
operator $A$ is with constant coefficients, then $s_k(x)=0$
identically for each $k$, consequently, $\theta =0$ and $\kappa =0$
identically and $g_1(z_a(x))=g_1(\alpha )= \phi _{g_1}( Im (z
_a(x)))$ for each $x$, when $f$ has a right linear derivative by
$z$. Arbitrary integration terms in $(11,14)$ can be specified from
the boundary conditions.
\par Finally, the restriction from the domain in ${\cal A}_v$
onto the initial domain of real variables in the real shadow and the
extraction of $\pi ^v _r\circ f\in {\cal A}_r$ with the help of
Formulas 1$(1-3)$ gives the reduction of a solution from ${\cal
A}_v$ to ${\cal A}_r$, where $\pi ^v_r : {\cal A}_v\to {\cal A}_r$
is the $\bf R$-linear projection operator defined as the sum of
projection operators $\pi _0+...+\pi _{2^r-1}$ given by Formulas
3$(P1,P2)$ on ${\bf R}i_j$ for $j=0,...,2^r-1$.

\par {\bf 28. Laplace's operator.}
When \par $(1)$ $A_0=\Delta _n=\sum_{j=1}^n\partial ^2/\partial
z_j^2$, \\ is Laplace's operator, then \par $(2)$ $\Upsilon f(z) =
\sum_{j=1}^n (\partial f(z) /\partial z_j) i_j^*$, so that \par
$(3)$ $\Delta _n = \Upsilon \Upsilon ^* = - \Upsilon \Upsilon $,
$\Upsilon ^1 = - \Upsilon $, \\ where $2\le n\le 2^r-1$,
$z_1,...,z_n\in {\bf R}$, in accordance with \S 2. Consider the
fundamental solution $\Psi $ of the following equation
\par $(4)$ $\Xi \Psi (z_1,...,z_n)=
\delta (z_1,...,z_n)$ with $\Xi =\Delta _n$ satisfies the identity:
\par $(5)$ $\Psi = - (\Upsilon \Psi )* (\Upsilon \Psi )$ (see the convolution
of generalized functions and this formula in \S \S 9, 19.4). \\
We seek the real fundamental solution $\Psi =\Psi _n$, since the
Laplace operator is real. The Fourier transform with the generator
${\bf i}$ (see \S 33 \cite{lumdltcdla}) by real variables
$z_1,...,z_n$ gives
\par $(6)$ $F(\Psi _n)(x) = - [F(\Upsilon \Psi _n)(x)]^2 = - [ -
{\bf i} (\sum_{j=1}^n x_ji_j^*) F(\Psi _n)(x)]^2,$ since
\par $F(\Upsilon \Psi _n)(x) = \sum_{j=1}^n F(\partial \Psi _n /\partial
z_j) i_j^* = - {\bf i} (\sum_{j=1}^n x_ji_j^*) F(\Psi _n)(x),$ \\
where $x=(x_1,...,x_n)$, $x_1,...,x_n\in {\bf R}$ (see also \S 2).
Thus we get the identity
\par $(7)$ $F(\Psi _n)(x) = - (\sum_{j=1}^nx_j^2) [F(\Psi _n)(x)]^2$ or
\par $(8)$ $F(\Psi _n)(x) = - (1/(\sum_{j=1}^n x_j^2))$ for $n\ge 3$
is the regular generalized function (functional), while
\par $(9)$ $F(\Psi _2)(x) = - {\cal P} (1/(\sum_{j=1}^2 x_j^2))$ for $n=2$. \\
We recall that the generalized function ${\cal P} (1/(\sum_{j=1}^2
x_j^2))$ on $\phi \in {\cal D}({\bf R}^2,{\bf R})$ is defined as the
regularization: $$({\cal P} (1/(\sum_{j=1}^2 x_j^2)),\phi ) =
\int_{|x|<1} [\phi (x) - \phi (0)] |x|^{-2}  dx + \int_{|x|>1} \phi
(x)|x|^{-2}dx ,$$ where $x=(x_1,x_2)$, $|x|^2 =x_1^2+x_2^2$, $x_1,
x_2\in{\bf R}$.
\par The inverse Fourier transform $(F^{-1}g)(x)=(2\pi )^{-n}(Fg)(-x)$
of the functions $1/(\sum_{j=1}^n z_j^2)$ for $n\ge 3$
and ${\cal P}(1/(\sum_{j=1}^2 z_j^2))$ for $n=2$ in the class of the
generalized functions is known (see \cite{gelshil} and \S \S 9.7 and
11.8 \cite{vladumf}) and gives
\par $(10)$ $\Psi _n(z_1,...,z_n) = C_n (\sum_{j=1}^n z_j^2)^{1-n/2}$
for $3\le n$, where $C_n = - 1/[(n-2) \sigma _n]$, $\sigma _n = 4\pi
^{n/2}/\Gamma ((n/2)-1)$ denotes the surface of the unit sphere in
${\bf R}^n$, $\Gamma (x)$ denotes Euler's gamma-function, while
\par $(11)$ $\Psi _2(z_1,z_2) = C_2 \ln (\sum_{j=1}^2 z_j^2)$ for
$n=2$, where $C_2= 1/(4\pi )$. \\ Thus the technique of convolutions
over the Cayley-Dickson algebra has permitted to get the solution of
the Laplace operator. \par Another method is with the line
integration over the Cayley-Dickson algebras. In accordance with
Formula 26$(9)$ we get
\par $\Psi _n(z_1,...,z_n) = - {\cal I}_{\Upsilon } {\cal
I}_{\Upsilon } \delta $. \par Laplace's operator and the
delta-function are invariant under any orthogonal transform $T\in
{\sf O}_n({\bf R})$ of ${\bf R}^n$. Therefore, a fundamental
solution $\Psi _n$ also is invariant relative to the orthogonal
group ${\sf O}_n({\bf R})$. That is $\Psi _n$ depends on $|z|$ and
is independent of spherical angles in spherical system of
coordinates. Thus we choose the corresponding branches of the
anti-derivative ${\cal I}_{\Upsilon } {\cal I}_{\Upsilon } \delta $.
The volume element in the Euclidean space ${\bf R}^n$ can be written
as $\lambda (dz) = x^{n-1}dx ds$, where $x=|z|$ and $ds$ is the
surface element (measure) on the unit sphere $S^{n-1}$. For each
orthogonal transform its Jacobian is unit.
\par One can take the family of test functions $\eta _{\epsilon } =
\frac{1}{(2\pi )^{n/2}\epsilon ^n} \exp \{ -
(z_1^2+...+z_n^2)/(2\epsilon ^2) \} $ tending to the delta-function,
when $\epsilon >0$ tends to zero. These functions can be written in
the $z$-representation due to Formulas 1$(1-3)$. On the other hand,
for each $z$-analytic function $\eta $ with real expansion
coefficients into a power series each line integral over the
Cayley-Dickson algebra ${\cal A}_v$ restricted on any complex plane
${\bf R}\oplus M{\bf R}$ coincides with the usual complex integral,
where $M$ is a purely imaginary Cayley-Dickson number. Therefore,
\par $\int_a^t [\int _a^y\eta (z)dz]y^kdy = \frac{1}{k+1}\int_a^t
(t^{k+1}-z^{k+1})\eta (z)dz$ for $k\ne -1$ and \par $\int_a^t [\int
_a^y\eta (z)dz]\frac{1}{y}dy = \int_a^t (\ln (t) - \ln (z))\eta
(z)dz$.
\par For the characteristic function $\chi _{B({\bf R}^n,0,x)}$ of
the ball $B({\bf R}^n,0,x)$ of radius $x>0$ with the center at zero
in the Euclidean space ${\bf R}^n$ embedded into the Cayley-Dickson
algebra ${\cal A}_v$ one can take a sequence of test functions
$\mbox{}_l\omega ^1$ converging to the regular generalized function
$\chi _{B({\bf R}^n,0,x)}$, when $l$ tends to the infinity,
consequently, $\lim_l \int_{{\bf R}^n}\mbox{}_l\omega ^1(z)\lambda
(dz) = \sigma _nx^n$. Integrating twice with the anti-derivative
operator these test functions $\eta _{\epsilon }$ in accordance with
Example 4.1 and Theorems 19 and 23 and taking the limit with
$\epsilon $ tending to zero from the right one gets Formulas
$(10,11)$.
\par This can also be deduced with the help of
the Fourier transform with the generator $\bf i$:
\par $(12)$ $F(\Psi _n)(x) =
F( - {\cal I}_{\Upsilon } {\cal I}_{\Upsilon } \delta ) =
(\sum_{j=1}^nx_j^2)^{-1} F(\delta )=(\sum_{j=1}^nx_j^2)^{-1}$. \\
Applying the inverse Fourier transform to both sides of Equation
$(12)$ we again get Formulas $(10,11)$.
\par {\bf 29. The hyperbolic operators with constant coefficients.}
\par Consider now the hyperbolic operator \par $(1)$ $A_0=L_{p,q} = \sum_{j=1}^p \partial ^2/\partial
z_j^2 - \sum_{j=p+1}^n \partial ^2/\partial z_j^2$, \\ where
$p+q=n$, $1\le p$ and $1\le q$, $(p,q)$ is the signature of this
operator, $z_1,...,z_n\in {\bf R}$.  Take two operators $\Upsilon $
and $\Upsilon _1$ with constant ${\cal A}_v$ coefficients so that
\par $(2)$ $\Upsilon f(z) = \sum_{j=1}^p (\partial f(z) /\partial
z_j) i_{2j}^* + \sum_{j=p+1}^n (\partial f(z) /\partial z_j) [ i_1^*
i_{2j}^*]$ and \par $(3)$ $\Upsilon _1 f(z) = \sum_{j=1}^p (\partial
f(z) /\partial z_j) i_{2j}^* + \sum_{j=p+1}^n (\partial f(z)
/\partial z_j) [ i_1 i_{2j}^*]$, so that
\par $(4)$ $L_{p,q} = \Upsilon \Upsilon _1$, \\
where $2\le n\le 2^{v-r}-1$, $r=1<v$, in accordance with Formulas
2$(7-9)$. Then the fundamental solution $\Psi $ of the partial
differential equation
\par $\Xi \Psi (z_1,...,z_n)= \delta (z_1,...,z_n)$ with $\Xi
=L_{p,q}$ satisfies the identity:
\par $(5)$ $\Psi = (\Upsilon ^* \Psi )* (\Upsilon _1^*\Psi )$. \\ We seek
the real fundamental solution $\Psi =\Psi _{p,q}$, since the
hyperbolic operator $L_{p,q}$ is real. Using the Fourier transform
with the generator ${\bf i}$ by real variables $z_1,...,z_n$ we
infer that
\par $(6)$ $F(\Psi _{p,q})(x) = [F(\Upsilon ^*\Psi _{p,q})(x)] [F(\Upsilon _1^*\Psi _{p,q})(x)]$  \\
$= [ - {\bf i} (\sum_{j=1}^p x_ji_{2j} + \sum_{j=p+1}^n x_j i_{2j}
i_1) F(\Psi _{p,q})(x)] [ - {\bf i} (\sum_{j=1}^p x_j i_{2j} +
\sum_{j=p+1}^n x_j i_{2j} i_1^*) F(\Psi _{p,q})(x)],$ since
\par $F(\Upsilon ^*\Psi _{p,q})(x) = \sum_{j=1}^p F(\partial \Psi _{p,q}/\partial
z_j) i_{2j} + \sum_{j=p+1}^n F(\partial \Psi _{p,q}/\partial
z_j) i_{2j} i_1$\\ $ = - {\bf i} (\sum_{j=1}^p x_j i_{2j} + \sum_{j=p+1}^n x_ji_{2j}i_1) F(\Psi _{p,q})(x)$ \\
and analogously for $\Upsilon _1^*$, where $x=(x_1,...,x_n)$,
$x_1,...,x_n\in {\bf R}$ (see also \S 2). For the function \par
$(7)$ $P(x) = \sum_{j=1}^p x_j^2 - \sum_{j=p+1}^n x_j^2$ with $p\ge
1$ and $q\ge 1$ the generalized functions $(P(x)+{\bf i}0)^{\lambda
}$ and $(P(x)-{\bf i}0)^{\lambda }$ are defined for any complex
number $\lambda \in {\bf C} = {\bf R}\oplus {\bf i}{\bf R}$ (see
Chapter 3 in \cite{gelshil}). The function $P^{\lambda }$ has the
cone surface $P(z_1,...,z_n)=0$ of zeros, so that for the correct
definition of generalized functions corresponding to $P^{\lambda }$
the generalized functions
\par \par $(8)$ $(P(x)+ c{\bf i}0)^{\lambda }=\lim_{0< c\epsilon ,
\epsilon \to 0 } (P(x)^2 + \epsilon ^2)^{\lambda /2} \exp ({\bf i}
\lambda ~ arg (P(x) + {\bf i}c\epsilon ))$ \\ with either $c=-1$ or
$c=1$ were introduced. Therefore, the identity
\par $(9)$ $F(\Psi _{p,q})(x) = - (\sum_{j=1}^px_j^2-\sum_{j=p+1}^n x_j^2) [F(\Psi _{p,q})(x)]^2$ or
\par $(10)$ $F(\Psi ) = - 1/(P(x) + c {\bf i} 0)$ follows, where
$c=-1$ or $c=1$.
\par The inverse Fourier transform in the class of the
generalized functions is: \par $(11)$ $F^{-1}((P(x)+c{\bf
i}0)^{\lambda })(z_1,...,z_n) = \exp (- \pi cq{\bf i}/2) 2^{2\lambda
+n}2^{2\lambda +n} \pi ^{n/2} \Gamma (\lambda +n/2)(Q(z_1,...,z_n) -
c{\bf i}0)^{- \lambda - n/2)}/[\Gamma (-\lambda )|D|^{1/2}]$ \\ for
each $\lambda \in {\bf C}$ and $n\ge 3$ (see \S IV.2.6
\cite{gelshil}), where $D=\det (g_{j,k})$ denotes a discriminant of
the quadratic form $P(x)=\sum_{j,k=1}^n g_{j,k}x_jx_k$, while $Q(y)=
\sum_{j,k=1}^n g^{j,k}x_jx_k$ is the dual quadratic form so that
$\sum_{k=1}^n g^{j,k}g_{k,l}=\delta ^j_l$ for all $j, l$; $\delta
^j_l=1$ for $j=l$ and $\delta ^j_l=0$ for $j\ne l$.  In the
particular case of $n=2$ the inverse Fourier transform is given by
the formula:
\par $(12)$ $F^{-1}((P(x)+c{\bf i}0)^{-1})(z_1,z_2) = - 4^{-1}|D|^{-1/2} \exp (-
\pi cq{\bf i}/2) \ln (Q(z_1,...,z_n) - c{\bf i}0).$ \par Making the
inverse Fourier transform $F^{-1}$ of the function $- 1/(P(x) + {\bf
i} 0)$ in this particular case of $\lambda =-1$ we get two complex
conjugated fundamental solutions
\par $(13)$ $\Psi _{p,q} (z_1,...,z_n) = - \exp (\pi cq{\bf i}/2)\Gamma ((n/2) -1) (P(z_1,...,z_n) +
c{\bf i}0)^{1-(n/2)}/(4\pi ^{n/2})$ for $3\le n$ and $1\le p$ and
$1\le q$ with $n=p+q$, while
\par $(14)$ $\Psi _{1,1}(z_1,z_2) = 4^{-1} \exp (\pi cq{\bf i}/2)
\ln (P(z_1,z_2) + c{\bf i}0) $ for $n=2$, where either $c= 1$ or
$c=-1$.
\par Another approach consists in using the anti-derivative
operators. The hyperbolic operator $L_{p,q}$ and the delta-function
are invariant under the Lie group ${\sf O}_{p,q}({\bf R})$ or all
linear transforms of the Euclidean space ${\bf R}^n$, $n=p+q$,
preserving the scalar product $(x,y)_{p,q}=\sum_{j=1}^px_jy_j -
\sum_{j=p+1}^{p+q} x_jy_j$ invariant. Thus $\Psi _{p,q}$ can be
written as a composition $\xi (P(x))$ of some function $\xi (y)$ and
of $P(x)$ given by Formula $(7)$. Therefore, we take the
corresponding branch of the anti-derivative in the form ${\cal
I}_{\Upsilon } {\cal I}_{\Upsilon _1} \delta = \xi (P(x))$. Applying
the Fourier transform with the generator $\bf i$ we infer that
\par $(15)$ $F(\Psi _{p,q})(x) =
F({\cal I}_{\Upsilon } {\cal I}_{\Upsilon _1} \delta ) = (P(x)+c{\bf
i}0)^{-1}F(\delta )=(P(x)+c{\bf i}0)^{-1}$. \\ Applying the inverse
Fourier transform to both sides of Equation $(15)$ one gets Formulas
$(13,14)$.
\par Thus the results of \S \S 2-25 over the Cayley-Dickson algebra ${\cal
A}_v$ lead to the fundamental solution of the hyperbolic operator
$L_{p,q}$. This means that the approach of \S \S 2-25 over the
Cayley-Dickson algebras leads to the effective solution of any
hyperbolic partial differential equation with constant coefficients.
Thus Formulas of \S \S 2, 8 with the known function $\Psi =\Psi _n$
from Formulas 28$(10,11)$ and $(13,14)$ of this section give the
fundamental solution of any first and second order linear partial
differential equation with variable $z$-differentiable ${\cal
A}_v$-valued coefficients, $z\in U\subset {\cal A}_v$.

\par {\bf 30. Example. The heat kernel.}  Each function of the type $f(z)
= P_n(z)\exp (-t |z|^2 )$ with a marked positive parameter can be
written in the $z$-representation due to Formulas 1$(1-3)$, where
$P_n(z)$ denotes the polynomial by $z$ of degree $n$. Therefore,
$f(z)$ in the $z$-representation is $z$-differentiable,
consequently, infinite $z$-differentiable (see
\cite{ludoyst,ludfov}) and $$\lim_{|z|\to \infty }
f^{(m)}(z).(h_1,...,h_m)(1+|z|^k) =0$$ for each $0\le m, k \in {\bf
Z}$ and Cayley-Dickson numbers $h_1,...,h_m\in {\cal A}_v.$
Therefore, the space ${\sf E}$ of infinite $z$-differentiable
tending to zero at infinity functions together with their
derivatives multiplied on the weight factor $(1+|z|^k)$ is
infinite-dimensional. Thus it is worthwhile to consider the
topologically adjoint space ${{\sf E}'}_q$ of ${\bf R}$-linear
${\cal A}_v$-additive continuous ${\cal A}_v$-valued functionals on
${\sf E}$. Elements of ${{\sf E}'}_q$ are also called the
generalized functions. A function or a generalized function is
called finite if its support is a bounded set.
\par The heat partial differential equation reads as
\par $(1)$ $\partial v(z)/\partial z_0 = a^2 \Delta v(z) + f(z),$ \\
where $z=z_0i_0+...+z_mi_m$, $z_0,...,z_m\in {\bf R}$, $1\le  m \le
2^v-1$, $2\le v$, where $a>0$, $f(z)$ is a real-valued generalized
finite function so that $f(z)$ is zero for $z_0<0$ (see \S 16
\cite{vladumf}). We shall seek the generalized solution $\cal E$ of
this equation using the technique given above. The generalized
function $v= {\cal E} * f$ is the solution of $(1)$, where
\par $(2)$ $\partial {\cal E}(z)/\partial z_0 - a^2 \Delta {\cal E}(z) = \delta (z)$,
$$(3)\quad ({\cal E}*f)(x) = \int_0^{x_0}\int_{\bf R^m} {\cal
E}(x-z) f(z) dz_0....dz_m.$$ As usually $\delta $ denotes the
$\delta $ generalized function so that $$(4)\quad ({\delta }
*f)(x) = \int_0^{x_0}\int_{\bf R^m} {\delta }(x-z) f(z)
dz_0....dz_m= f(x)$$   for each continuous bounded function $f$.  If
$f$ is (super-)differentiable and bounded in each domain $\{ z: 0\le
z_0 \le T \} $ for $0<T<\infty $, $f(z)=0$ for $z_0<0$, then the
solution $v$ is also (super-)differentiable in the domain $z_0>0$ as
it will be seen from the formulas given below. Let us seek the
generalized solution $\cal E$ in the form ${\cal E} (z) = w(z_0)
e^{u(z)}$, where $w$ and $u$ are unknown real-valued functions to be
determined. Calculating derivatives of ${\cal E}$ and substituting
into Equation $(2)$ one gets:
\par $(5)$ $e^{u(z)} \{ w'(z_0) + w(z_0) \partial u(z)/\partial z_0
\} - a^2 e^{u(z)} w(z_0) \sum_{j=1}^m [ (\partial u(z)/\partial
z_j)^2 + \partial ^2 u(z)/\partial z_j^2 ] = \delta (z) ,$
consequently,
\par $(6)$ $(1/w(z_0)) w'(z_0) = - \partial u(z)/\partial z_0
+ a^2 \sum_{j=1}^m [ (\partial u(z)/\partial z_j)^2 +
\partial ^2 u(z)/\partial z_j^2 ] + e^{-u(z)} (1/w(z_0)) \delta (z)
.$ \\ Take now any sequence of continuous non-negative functions
$\eta _n$ with compact supports $U_n$ such that $U_{n+1}\subset U_n$
for each $n$, with $\bigcap_n U_n = \{ 0 \} ,$  $$(7)\quad
\int_{{\bf R}^{m+1}} \eta _n(z) dz_0 ... dz_m = 1$$ for all $n$,
tending to $\delta $ on the space of continuous functions $p(z)$ on
${\bf R}i_0\oplus ... \oplus {\bf R} i_m$ with the converging
integral $\int_{{\bf R}^{m+1}} |p(z)|^2 dz_0...dz_m<\infty $:
\par $(8)$ $\lim_{n\to \infty } \int_{{\bf R}^{m+1}} p(z) \eta _n(z)
dz_0 ... dz_m = p(0) .$ \par Therefore, we get that on ${\bf
R}i_0\oplus ... \oplus {\bf R} i_m\setminus \{ 0 \} $ for $z_0>0$
the following equation $$(9)\quad (1/w(z_0)) w'(z_0) = - \partial
u(z)/\partial z_0 + a^2 \sum_{j=1}^m [ (\partial u(z)/\partial
z_j)^2 +
\partial ^2 u(z)/\partial z_j^2 ]$$   need to be satisfied.
The left side of $(9)$ is independent of $z-z_0$, hence the right
side is also independent of $z-z_0$. The partial differential
operator
\par $\{ \partial u(z)/\partial z_0 + a^2 \sum_{j=1}^m [ (\partial
u/\partial z_j)^2 + \partial ^2 u(z)/\partial z_j^2 ] \} $ acting on
$u$ is of the second order.  For each Cayley-Dickson number $z\in
{\cal A}_r$ the identity $z^2 = 2z ~ Re (z) - |z|^2$ is satisfied,
particularly, $M^2 = - |M|^2$ for each purely imaginary number $M\in
{\cal A}_r$.  Therefore, a function $u$ may be only a polynomial by
real variables $z_1,...,z_m$ of degree not higher than two. On the
other hand, the Laplace operator $\Delta $ and the $\delta $
function are invariant relative to all elements $C$ of the
orthogonal group ${\sf O}_m({\bf R})$ acting on variables
$z_1,...,z_m$. Each ${\sf O}_m({\bf R})$ invariant real polynomial
$P$ of the second order has the form $\alpha (z_1^2+...+z_m^2) +
\beta ,$ where $\alpha $ and $\beta $ are two constants independent
of $z_1,...,z_m$. Thus $u$ as the polynomial of $z_1,...,z_m$ may
depend on $|z-z_0|^2 = z_1^2+...+z_m^2$ only. The latter sum of
squares can be written in the $z$-representation with the help of
Formulas 1$(1-3)$. This means that ${\cal E}$ has the form: \par
$(10)$ ${\cal E} (z) =w(z_0) \exp \{ \alpha (z_0) (z_1^2+...+z_m^2)
+ \beta (z_0) \} $ and Equation $(9)$ simplifies:
\par $(11)$ $(1/w(z_0)) w'(z_0) = - [d\alpha (z_0)/d
z_0](z_1^2+...+z_m^2) - [d\beta (z_0)/dz_0]$ $$ + a^2 \alpha (z_0)
\{ 2m + \alpha (z_0) \sum_{j=1}^m 4z_j^2 \} .$$ We  can denote
$w(z_0)e^{\beta (z_0)}$ by $w(z_0)$ again and take without loss of
generality that $\beta =0$. The left side of $(11)$ is independent
of $z_1,...,z_m$, hence terms with $|z-z_0|^2$ in $(11)$ are
canceling: $\alpha ^{-2}(z_0) [d\alpha (z_0)/dz_0] = 4 a^2$. The
latter differential equation gives $\alpha (z_0) =  -
1/(c_0+4a^2z_0)$, where $c_0$ is the real constant. Substituting
this $\alpha $ into $(11)$ one gets:
\par $(12)$ $(1/w(z_0)) w'(z_0) = a^2 \alpha (z_0) 2m   $. \\ Together with
Condition $(2)$ this gives $C_0=0$ and the heat kernel ${\cal E}$:
\par $(13)$ ${\cal E} (z) = \theta (z_0) [2a(\pi z_0)^{1/2}]^{-m}
\exp \{ - |z-z_0|^2 /[4a^2 z_0] \} $ \\ and the solution \par $(14)$
$v= {\cal E} * f$, \\ where $\theta (z_0)=1$ for $z_0\ge 0$ and
$\theta (z_0)=0$ for $z_0 <0$. \par If use anti-derivation operators
the solution has the form 26$(6-8)$ supposing that a solution ${\cal
E}$ is real-valued on real-valued test functions $\omega $, $[{\cal
E},\omega )\in {\bf R}$, where $\Upsilon _1=\Upsilon $, $\Upsilon
\Upsilon = - a^2 \Delta _m$ (see \S ) and $\Upsilon _2=\partial
/\partial z_0$. Therefore, \par $(14)$ $ a^2{\cal E} = - {\cal
I}_{\Upsilon } {\cal I}_{\Upsilon } (\partial {\cal E}/\partial z_0
- \delta )$. \\ Making the Fourier transform $F=F_{z_1,...,z_m}$ by
the variables $z_1,...,z_m$ with the generator ${\bf i}$ of both
sides of Equation $(14)$ one gets for suitable branches of the
anti-derivatives
\par $(15)$ $a^2 F({\cal E})(z_0,x_1,...,x_m) = [a^2\sum_{j=1}^m
x_j^2]^{-1}(\partial F({\cal E})/\partial z_0 - \delta (z_0))$. \\
Solving the latter ordinary differential equation one finds $F({\cal
E})(z_0,x_1,...,x_m)$ and making the inverse Fourier transform by
the variables $x_1,...,x_m\in {\bf R}$ one gets Formula $(13)$.
\par {\bf 31. Example. Wave operator.} In this section the
fundamental solution ${\cal E}={\cal E}_n$ of the wave operator is
considered:
\par $(1)$ ${\overline {\coprod }}  {\cal E}(t,x) = \delta (t,x)$, where
\par $(2)$ ${\overline {\coprod }} f = \partial ^2f/\partial t^2 -  \Delta f$
denotes the wave (d'Alambert) operator with
$$(3)\quad \Delta f = \sum_{j=1}^n \partial ^2f/\partial x_j^2,$$
where $t\ge 0$, $x_1,...,x_n\in {\bf R}$. We make the change of
variables putting $t=z_2$, $x_j=z_{2j+2}$ for each $j=1,...,n$,
$z=z_0i_0+z_2i_2+...+z_{2^v-2}i_{2^v-2} \in {\cal A}_{1,v}$,
$z_0,...,z_{2^v-1}\in {\bf R}$, $r=1$. We consider the case $n=3$
and $v=4$ so that ${\cal A}_{1,4}$ is isomorphic with the octonion
algebra ${\cal A}_3={\bf O}$. Let us seek ${\cal E}$ in the class of
the generalized functions in the form ${\cal E}(z) = \theta (z_2)
f(z)$, where $\theta $ and $f$ are some generalized functions to be
calculated, $f$ may depend only on $z_2,z_4,...,z_{2n+2}$.
D'Alambert's operator ${\overline {\coprod }} $ is invariant
relative to any $\bf R$-linear transformations $A$ from the Lie
group ${\sf O}_{1,n}({\bf R})$. Elements of the group ${\sf
O}_{1,n}({\bf R})$ are characterized by the condition $A^tGA=G$,
where $G$ denotes the square $(n+1)\times (n+1)$ diagonal matrix
$G=diag (1,-1,...,-1)$, the transposed matrix $A$ is denoted by
$A^t$. This means that the wave operator ${\overline {\coprod }} $
is invariant under change of variables $\xi = (z_2,z_4,...,z_{2n+2})
A$ for any $A\in {\sf O}_{1,n}({\bf R})$. Making the differentiation
of ${\cal E}$ one gets the differential equation:
\par $(4)$ ${\overline {\coprod }} {\cal E} = (\partial ^2\theta /\partial z_2^2)f+
2(\partial \theta /\partial z_2) (\partial f/\partial z_2)+ \theta
{\overline {\coprod }} f= \delta (z)$. \par The $\delta $-function
$\delta (z)$ is also invariant relative to all transformations of
the Lie group ${\sf O}_{1,n}({\bf R})$, since $\delta g=g(0)$ for
each continuous function $g$ with $\int_{{\bf R}^{n+1}} |g(z)|^2
dz_2 ... dz_{2n+2} <\infty $. On the other hand, Equation $(4)$
implies that
\par $(5)$ $\partial ^2 \theta /\partial z_2^2 = - [ 2(\partial \theta /\partial z_2)
(\partial f/\partial z_2) + \theta {\overline {\coprod }} f - \delta (z) ]/f(z)$ \\
for each $z\in {\cal A}_{r,v}$, when $f(z)\ne 0$. The left side of
Equation $(5)$ may depend only on $z_2$, consequently, the right
side of $(5)$ is independent of $z_4, ..., z_{2n+2}$.  In view of
Formulas 29$(2,3)$ with $p=1$ and $q=n$ we get the operators $\sigma
= \Upsilon $ and $\sigma _1=\Upsilon _1$ with $Q=0$ up to the
enumeration of the variables. Therefore, one gets the functions
$\Psi _{1,n}$ (see Formulas 29$(13,14)$) over the Cayley-Dickson
algebra ${\cal A}_v$. But due to the ${\sf O}_{1,n}({\bf R})$
invariance of the generalized function ${\cal E}$ we infer that it
is necessary to take the ${\sf O}_{1,n}({\bf R})$ invariant
polynomial $P(y)= (y_2^2 - \sum_{j=1}^n y_{2j+2}^2)$. Thus we put
${\cal E} = \theta (z_2) f(z)$ with $f(z)=u(z_2^2 - \sum_{j=1}^n
z_{2j+2}^2)$, where $u$ is some generalized function. Substituting
$u$ instead of $f$ into $(5)$ one gets the simplified differential
equation. If suppose that $\partial \theta /\partial z_2=0$ for
$z_2>0$, then $\partial ^2\theta /\partial z_2^2=0$ and Equation
$(5)$ leads to the differential equation
\par $(7)$ $4 u^{(2)}.(1,1) (\eta ) \eta  - 4 u'.1 (\eta ) =\delta (z)/c $, \\
where $\eta = \eta (z) =z_2^2-z_4^2-z_6^2-z_8^2$, $\theta (z_2) = c
= const$ for $z_2>0$. Choose any sequence of $z$-differentiable
functions $g_n(z)$ with compact supports converging to the $\delta
$-function as in \S \S 24 and 30, when $n$ tends to the infinity.
Each function $g_n(z)/\eta (z)$ has poles of the first order at
points $z_2=[z_4^2+z_6^2+z_8^2]^{1/2}$ and $z_2= -
[z_4^2+z_6^2+z_8^2]^{1/2}$. Making the substitution $p=u'.1$ in
$(7)$ and Formula 3$(10)$ \cite{ludeoc} with the right side $Q(z) =
g_n(z) /\eta (z)$ and $b(z) = - 1/\eta (z)$ we obtain the integral
expression for the solution $p_n$ of the differential equation
\par $(8)$ $ {p'}_n.1 (\eta ) - p_n(\eta )/\eta  = g_n(z)/(4c\eta
)$. \\ To evaluate the appearing integrals it is possible to use
Jordan's Lemma (see \S 2 in \cite{lutsltjms}) over the octonion
algebra isomorphic with ${\cal A}_{1,4}$. The evaluation of the
integrals (see \S 3 also) with the given functions can be reduced to
the complex case, when $\alpha $ and $x$ belong to the same complex
plane ${\bf C}_M$. Applying Jordan's lemma one deduces the
expression for $p_n$ and the limit function $p(\eta ) = \delta
'(\eta )/(2\pi c) + K$, where $K$ is a constant, since $\eta \in
{\bf R}$. Therefore, $u(\eta ) = \int p(\eta )d\eta $. Thus one
infers the fundamental solution
\par $(9)$ ${\cal E}_3(z) = \theta (t) \delta (t^2-|x|^2) /(2\pi )$
\\ and the generalized solution ${\cal E}_3*s$ of the wave equation
\par ${\overline {\coprod }} f=s$, where $s=s(z)$ is a generalized function or particularly a
$z$-differentiable function. The delta generalized function $\delta
(P)$ of the quadratic form
$P(x)=x_1^2+...+x_p^2-x_{p+1}^2-...-x_{p+q}^2$ is described in
details in \S IV.2 \cite{gelshil}.

\par {\bf 32. Helmholtz' operator.}
\par When $\beta \ne 0$ with $Re (t_j^*\beta )=0$ for each $j$ (see \S 8),
for example, $\beta =\beta _0+ i_k\beta _1$ with real $\beta _0$ and
$\beta _1$ and $k>n$, then
\par $(1)$ $A_0= \Delta _n + |\beta |^2$ is Helmholtz' operator.
The corresponding operator $\Upsilon $ is given by Formula 28$(2)$.
\par For an arbitrary real non-degenerate quadratic form $P(x)$
generalized functions $(c^2+P + b{\bf i}0)^{\lambda }$ with $c>0$,
$b=1$ or $b=-1$, $\lambda \in {\bf C} = {\bf R}\oplus {\bf R}{\bf
i}$, are defined as:
\par $(2)$ $(c^2+P + b{\bf i}0)^{\lambda } =\lim_{0<\epsilon \to 0}
(c^2+ P + b \epsilon {\bf i}P_1)^{\lambda },$ \\ where $P_1$ is a
positive definite quadratic form.
\par Some special functions are useful for such equations.
Bessel's functions are solutions of the differential equation
\par $(S1)$ $z^2d^2 u/dz^2+zdu/dz + (z^2-\lambda ^2)u=0$, \\ where $\lambda $
and $z$ are complex. Bessel's function of the first kind is given by
the series:
$$(S2)\quad J_{\lambda }(z)=\sum_{m=0}^{\infty
}(-1)^m(z/2)^{2m+\lambda }/[m!\Gamma (\lambda +m+1)],$$ where $z$
and $\lambda \in {\bf C}_{\bf i}$. Then \par $(S3)$ $I_{\lambda }(z)
= \exp (-\pi \lambda {\bf i}/2) J_{\lambda }({\bf i}z)$ \\
is called Bessel's function of the imaginary argument. Other needed
functions for non-integer $\lambda $ are:
\par $(S4)$ $N_{\lambda }(z)=[J_{\lambda }(z) \cos (\pi \lambda ) -
J_{-\lambda }(z)]/\sin (\pi z)$,
\par $(S5)$ $H^{(1)}(z) = J_{\lambda }(z) + {\bf i} N_{\lambda }(z)$,
\par $(S6)$ $H^{(2)}(z) = J_{\lambda }(z) - {\bf i} N_{\lambda }(z)$,
\par $(S7)$ $K_{\lambda }(z) = \pi [I_{-\lambda }(z) - I_{\lambda
}(z)]/[2\sin (\pi \lambda )]= \pi \exp ({\bf i}\pi (\lambda +1)/2) H^{(1)}_{\lambda }({\bf i}z)/2$ \\
with the complex variable $z$ and non-integer complex parameter
$\lambda $. For integer $\lambda $ values of these functions
$(S4-7)$ are defined as limits by $\lambda \in {\bf C}_{\bf
i}\setminus {\bf Z}$. The functions $H^{(1)}$ and $H^{(2)}$ are also
solutions of Bessel's differential equation $(S1)$ and they are
called Hankel's functions of the first and the second kind
respectively, $K_{\lambda }(z)$ is known as Mcdonald's function. The
functions $I_{\lambda }(z)$ and $K_{\lambda }(z)$ are linearly
independent solutions of the differential equation:
\par $(S8)$ $z^2d^2 u/dz^2 + zdu/dz - (z^2 + \lambda ^2)u=0$ \\
(see about special functions in \cite{lavrsch,nikuvar}).
\par Analogously to \S 28 using Formulas 19.4$(C4-C7)$ and 14$(3,4)$
or Theorem 23 for a fundamental solution $\Psi _n$ of the equation
\par $(3)$ $A_0\Psi _n=\delta $, \\ where $A_0$ is Helmholtz'
operator, we get the identity
\par $(4)$ $F(\Psi _n)(x) = [ c^2 - (\sum_{j=1}^nx_j^2)] [F(\Psi _n)(x)]^2$ or
\par $(5)$ $F(\Psi _n)(x) = 1/[c^2 - (\sum_{j=1}^n x_j^2) + b {\bf i}0]$, \\
where $c>0$, $c=|\beta |$. The Fourier transform of the generalized
function $(c^2+P(x) + b{\bf i}0)^{\lambda }$ by the real variables
$x=(x_1,...,x_n)$ with the generator ${\bf i}$ is:
\par $(6)$ $F[(c^2+P(x) + b{\bf i}0)^{\lambda }](y) =$ \\ $ 2^{\lambda +1} (2\pi
)^{n/2} c^{\lambda +(n/2)} K_{\lambda +(n/2)} [c (Q(y) -b{\bf
i}0)]/[
\Gamma (-\lambda ) D^{1/2} (Q(y)-b{\bf i}0)^{(\lambda /2 + n/4)}],$ \\
where $D=\det (g_{j,k})$ denotes a discriminant of the quadratic
form $P(x)=\sum_{j,k=1}^n g_{j,k}x_jx_k$, $Q(y)= \sum_{j,k=1}^n
g^{j,k}x_jx_k$ is the dual quadratic form so that $\sum_{k=1}^n
g^{j,k}g_{k,l}=\delta ^j_l$, $\delta ^j_l=1$ for $j=l$ and $\delta
^j_l=0$ for $j\ne l$ (see \S IV.8.2 \cite{gelshil}). Mention that
$D^{1/2}=|D|^{1/2} \exp (-q\pi {\bf i}/2)$ if the canonical
representation of the quadratic form $P$ has $q$ negative terms.
\par Another formula is:
\par $(7)$ $F[(c^2+P(x) + b{\bf i}0)^{\lambda }](y) =$ \\ $ 2^{\lambda +1} (2\pi
)^{n/2} \exp ( - bq\pi {\bf i}/2) c^{\lambda +(n/2)} \{ K_{\lambda
+(n/2)} [c (Q_+(y))^{1/2}]/[ \Gamma (-\lambda ) |D|^{1/2}
(Q_+(y))^{\lambda /2+n/4}]$ \\ $ + (b \pi {\bf i}/2) H^{(j(b))}_{-
\lambda - (n/2)} [c (Q_-(y))^{\lambda /2+n/4}]/[
\Gamma (-\lambda ) |D|^{1/2} (Q_-(y))^{\lambda /2+n/4}] \} ,$ \\
where $j(1)=1$, $j(-1)=2$, \par $(P_+^{\lambda },\phi ) =\int_{P>0}
P^{\lambda }\phi dx_1...dx_n,$ \par $(P_-^{\lambda },\phi )
=\int_{P<0} P^{\lambda }\phi dx_1...dx_n$. \\ The functions
$(P+b{\bf i}0)^{\lambda }$ and $P_+^{\lambda }$ and $P_-^{\lambda }$
by the complex variable $\lambda $ are regularized as it is
described in \cite{gelshil} with the help of their Laurent series in
neighborhoods of singular isolated points $\lambda $ such that after
the regularization only the regular part of the Laurent series
remains. The functions $(P+b{\bf i}0)^{\lambda }$ with $b=1$ or
$b=-1$ have only simple poles at the points $\lambda = -n/2$,
-(n/2)-1,...,$-(n/2)-k$,..., where $k=1,2,...$ is a natural number.
Using formula $(6)$ with $\lambda = -1$ and $P(x)= -
(x_1^2+...+x_n^2)$ one gets the fundamental solution $\Psi _n$,
where $(F^{-1}g)(x)= (2\pi )^{-n} (Fg)(-x)$. Particularly, $\Psi
_3(x) = - \exp (bc{\bf i} |x|)/(4\pi |x|)$, $\Psi _2(x) = - {\bf i}
H_0^{(1)}(c|x|)/4$ or its complex conjugate ${\bf i}
H_0^{(2)}(c|x|)/4$, where $H_0^{(j)}$ denotes Hankel's function,
$j=1, 2$.

\par {\bf 33. Klein-Gordon's operator.}
\par Consider $\beta $ and $t_j$ as in \S 8 with $Re (t_j^*\beta
)=0$  for each $j$, $c=|\beta |>0$.  Take the operator
\par $(1)$ $A_0=L_{p,q} +c^2$, \\ where $L_{p,q}$ is the hyperbolic
operator as in \S 29. For $p=1$ and $q=3$ the operator $A_0$ is
called Klein-Gordon's operator. From Formulas 32$(C4-C7)$ and
14$(3,4)$ or Theorem 23 we infer that
\par $(2)$ $F(\Psi _n)(x) = [ c^2 - (\sum_{j=1}^px_j^2-\sum_{j=p+1}^nx_j^2)] [F(\Psi _n)(x)]^2$ or
\par $(3)$ $F(\Psi _n)(x) = 1/[c^2 - (\sum_{j=1}^p x_j^2-\sum_{j=p+1}^nx_j^2) + b {\bf i}0]$. \\
Then Formulas 19.4$(6)$ or 32$(7)$ with $\lambda =-1$ and $P(x) = -
x_1^2-...-x_p^2 + x_{p+1}^2+...+x_n^2$, $n=p+q$, give the
fundamental solution $\Psi _{p,q}$ of the equation
\par $(4)$ $(L_{p,q}+c^2)\Psi _{p,q} =\delta $, \\
where $(F^{-1}g)(x)= (2\pi )^{-n} (Fg)(-x)$. There are two $\bf
R$-linearly independent fundamental solutions, so their $\bf
R$-linear combination with real coefficients $\alpha _1$ and $\alpha
_2$ such that $\alpha _1+\alpha _2=1$ is also a fundamental
solution.

\par {\bf 34. Remark.} Certainly, more general partial differential equations as 30$(1)$,
but with $\partial ^lv/\partial z_0^l$, $l\ge 2$, instead of
$\partial v/\partial z_0$ can be considered. It is worth to mention,
that alternative deductions using Formulas 7$(1)$ and 27$(11,14)$
can be used instead of 8$(1)$ and 19.4$(C1-C7)$ in \S \S 30 and 31
providing also $u(z) = \alpha (z_1^2+...+z_m^2) + \beta $ and
$f(z)=u(z_2^2 - \sum_{j=1}^n z_{2j+2}^2)$ with the help of the
symmetry Lie groups ${\sf O}_m({\bf R})$ and ${\sf O}_{1,n}({\bf
R})$. Indeed, Functions $P(x)^{\lambda }$ satisfy Condition 7$(P3)$
for any real $\lambda $, where $P(z) =
z_1^2+...+z_p^2-z_{p+1}^2-...-z_n^2$, $1\le p \le n\le 2^{v-1}-1$,
$1\le v$, since $(dP(z)^{\lambda }/dz).h = P(z)^{\lambda -1} \lambda
2 Re (\eta (z)h)$, where $z\in {\cal A}_v$, $\eta (z) =
z_1i_1+...+z_ni_n$ for $p=n$, while $\eta (z)= z_1i_2+...+z_pi_{2p}
+ z_{p+1} i_1^* i_{2(p+1)}+...+z_ni_1^*i_{2n}$ for $p<n$. The
function $\eta (z)$ can be written in the $z$-representation due to
Formulas 1$(1-3)$.
\par Formally the case of the hyperbolic operators
$L_{p,q}+c^2$ and their fundamental solutions is obtained from the
elliptic operators $\Delta _n +c^2$ with $c\ge 0$ by the change of
variables $(x_1,...,x_n)\mapsto (x_1,...,x_p,x_{p+1}{\bf
i},...,x_{p+q}{\bf i})$, where $n=p+q$, $x_j\in {\bf R}$ and
$x_j{\bf i}\in {\bf C}={\bf R}\oplus {\bf R}{\bf i}$ for each $j$,
since the quadratic forms $P$ may be with complex coefficients and
their Fourier transforms can be considered as in \cite{gelshil}. At
the same operators $\sigma $ or $\Upsilon $ for these particular
operators $L_{p,q}+c^2$ and $\Delta _n+c^2$ can be written over the
complexified algebra $({\cal A}_r)_{\bf C}$ instead of the
Cayley-Dickson algebra ${\cal A}_v$, $2\le r<v$ (see \S 2 above).
For this we take in Formula 2$(8)$ ${\bf i}$ instead of $w_j^*$ so
that \par $(1)$ $\sigma f(z) = \sum_{j=1}^p (\partial f(z)/\partial
z_j) i_j^* + \sum_{j=p+1}^n (\partial f(z)/\partial z_j) {\bf i}
i_j^*$, consequently,
\par $(2)$ $ (c +\sigma ) (c - \sigma ) f = (L_{p,q}+c^2)f$ and
\par $(3)$ $(c+\sigma ) (c+ \sigma )^* f = (c+\sigma )^* (c+\sigma ) f = (\Delta _n +c^2)f$.
\par Let \par $(4)$ $\Xi _{c,p,q} = L_{p,q} +c^2$, where $c\in {\bf R}$, $L_{n,0}=\Delta
_n$, $1\le p\le n$, $q=n-p$, and let $\Psi _{c,p,q}$ be a
fundamental solution of the equation
\par  $(5)$ $\Xi _{c,p,q}\Psi _{c,p,q} = \delta $. \\ Then due to
Identities $(2,3)$ a fundamental solution ${\cal E} = {\cal
E}_{\beta + \sigma }$ of the equation
\par $(6)$ $(\sigma + \beta ) {\cal E} = \delta $ can be written in the form:
\par $(7)$ ${\cal E}_{\beta +\sigma } = (\sigma + \beta )^* \Psi _{c,n,0}$,
where $\beta \in {\cal A}_r$, $|\beta |=c$, $Re (t_j^*\beta )=0$ for
each $j$, $t_j=i_j^*$ for $1\le j\le p$, $t_j={\bf i}i_j^*$ for
$p<j\le n$. Moreover,
\par $(8)$ ${\cal E}_{c +\sigma } = (c - \sigma ) \Psi _{c,p,q}$.
\par Therefore, we infer a solution of the equation
\par $(9)$ $(\sigma + \beta )f = g $ in ${\cal D}({\bf R}^m,{\cal A}_r)$ or in
the space ${\cal D}({\bf R}^m,{\cal A}_r)^*_l$: \par $(10)$  $f =
{\cal E}_{\beta +\sigma }*g$.  From $(2,3)$ we deduce a fundamental
solution $\cal V$ of the equation
\par $(10)$ $A_0{\cal V} = (\sigma + \beta ) (\sigma _1 + \beta _1) {\cal V}= \delta $
in the convolution form: \par $(11)$
${\cal V} = {\cal E}_{\sigma +\beta } * {\cal E}_{\sigma _1+\beta _1},$ \\
since
\par $(12)$ $A_0 {\cal V} = (\sigma +\beta ) (\sigma _1+\beta _1) (((\sigma +\beta )^*
\Psi _{c,n,0})* ((\sigma _1+\beta _1)^*\Psi _{c,n,0})$
\par $ = (((\sigma +\beta ) (\sigma + \beta _1)^* \Psi _{c,n,0})* (((\sigma _1
+\beta _1) (\sigma _1 + \beta _1)^* \Psi _{c,n,0})) = \delta *\delta
=\delta $. Particularly, \par $\Psi _{c,p,q} = - {\cal E}_{c+\sigma
}*{\cal E}_{- c + \sigma } = ((c-\sigma )\Psi _{c,n,0})*((c+\sigma
)\Psi _{c,n,0})$, that can be lightly verified after the Fourier
transform by real variables with the generator ${\bf i}$, since by
Formulas $(1,2)$ the operator $\sigma $ and its anti-derivative
operator ${\cal I}_{\sigma }$ correspond to the signature $(p,q)$
and $F(\sigma \sigma ^* g) = - |z|^2 F(g)$ for any $g\in {\cal
D}({\bf R}^n,{\cal A}_v)^*_l$.
\par Knowing a fundamental solution it is possible to consider then
a boundary problem (see also \cite{hoermpd,vladumf}).

\par {\bf 35. Partial differential equations of order higher than
two.}  \par The fundamental solution $\Psi _{\Upsilon ^m +\beta }$
of the equation \par $(1)$ $\Xi _{2m, \beta }\Psi _{\Upsilon ^m
+\beta } = \delta ,$ \\ where \par $(2)$ $\Xi _{2m,\beta } =
(\Upsilon ^m + \beta )(\Upsilon ^m + \beta )^*$ \\  can be obtained
using decompositions with the help of operators of the first order
$\Upsilon _k + \beta _k $ by induction, if such decomposition exists
(see \S \S 10-14 above). Suppose that this decomposition is
constructed
\par $(3)$ $(\Upsilon ^m+ \beta ) f(z) = (\Upsilon _m + \beta
_m)[...[(\Upsilon _2+\beta _2) [(\Upsilon _m+\beta _m )f(z)]]...]$,
then the fundamental solution can be written as the iterated
convolution \par $(4)$ $\Psi _{\Upsilon ^m+\beta } = [...[[(\Upsilon
_m+\beta _m)^*{\cal E}_m]*[(\Upsilon _{m-1}+\beta _{m-1})^*{\cal
E}_{m-1}]]*...]*[(\Upsilon _1+\beta _1)^*{\cal E}_1] ,$
\\ where ${\cal E}_j$ denotes the fundamental solution of the elliptic second order
partial differential equation \par $(5)$ $A_j {\cal E}_j = \delta $,
with \par $(6)$ $A_j=(\Upsilon _j+\beta _j)(\Upsilon _j+\beta
_j)^*$.
\\ The fundamental solutions ${\cal E}_j$ were written above in \S
\S 2 and 28. Indeed, using Equalities 4$(7-9)$ by induction we have
\par $(7)$ $\sum_{s} (...(a_{m,k_m}^*
a_{m-1,k_{m-1}}^*)a_{m-2,k_{m-2}}^*)...)a_{1,k_1}^*)a_{1,l_1})a_{2,l_2})...)a_{m,l_m}$
\\  $= Re (a_{m,k_m}^*a_{m,l_m}) ... Re (a_{1,k_1}^*a_{1,l_1})$, \\ where
$\sum_s$ denotes the sum by all pairwise transpositions
$(k_1,l_1)$,...,$(k_m,l_m)$, $ ~ a_{j,k}\in {\cal A}_v$. Therefore,
\par $(8)$ $\Xi _{2m, \beta }\Psi _{\Upsilon ^m+\beta } =$ \\ $
[...[[(\Upsilon _m+\beta _m)(\Upsilon _m+\beta _m)^*{\cal
E}_m]*[(\Upsilon _{m-1}+\beta _{m-1})(\Upsilon _{m-1}+\beta
_{m-1})^*{\cal E}_{m-1}]]*...]*[(\Upsilon _1+\beta _1)(\Upsilon
_1+\beta _1)^*{\cal E}_1]$ $ = [...[\delta * \delta ]*...]*\delta =
\delta .$ \par Vice versa if the fundamental solution $\Psi
_{\Upsilon ^m + \beta }$ is known, then we get the fundamental
solution ${\cal E}^m_{\beta }$ of the equation
\par $(9)$ $(\Upsilon ^m + \beta ) {\cal E}^m_{\beta } = \delta $
as \par $(10)$ ${\cal E}^m_{\beta } = (\Upsilon ^m + \beta )^* \Psi
_{\Upsilon ^m+\beta }$ in accordance with $(2,7)$. Moreover, the
equation
\par $(11)$ $A_{m+k} f = g $ with $A_{m+k} =
(\Upsilon ^m_1 + \beta _1) (\Upsilon ^k_2 + \beta _2)$
\\ in ${\cal D}({\bf R}^n,{\cal A}_v)$ or in
the space ${\cal D}({\bf R}^n,{\cal A}_v)^*_l$, where $n$ is a
number of real variables, $2\le n\le 2^v$, has the fundamental
solution ${\cal V}_{m+k}$:
\par $(12)$ ${\cal V}_{m+k} = {\cal E}_{\Upsilon ^m_1 + \beta _1} *
{\cal E}_{\Upsilon ^k_2 + \beta _2}$, where \par $(13)$ ${\cal
E}_{\Upsilon ^m_1 + \beta _1} = (\Upsilon ^m_1 + \beta _1)^*\Psi
_{\Upsilon ^m_1+\beta _1}$ \\ denotes the fundamental solution of
the equation
\par $(14)$ $(\Upsilon ^m_1 + \beta _1){\cal
E}_{\Upsilon ^m_1 + \beta _1} = \delta $, consequently, \par $(15)$
$f= {\cal V}_{m+k} *g$ is the solution of Equation $(11)$.
\par For example, the fourth order partial differential operator
\par $A_4 f(z)=\sum_{j=1}^p \partial ^4f(z)/\partial z_j^4 -
\sum_{j=p+1}^n
\partial ^4f(z)/\partial z_j^4$ \\ can be decomposed as the composition
of two operators of the second order $\Upsilon ^2$ and $\Upsilon
^2_1$ formally as $\sigma $ and $\sigma _1$ in 2$(8,9)$ with the
substitution of $\partial f/\partial z_{2^rj}$ on $\partial ^2
f/\partial z_{2^rj}^2$ so that in accordance with Theorem 10 this
operator $A_4$ can be presented in the form given by Formulas
$(2,3,11)$.
\par On the other hand, fundamental solutions of $\Delta _n^k$ and
$L_{p,q}^k$ and $A_2^k$ for certain other second order partial
differential operators are known. So combining them with operators
of the form $\Upsilon ^{m_1}_1...\Upsilon ^{m_k}_k$ permits to
consider fundamental solutions of many partial differential
operators of order higher than two as well. \par  Thus knowing
fundamental solutions of the corresponding first or second order
operators permits to write fundamental solutions of higher order
partial differential operators considered above with the help of the
iterated convolutions in a definite order prescribed by the
induction process.

\par {\bf 36. Non-linear partial differential equations.}
\par We consider the differential equation
\par $(1)$ $(\Upsilon ^m + \beta + {\hat f}(y) \Upsilon ) y =g$, \\
where $\Upsilon ^m+\beta $ is a partial differential operator as in
Formula 10$(13)$ of order $m$, $f(y)$ is a ${\cal A}_v$
differentiable function, $y=y(z)$ is an unknown function, ${\hat
f}(y)\Upsilon y := \sum_{j=0}^n[{\hat f}(y).(\partial y(z)/\partial
z_j)]\phi _j^*(z)$. Suppose that a fundamental solution ${\cal
E}_{\Upsilon ^m + \beta }$ of Equation 35$(9)$ for the operator
$(\Upsilon ^m + \beta)$ is known. Find at first a fundamental
solution $y={\cal V}$ of $(1)$ with $g=\delta $. Then
\par $(2)$ $(\Upsilon ^m + \beta){\cal V} = \delta - \mu $, \\ where
$\mu (z)= {\hat f}(y(z))\Upsilon y(z)$. The anti-derivative gives
\par $(3)$ $w(y(z)) = \mbox{}_{\Upsilon }\int ({\hat f}(y) \Upsilon y)dz =
\int f(y(x))dy(x) = \int f(y)dy =\mbox{}_{\Upsilon }\int \mu (z)dz,$
\\ then
\par $(4)$ $y=w^{-1}(\mbox{}_{\Upsilon }\int_{\gamma ^{\alpha }|_{[a,t_z]}}
 \mu (x)dx ) ,$ \\ where $w^{-1}$
denotes the inverse function. On the other hand, \par $(5)$ $y={\cal
E}_{\Upsilon ^m + \beta } * (\delta - \mu ) = {\cal I}_{\Upsilon ^m
+\beta }(\delta -\mu )$, when $\Upsilon ^m$ is either of the first
order for $m=1$ or is expressed as a composition of operators of the
first order, \par $(5.1)$ $\Upsilon ^m + \beta = (\Upsilon _1+\beta
^1)...(\Upsilon _m+\beta ^m)$ so that \par $(5.2)$ ${\cal
I}_{\Upsilon ^m +\beta } = {\cal I}_{\Upsilon _m +\beta ^m}...{\cal
I}_{\Upsilon _1 +\beta ^1}$, \\ consequently, $(4,5)$ imply the
equation:
\par $(6)$ ${\cal E}_{\Upsilon ^m + \beta } * (\delta - \mu ) =
{\cal I}_{\Upsilon ^m + \beta }(\delta -\mu )=
w^{-1}(\mbox{}_{\Upsilon }\int_{\gamma ^{\alpha }|_{[a,t_z]}}
 \mu (x)dx )$ or
\par $(7)$ $w({\cal E}_{\Upsilon ^m + \beta } * (\delta - \mu )) =
\mbox{}_{\Upsilon }\int_{\gamma ^{\alpha }|_{[a,t_z]}}
 \mu (x)dx  $. \\ We have the identity \par
$(\partial ({\cal E}* {\Psi })/\partial z_p),\phi ) = - (({\cal E}*
{\Psi }),\partial \phi /\partial z_p)  =(({\cal E}* (\partial {\Psi
}/\partial z_p)) ,\phi )$ \\ with a generalized function $\Psi $.
\par Therefore, differentiating $(6)$ by $z_0$,...,$z_n$, we infer
that: \par $(8)$ $\{ \sum_{j=0}^n [{\cal E}_{\Upsilon ^m + \beta } *
(\partial \mu (z)/\partial z_j)] \phi _j^*(z) \} + $ \par
$(n+1)^{-1} \sum_{j,k=0}^n [(dw^{-1}(\xi )/d\xi ).({\hat \theta
}(z).(\partial \nu _j(z)/\partial z_k))]\phi _k^*(z) = \Upsilon
{\cal E}_{\Upsilon ^m + \beta }$ or
\par $(9)$ $\Upsilon ({\cal I}_{\Upsilon ^m + \beta }\mu ) +
$ \par $(n+1)^{-1}\sum_{j,k=0}^n [(dw^{-1}(\xi )/d\xi ).({\hat
\theta }(z).(\partial \nu _j(z)/\partial z_k))]\phi _k^*(z) =
\Upsilon ({\cal I}_{\Upsilon ^m +\beta }\delta )$, \\ where
\par $(9.1)$ $[((dg(z)/dz).\phi _j(z))\otimes \phi
_j^*(z), \omega ^{\otimes 2}) =[\mu (z)\otimes 1,\omega ^{\otimes
2}) $ \\ for each real-valued test function $\omega $ and each $j$,
$ ~ {\hat \theta } (z)=dg(z)/dz$, $\xi = \mbox{}_{\Upsilon
}\int_{\gamma ^{\alpha }|_{[a,t_z]}} \mu (x)dx $ (see also \S \S 4,
5, 17, 23 and 26). If $\Upsilon $ and $\beta $ are independent of
$z_j$, i.e. $\phi _j(z)=0$ is zero identically on $U$, then
$\partial ({\cal I}_{\Upsilon ^m+\beta }\mu )/\partial z_j={\cal
I}_{\Upsilon ^m+\beta }(\partial \mu /\partial z_j)$ (see also Note
23.1). Otherwise the derivative $\partial ({\cal I}_{\Upsilon
^m+\beta }\mu )/\partial z_j$ is given by Formulas 4$(6)$ and
23$(8)$ and 26$(1)$.
 The function $w$ is known from
$(3)$ after the line integration by the variable $y$, so Equation
$(8)$ is linear by $(\partial \mu (z)/\partial z_j)$, $j=0,...,n$.
It can be solved as in \cite{ludeoc}. Calculating $\mu $ from $(8)$
or $(9)$  we get the fundamental solution:  \par $(10)$  ${\cal V}
={\cal E}_{\Upsilon ^m + \beta } * (\delta - \mu )$ \\ and the
(particular) solution of $(1)$ is: \par $(11)$ $y={\cal V}* g$.
\par When $[{\cal E}_{\Upsilon ^m +
\beta }, \omega )$ is real for each real-valued test function
$\omega $ or $\Upsilon f=(\partial f/\partial z_0)\phi _0(z)$ with a
real-valued function $\phi _0(z)$ and the inverse relative to the
convolution generalized function is known ${\cal E}_{\Upsilon ^m +
\beta }^{-1}$ such that
\par $(12)$ ${\cal E}_{\Upsilon ^m + \beta }^{-1}*{\cal E}_{\Upsilon
^m + \beta }=\delta $, \\ then Equation $(8)$ simplifies:
\par $(13)$ $\Upsilon \mu (z) + (n+1)^{-1}\sum_{j,k=0}^n  {\cal E}_{\Upsilon ^m + \beta }^{-1}*
[[(dw^{-1}(\xi )/d\xi ).({\hat \theta }(z).(\partial \nu
_j(z)/\partial z_k))]\phi _k^*(z)] = \sum_{j=0}^n [{\cal
E}_{\Upsilon ^m + \beta }^{-1}*(\partial {\cal E}_{\Upsilon ^m +
\beta }/\partial z_j)]\phi _j^*(z)$, consequently,
\par $(14)$ $\Upsilon \mu (z) + \nu (\mu ) = b(z) ,$ \\ where $\nu
(\mu ) := (n+1)^{-1}\sum_{j,k=0}^n{\cal E}_{\Upsilon ^m + \beta
}^{-1}* [[(dw^{-1}(\xi )/d\xi ).({\hat \theta }(z).(\partial \nu
_j(z)/\partial z_k))]\phi _k^*(z)]$ and $b(z) = {\cal E}_{\Upsilon
^m + \beta }^{-1}*(\Upsilon {\cal E}_{\Upsilon ^m + \beta })$.
\par If equation $(1)$ is solved, then it provides a solution of
more general equation:
\par $(15)$ $(\Upsilon ^m + \beta + {\hat f} ((\Upsilon )^{k-1}\xi )
(\Upsilon )^k) \xi  =g$ \\
finding $\xi $ from the equation $(\Upsilon )^{k-1}\xi  =y$, where
$(\Upsilon )^k$ denotes the $k$-th power of the operator $\Upsilon
$. \par If $\phi _j(z)=i_j\psi _j(z)$ for each $j$, then functions
$\{ \nu _j(z): ~ j \} $ can be chosen real-valued or real-valued
generalized functions on real valued test functions (see System
4$(10,11)$ and \S 23). In such case the equality \par
$\sum_{j,k=0}^n [(dw^{-1}(\xi )/d\xi ).({\hat \theta }(z).(\partial
\nu _j(z)/\partial z_k))]\phi _k^*(z) $\par $= \sum_{j,k=0}^n
[(dw^{-1}(\xi )/d\xi ).(\mu (z)(\partial \nu _j(z)/\partial
z_k)i_j)]i_k^*\psi _k(z)]$\\ is satisfied.
 For $\Upsilon f=(\partial f/\partial z_0)\phi _0(z)$ with a
real-valued function $\phi _0(z)$ these equations simplify, since
${\hat \theta }.h = \mu (z)h$ for each $h\in {\bf R}$ and $z\in U$
and $(n+1)^{-1}\sum_{j,k=0}^n{\cal E}_{\Upsilon ^m + \beta }^{-1}*
[[(dw^{-1}(\xi )/d\xi ).({\hat \theta }(z).(\partial \nu
_j(z)/\partial z_k))]\phi _k^*(z)]={\cal E}_{\Upsilon ^m + \beta
}^{-1}* [(dw^{-1}(\xi )/d\xi ).\mu (z)]$.
\par Thus the results of this paper over the Cayley-Dickson algebras
enrich the technique of integration of partial differential
equations in comparison with the complex field.
\par It is planned to present in the next paper solutions of some types
of non-linear partial differential equations with the help of
non-linear mappings and non-commutative line integration over the
Cayley-Dickson algebras.

\par Department of Applied Mathematics,
\par Moscow State Technical University MIREA,
av. Vernadsky 78,
\par Moscow, Russia
\par e-mail: sludkowski@mail.ru
\end{document}